\documentclass[11pt]{article}
\usepackage[utf8]{inputenc}
\usepackage{blindtext}
\usepackage{soul}
\usepackage{latexsym} 
\usepackage{amsmath,amssymb,amsfonts,amsthm} 
\usepackage{graphics,psfrag} 
\usepackage{fancyhdr,lastpage} 
\usepackage{color,verbatim} 
\usepackage[hidelinks]{hyperref} 
\usepackage{hyperref} 
\usepackage[table]{xcolor} 
\usepackage{metalogo}
\usepackage{etoolbox}
\usepackage{cite}
\usepackage{mathrsfs}
\usepackage[affil-it]{authblk}

\usepackage{xfrac}
\usepackage{bbm}
\patchcmd{\thebibliography}{\section*{\refname}}{}{}{}

\newtheorem{theorem}{Theorem}[section]
\newtheorem{corollary}{Corollary}[theorem]
\newtheorem{lemma}[theorem]{Lemma}

\newtheorem{definition}[theorem]{Definition}
\newtheorem{example}[theorem]{Example}
\newtheorem{remark}[theorem]{Remark}

\hypersetup{
    colorlinks=false,
    pdfborder={0 0 0},
}

\usepackage{calc}

\usepackage[a4paper,top=1in,bottom=0.8in,right=0.8in,left=0.8in]{geometry}
\let\oldbibliography\thebibliography
\renewcommand{\thebibliography}[1]{%
  \oldbibliography{#1}%
  \setlength{\itemsep}{-1.5mm}%
}

\def\R{\mathbb{R}}

\def\N{\mathbb{N}}
\def\P{\mathbb{P}}
\def\Z{\mathbb{Z}}
\def\E{\mathbb{E}}

\newcommand{\be}{\begin{equation}}
\newcommand{\ee}{\end{equation}}
\newcommand{\bea}{\begin{eqnarray}}
\newcommand{\eea}{\end{eqnarray}}
\newcommand{\beann}{\begin{eqnarray*}}
\newcommand{\eeann}{\end{eqnarray*}}
\newcommand{\benn}{\begin{equation*}}
\newcommand{\eenn}{\end{equation*}}

\newcommand{\cA}{{\mathcal A}}  
\newcommand{\cB}{{\mathcal B}}  
\newcommand{\cF}{{\mathcal F}}  
\newcommand{\cK}{{\mathcal K}}  
\newcommand{\cL}{{\mathcal L}}  
\newcommand{\cM}{{\mathcal M}}  
\newcommand{\cP}{{\mathcal P}}  
\newcommand{\cS}{{\mathcal S}}  
\newcommand{\cT}{{\mathcal T}}  
\newcommand{\cU}{{\mathcal U}}  
\newcommand{\cW}{{\mathcal W}}  

\newcommand{\ie}{i.e.\ }

\newcommand{\Space}{\mathbf{Z}}

\newcommand{\CbFunct}{C_b(\Space)}

\newcommand{\Graphon}{\mathcal{W}_1}

\newcommand{\UGraphon}{\widetilde{\mathcal{W}}_1}

\newcommand{\UGraphond}{\widetilde{\mathcal{W}}_{1,d}}

\newcommand{\InvRelabel}{S_{[0,1]}}
\newcommand{\Relabel}{\bar{S}_{[0,1]}}

\newcommand{\NcutR}[1]{\Vert#1\Vert_{\square,\R}}

\newcommand{\NcutRSymbol}{\Vert\cdot\Vert_{\square,\R}}

\newcommand{\simd}{\sim_{d}}

\newcommand{\dd}{\delta_{\square}}

\newcommand{\F}{\mathcal{F}}

\renewcommand{\P}{\mathbb{P}}

\newcommand{\un}{\mathbbm{1}}

\newcommand{\drv}{\mathrm{d}}
\newcommand{\rd}{\mathrm{d}}

\fancypagestyle{plain}{%
  \fancyhf{}%
  \fancyfoot[C]{\footnotesize Page \thepage\ of \pageref{LastPage}}%
}

\title{Probability graphons and P-variables: two equivalent viewpoints for dense weighted graph limits}
\author{Giulio Zucal\thanks{giulio.zucal@mis.mpg.de}}
\affil[1]{Max Planck Institute for Mathematics in the Sciences, Leipzig, Germany}

\date{\today}

\begin{document}
\maketitle

\begin{abstract}
We develop further the graph limit theory for dense weighted graph sequences. In particular, we consider probability graphons, which have recently appeared in graph limit theory as continuum representations of weighted graphs, and we introduce $P-$variables, which also appear in the context of the Aldous-Hoover theorem for exchangeable infinite random arrays, as an alternative continuum representation for weighted graphs. In particular, we explain how $P-$variables are related to probability graphons in a similar way in which random variables are related to probability measures. We define a metric for $P-$variables (inspired by action convergence in the graph limit theory of sparse graph sequences) and show that convergence of $P-$variables in this metric is equivalent to probability graphons convergence. We exploit this equivalence to translate several results from the theory of probability graphons to $P-$variables. In addition, we prove several properties of $P-$variables convergence, thus showing new properties also for probability graphons convergence and demonstrating the power of the connection between probability graphons and $P-$variables. Furthermore, we show how $P-$variables convergence can be easily modified and generalised to cover other combinatorial structures such as bipartite graphs and hypergraphs.

    \vspace{0.2cm}
\noindent {\bf Keywords:}  Graph limits, large networks, probability graphons, edge-decorated graphs, dense weighted graph sequences, random matrices, action convergence, hypergraphs

    \vspace{0.2cm}
\noindent {\bf  Mathematics Subject Classification Number:}  05C80 (Random Graphs)   60B20 (Random matrices) 60B10 (Convergence of measures)
\end{abstract}

\section{Introduction}
\subsection{Motivation and background}

The analysis of complex networks plays a central role in many disciplines, including urban systems \cite{barthelemy_2016_structure}, epidemiology \cite{pastor-satorras_2015_epidemic}, electrical power grids \cite{pagani_2013_power}, economics \cite{hausmann_2013_atlas, battiston_2012_debtrank} and neurobiology \cite{fornito_2016_fundamentals}. Graphs are the natural mathematical structures to represent complex networks allowing the sharing of concepts and tools across different areas. However, the graphs considered in these applications are so large that combinatorial techniques are often not feasible. Moreover, for large networks such as the brain or the Internet, only approximate information about quantities like the number of graph nodes is available. For this reason, graph limit theory emerged as an effective alternative to study large networks. Graph limit theory deals with the convergence of graph sequences and their limit objects \cite{LovaszGraphLimits}. In fact, in graph limit theory the limits are continuum objects which still encode important information forgetting superfluous information.

Graphons \cite{BORGS20081801,Lovsz2007SzemerdisLF, LOVASZ2006933,borgs2011convergentAnnals} are continuum limits of dense simple graphs sequences, i.e.\ when the simple graphs considered in the sequence have asymptotically number of edges proportional to the square of the number of vertices. On the other hand, convergence notions for graph sequences with a uniform bound on the degrees have also been thoroughly studied. In this case, the most famous convergences are Benjamini-Schramm convergence \cite{BenjaminiLimit}, and its stronger version called local-global convergence \cite{local-global1,Hatami2014LimitsOL}. We refer the reader also to the monograph \cite{LovaszGraphLimits}.

Graph limits have found applications in several domains as nonparametric statistics \cite{wolfe2013nonparametricgraphonestimation,BorgsNonparaStat}, network dynamics \cite{Kuehn_2020GraphlimitDynam1,Kuehn_2019GraphlimitDynam2,bramburger2023pattern,MedvedevGraphLimitDynamics2,MedvedevGraphLimitDynamics1,gkogkas2022mean,bick2024dynamicalGraphLimi,ayi2023graphlimitinteractingparticle,ayi2024meanfieldlimitnonexchangeablemultiagent,kuehn2022vlasov,jabin2021meanfield,luccon2020quenched,bet2024weakly} and Graph Neural Networks (GNNs) \cite{keriven2024functions,NEURIPS2023_8154c89c,maskey2022generalization,maskey2023transferability,le2023poincar,levie2024graphon}. 

Recently, graph limits theory has developed in several new directions, many motivated by applications. 

In particular, as in applications networks are sparse and heterogeneous, the development of convergence notions for graph sequences of intermediate density (that are neither dense nor of uniformly bounded degree) attracted a lot of interest \cite{BORGS20081801,borgs2011convergentAnnals,MarkovSpaces,KUNSZENTIKOVACS20191,frenkel2018convergence,backhausz2018action,veitch2015classrandomgraphsarising,janson2016graphons,caron2017sparse,borgs2018sparse,borgs2020identifiability,borgs2019sampling,10.1214/18-AOS1778,JANSON2022103549}. Among others, action convergence \cite{backhausz2018action} has the advantage of bringing graphons convergence, local-global convergence and the convergence of intermediate density sequences of simple graphs into the same mathematical framework. Related works are \cite{ArankaAction2022,MeasTheorActionZucal,zucal2023action}.

Another interesting direction motivated by applications is the extension of graph limits to more general combinatorial objects such as hypergraphs and simplicial complexes \cite{hypergrELEK20121731,HypergraphsSzegedy2,HypergraphonsZhao,zucal2023action}. Recently, higher-order interactions (interactions beyond pairwise) and the phenomena they generate attracted a lot of interest in the study of complex networks \cite{HypergraphsNetworks,HigerOrdIntBook,carletti2020dynamical,majhi2022dynamics,MHJ,HypergraphsDynamics,Bohle_2021,JOST2019870,JostMulasBook}. Hypergraphs and simplicial complexes are the natural mathematical structures for the representation of higher-order interactions.

Recently the limit theory of matroids also emerged \cite{SubModLovaszquotientconvergence,berczi2024cycle} and graph limit theory has also been employed to obtain new results about random matrix theory and a deeper connection seems to be emerging \cite{RandomMatricesGraphonsZhu,RandomMatrixGraphonmale2014,backhausz2018action}.

Related to random matrices, the development of a limit theory for dense weighted graphs (or more generally edge-decorated graphs) \cite{lovász2010limits, falgasravry2016multicolour, KUNSZENTIKOVACS2022109284,abraham2023probabilitygraphons,zucal2024probabilitygraphonsrightconvergence} and the associated limit objects called probability graphons attracted quite some interest. In fact, in applications, edges encode important additional features about how two vertices interact, such as how intense the interaction is for example. In this work, we develop this further, employing the unlabelled cut metric studied in \cite{abraham2023probabilitygraphons} and the  `right convergence' viewpoint developed in \cite{zucal2024probabilitygraphonsrightconvergence} to connect probability graphons convergence with action convergence \cite{backhausz2018action} and the Aldous-Hoover theorem for infinite exchangeable arrays\cite{aldous1981representations,aldous2010exchangeability,hoover1979relations,austin2008exchangeable,orbanz2014bayesian}.

In this work, we consider $P-$variables, mathematical objects related to probability graphons. The relationship between probability graphons and $P-$variables is similar to the relationship between probability measures and random variables. We explain this relationship in detail, explaining how to naturally associate a $P-$variable to a probability graphon and vice versa.
Moreover, we define a pseudo-metric on the space of $P-$variables and we show that this metric is topologically equivalent (on tight subsets) to the unlabelled cut metric for probability graphons introduced in \cite{abraham2023probabilitygraphons}. Exploiting this equivalence, we translate several results from the theory of probability graphons to $P-$variables and prove new results for $P-$variables and probability graphons.

We summarize and explain the results obtained in this work in more detail in the next section.

\subsection{Summary of the results}
Let $\Space$ be a measurable subset of the real line $\R$ and $[0,1]$ the unit interval endowed with the Lebesgue measure. A probability graphon is a measurable map $\widetilde{W}$ from the unit square $[0,1]\times [0,1]$ (or more generally from the product $\Omega_1\times \Omega_1$ of a probability space with itself) to the space of probability measures on $\Space$ (equipped with the topology induced by the weak convergence of measures). In \cite{abraham2023probabilitygraphons}, the convergence of probability graphons has been characterized with respect to the unlabelled cut-metrics $\delta_{\square}$ (see Definition \ref{def:ddcut} for example) and in \cite{zucal2024probabilitygraphonsrightconvergence}, the equivalent global/`right convergence' viewpoint on probability graphons has been developed. 

In this work, we also consider $P-$variables, i.e.\ measurable maps 
$$
W:\Omega_1 \times \Omega_1 \times \Omega_2 \rightarrow \Space\subset \R,
$$
where $\Omega_1\times \Omega_1\times \Omega_2$ denotes the product of the probability spaces $\Omega_1$ and $\Omega_2.$ These objects have already appeared in the celebrated Aldous-Hoover theorem for exchangeable infinite random arrays \cite{aldous1981representations,aldous2010exchangeability,hoover1979relations,austin2008exchangeable} which is known to be related to (real-valued) graphons \cite{austin2008exchangeable,diaconis2007graph,aldous2010exchangeability}.    

We describe how one can naturally identify weighted graphs with $P-$variables. In particular, one can consider a weighted graph $H$ on the vertex set $[n]$ with the adjacency matrix $A(H)=(A_{i,j})_{i,j\in [n]}.$ The set $[n]$ is a probability space if we endow it with the uniform measure. Therefore, one can consider the $P-$variable $H$ from $[n]\times [n]\times \Omega_2$ to $\R$ such that $H(i,j,x_{12})=A(H)_{i,j}$ for $i,j\in [n]$ and any $x_{12}\in \Omega_1.$ Weighted graphs can also be identified with probability graphons, see \cite{abraham2023probabilitygraphons, zucal2024probabilitygraphonsrightconvergence} for example. Therefore, any convergence notion for $P-$variables (or probability graphons) gives us naturally a convergence notion for weighted graphs.

Moreover, $P-$variables can be interpreted as random variables versions of probability graphons in the following way. For a $P-$variable $W$ as defined above and all $x_1,x_2\in \Omega_1$ we can consider the measurable function (or in probabilistic language random variable) $W(x_1,x_2,\cdot)$ from $\Omega_1$ to $\Space$ and its distribution $\widetilde W(x_1,x_2)$ that is a probability measure on $\Space$. In this way, we obtain a probability graphon $\widetilde W$ that associates to every $x_1,x_2\in \Omega_1$ a probability measure $\widetilde W(x_1,x_2)$ in a measurable way.

In Section \ref{Sec4:EquivPsetsProbGraphon}, we also explain how one can conversely associate to a probability graphon $\widetilde{W}$ a $P-$variable $W.$

We develop further the connection between probability graphons and $P-$variables. In particular, we define a pseudo-metric on the space of $P-$variables, the $P-$variables metric, using a similar idea to the construction of the action convergence metric in \cite{backhausz2018action}. Moreover, we show that (on tight subsets) $P-$variables metric for $P-$variables and the unlabelled cut metric for probability graphons are topologically equivalent, i.e.\

$$
\textit{probability graphon convergence and $P-$variables convergence are equivalent,}
$$
that is the main result of this work. We prove this result by developing an alternative formulation of $P-$variables convergence and using the global/`right convergence' viewpoint for probability graphons developed in \cite{zucal2024probabilitygraphonsrightconvergence}.

Moreover, we translate several results from the theory of probability graphons to the language of $P-$variables. In particular, we characterize $P-$variables identified by their metric and the compact sets of $P-$variables in their metric. Moreover, we prove that weighted graphs are dense in the space of $P-$variables, showing that sampling from a $P-$variable, as in Aldous-Hoover theorem  \cite{aldous1981representations,aldous2010exchangeability,hoover1979relations,austin2008exchangeable}, larger and larger weighted graphs one approximates the $P-$variable in the $P-$variables metric. We also characterize $P-$variables convergence in terms of  invariants as homomorphism densities, quotient sets and overlay functionals. Therefore, we develop the equivalent local/`left convergence' and the global/`right convergence' viewpoints for $P-$variables. See also \cite{LovaszGraphLimits,borgs2011convergentAnnals,Lovsz2007SzemerdisLF, LOVASZ2006933} and \cite{abraham2023probabilitygraphons,zucal2024probabilitygraphonsrightconvergence} for the definition of these invariants and the equivalent local/`left convergence' and the global/`right convergence' viewpoints in the (real-valued) graphons and probability graphons cases. In addition, we provide several examples of converging sequences of weighted graphs in the $P-$variables metric and their limits.

In addition, we also use the $P-$variables formulation to obtain new results for probability graphons. 
In fact, we show that many properties are preserved under $P-$variables convergence as symmetry, taking values in a closed set (being $\{0,1\}-$valued or being non-negative for example) or having a uniform bound on the $L^p-$norm for $p>1$. These results can be easily turned into results for probability graphons. The proofs of these results are very simple showing the power of $P-$variables convergence.

Moreover, using the $P-$variable formulation, we show that if a sequence of probability graphons $(\widetilde{W}_n)_n$ converges in the unlabelled cut metric to $\widetilde{W}$ and satisfies a suitable uniform bound then the sequence $(w_n)$ of contractions of $(\widetilde{W}_n)_n$ converges to $w$ the contraction of $\widetilde{W}$ in (real-valued) graphons convergence. For a probability graphon $\widetilde{W}$ with contraction we mean the function from $[0,1]\times [0,1]$ to $\R$ defined as  

\begin{equation*}
w(x_1,x_2)=\int_{\Space \subset\mathbb{R}}z\widetilde W(x_1,x_2,\mathrm{d}z).
\end{equation*}
for all $x_1,x_2\in[0,1],$ that is an averaging of the probability graphon $\widetilde W.$

To conclude we present generalizations and modifications of $P-$variables convergence for vectors of $P-$variables, bipartite graphs, graphs with self-loops, sparse graph sequences and hypergraphs. In particular, $P-$variables convergence provides a systematic and transparent way to define hypergraph limits that we will develop in future work explaining the connections with \cite{HypergraphonsZhao}, \cite{zucal2023action} and \cite{aldous1981representations,aldous2010exchangeability,hoover1979relations,austin2008exchangeable}.

\subsection{Organization}

In Section \ref{Sec1Notation}, the notation is fixed and basic notions about analysis, probability, and graph theory are introduced. In Section \ref{Sec2ProbGraphons}, we summarise the theory of probability graphons and the associated unlabelled cut-metric. In Section \ref{Sec3Psets}, we introduce $P-$variables and the associated $P-$variables metric. Moreover, we develop an equivalent formulation of the convergence in the $P-$variables metric. In Section \ref{Sec4:EquivPsetsProbGraphon}, we explain how to associate probability graphons to $P-$variables and vice versa. Moreover, we prove the main result of this work, the equivalence of unlabelled cut-metric convergence for probability graphons and $P-$variables convergence. In Section \ref{Sec5TranslatResultsExamp}, several results from the theory of probability graphons are translated into the language of $P-$variables and examples of sequences of weighted graphs and their limits are given. In Section \ref{Sec6ImpliesRealValuedGraphons}, we show that, if a sequence of $P-$variables has $L^p-$norm for $p>1$ uniformly bounded then $P-$variables convergence (and therefore also probability graphons convergence) implies the convergence of the contraction of the $P-$variable (or probability graphon) considered as a real-valued graphon. In Section \ref{Sec7OtherProp}, it is shown that many properties of $P-$variables (and probability graphons) are preserved under $P-$variables convergence. Finally, in Section \ref{Sec8GeneralizationsAndDirect} we present possible generalizations and modifications of the theory of $P-$variables for covering bipartite graphs, graphs with self-loops, sparse graphs and hypergraph sequences.

\section{Notation}\label{Sec1Notation}
This section introduces the notation and recalls basic notions from analysis, probability and graph theory.

We will denote with $(\Omega,\mathcal{A},\P)$ a probability space where $\mathcal{A}$ is a  $\sigma-$algebra and $\P$ is a probability measure on $(\Omega,\mathcal{A})$. When there is no risk of confusion, we will often denote a probability space $(\Omega,\mathcal{A},\P)$ only with $\Omega,$ omitting the probability measure and the $\sigma-$algebra. For two probability spaces $(\Omega_i, \mathcal{A}_i, \P_i)$, $i=1,2$ we will denote with $(\Omega_1\times \Omega_2,\mathcal{A}_1\times \mathcal{A}_2,\P_1\otimes\P_2)$ their product probability space. When there is no risk of confusion we will write $\Omega_1\times \Omega_2.$ We will denote with $\mathcal{P}(\Space,\mathcal{F})$ or shortened $\mathcal{P}(\Space)$ the set of probability measures on the measurable space $(\Space,\mathcal{F})$. In addition, we indicate the expectation of a real-valued measurable function (or in probabilistic language a random variable) $f$ on the probability space $(\Omega,\mathcal{F},\P)$ with $\E[f]$. We indicate the (possibly infinite) $L^p-$norm of a real-valued measurable function $f$ with  
$$
\|f\|_p=\left(\int_{\Omega}|f(\omega)|^pd\P(\omega)\right)^{\frac{1}{p}}=\left(\E[|f|^p]\right)^{\frac{1}{p}}
$$
if $1\leq p \leq +\infty$ and 
$$\|f\|_{\infty}=\sup_{\omega\in \Omega}|f(\omega)|,$$
where with $\sup$ here we mean the essential supremum with respect to $\P$.
If a measurable function $f$ has finite $L^p-$norm we say that $f$ is $p-$integrable (or has finite $p-$moment). We denote with $L^p(\Omega,\mathcal{F},\P)$ the usual Banach space of the real-valued measurable $p-$integrable functions (identified if they are equal almost everywhere) on $(\Omega,\mathcal{F},\P)$ equipped with the $L^p-$norm or equivalently, in probabilistic language, the space of random variables with finite $p-$moment. We will often use the shortened notations $L^p(\Omega)$ or $L^p$ when there is no risk of confusion. For a measurable set $S\subset \R$ we will also denote with $L_S^p(\Omega)$ the space of the $p-$integrable random variables taking values in $S$.

A $k-$dimensional random vector is a measurable function $\mathbf{f}$ from a probability space $(\Omega,\mathcal{F},\P)$ to $\R^k$ and we can naturally represent it as $$
\mathbf{f}=(f_1,\ldots,f_k)
,$$ where $f_1,\ldots,f_k$ are real-valued random variables on $(\Omega,\mathcal{F},\P)$. Thus, a real-valued random variable is a $1-$dimensional random vector. For a $k-$dimensional random vector $\mathbf{f}$, we denote with $\mathcal{L}(\mathbf{f})=\mathcal{L}(f_1,\ldots,f_k)$ its distribution (or law), i.e.\ the measure on $\R^k$ defined for every set $A$ in the Borel sigma-algebra of $\R^k$ (a measurable set of $\R^k$) as 
$$
\mathcal{L}(\mathbf{f})(A)=\P(\mathbf{f}^{-1}(A))
$$
 where $f^{-1}$ denotes the preimage of the function $f.$ 
\newline
Let $[0,1]$ be the unit interval. We will consider $[0,1]$ as the normed probability space endowed with the Euclidean norm, the Borel $\sigma-$algebra and the Lebesgue measure, which we denote with $\lambda,$ if not specified differently.

Given $n\in\N$, we denote by $[n]$ the set $\{1,\ldots,n\}$. The set $[n]$ can be endowed with the uniform measure to obtain a probability space. For finite probability spaces, the law of a random vector is easy to represent. We illustrate this with the following example.
\begin{example}
Let's consider the probability space $([n],\mathcal{U},\mathcal{D})$ where $\mathcal{U}$ is the uniform probability measure on $[n]$ and $\mathcal{D}$ is the discrete $\sigma$-algebra on $[n]$. Then for any $k-$dimensional random vector $$
\mathbf{f}=(f_1,\ldots,f_k)$$

the law $\mathcal{L}(\mathbf{f})$ is 
$$
\mathcal{L}(\mathbf{f})=\frac{1}{n}\sum^n_{i=1}\delta_{(f_1(i),\ldots,f_k(i))}
$$
where $\delta_{(x_1,\ldots,x_k)}$ is the Dirac measure centered in $(x_1,\ldots,x_k)\in \R^k$. \newline
\end{example}
Let $f : \Omega_1 \to \Omega_2$ be a measurable map between a probability space $(\Omega_1, \mathcal{A}_1, \P_1)$ and a measurable space $(\Omega_2,\cA_2).$ The pushforward measure of $\P_1$ through $f$ is the probability measure $f_{\#}\P_1$ defined for any $A\in \cA_2$ as $$
f_{\#}\P_1(A)=\P_1(f^{-1}(A)).$$

For a measurable function $g:\Omega_2  \rightarrow \R$ and a measurable set $X_2\subset \Omega_2$ we have the equality 
\begin{equation}\label{EqPushForwardInt}
\int_{X_1}g(f(\omega_1))\P_1(\mathrm{d}\omega_1)=\int_{X_2}g(\omega_2)f_{\#}\P_1(\mathrm{d}\omega_2)
\end{equation}
where $X_1=f^{-1}(X_2).$

Observe that for a random variable $f:(\Omega,\cA,\P) \rightarrow \R$ (or more generally a random vector) we have that its distribution $\mathcal{L}(f)$ is the pushforward of $\P$ through $f,$ i.e.\ $\mathcal{L}(f)=f_{\#}\P.$ 

A map $\varphi : \Omega_1 \to \Omega_2$ between two probability spaces $(\Omega_i, \mathcal{A}_i, \P_i)$, $i=1,2$, is measure-preserving if it is measurable and if for every $A\in \mathcal{A}_2$, $\P_2(A) = \P_1(\varphi^{-1}(A)),$ i.e.\ if $\varphi_{\#}\P_1=\P_2$.
In this case, for every measurable non-negative function $f: \Omega_2 \to \R$,
we have:
\begin{equation}
   \label{eq:re-label}
 \int_{\Omega_1} f(\varphi(x))\ \drv \P_1(x) = \int_{\Omega_2} f(x)\
 \drv \P_2(x).
\end{equation}
We denote by $  \InvRelabel$ the set of bijective measure-preserving maps from $[0,1]$ with the Lebesgue measure to itself, and by $\Relabel$ the set of   measure-preserving maps from $[0,1]$ with the Lebesgue measure to itself. 
\newline
Let $\Omega_1$ be a probability space. For a measurable subset $S\subset \Omega_1$ we denote with $\mathbbm{1}_{S}$ the indicator function of $S,$ i.e.\ the function taking value $1$ if $x\in S $ and $0$ if $x\in \Omega_1\setminus S. $
Let $\Omega_2$ be another probability space. For two measurable functions $f:\Omega_1\rightarrow \R$ and $g:\Omega_2 \rightarrow \R$ we denote with $f\otimes g$ the measurable function defined as $$
f\otimes g:\Omega_1 \times \Omega_2\rightarrow \R$$
$$(x_1,x_2)\mapsto (f\otimes g)(x_1,x_2)=f(x_1)g(x_2).$$
In particular, we have that $(f\otimes\mathbbm{1}_{\Omega_2})(x_1,x_2)=f(x_1)$ and $(\mathbbm{1}_{\Omega_1}\otimes g)(x_1,x_2)=g(x_2).$ Observe that $\mathcal{L}(f)=\mathcal{L}(f\otimes\mathbbm{1}_{\Omega_2})$ and $\cL(g)=\cL(\mathbbm{1}_{\Omega_1}\otimes g)$ and therefore we have also $\|f\|_p=\|f\otimes\mathbbm{1}_{\Omega_2}\|_p$ and $\|g\|_p=\|\mathbbm{1}_{\Omega_1}\otimes g\|_p.$ We will often implicitly use this fact throughout this work.
\newline

For a measure $\mu$ on a metric space $X,$ considered with its Borel $\sigma-$algebra (in this work we will typically consider some measurable subset of $\R^m$), and a real-valued measurable function $f$ defined on $X$, we denote by $\mu[f]=\mu(f)=\langle \mu, f \rangle=
\int f \, \rd \mu=\int_{X} f(x)\, \rd\mu$ the integral of $f$ with respect to $\mu$ when well defined. 

We say that a sequence  of signed measures $(\mu_n)_{n\in\N}$  weakly converges to some $\mu$ if for each function $f\in   \CbFunct$, we have $\lim_{n\to +\infty} \mu_n(f) = \mu(f)$. \newline

We recall here the Lévy-Prokhorov metric on the space of probability measures.
\begin{definition}[Lévy-Prokhorov metric]\label{LevyProk}
Let $(X,d)$ be a metric space. The \emph{Lévy-Prokhorov Metric} $d_{\mathcal{LP}}$ on the space of probability measures (or more generally on the space of measures) on $(X,d)$ is for $\eta_1,\eta_2$ measures on $X$
$$\begin{aligned}
d_{\mathcal{LP}}\left(\eta_{1}, \eta_{2}\right)=&\inf \left\{\varepsilon>0: \eta_{1}(U) \leq \eta_{2}\left(U^{\varepsilon}\right)+\varepsilon \text{ and } \right.\\
&\left.\eta_{2}(U) \leq \eta_{1}\left(U^{\varepsilon}\right)+\varepsilon  \text{ for all } U \in \mathcal{B}\right\},
\end{aligned}$$
where $\mathcal{B}$ is the Borel $\sigma$-algebra on $X$, $U^{\varepsilon}$ is the set of points that have distance smaller than $\varepsilon$ from $U$.
\end{definition}
\begin{remark}
It is well known that the Levy-Prokhorov metric metrizes weak convergence on the space of measures when the metric space is separable, see Theorem 6.8 in \cite{billingsley1968convergence} for example.
\end{remark}

\begin{remark}
In this work, we will mainly focus on the Levy-Prokhorov metric when the metric space $(X,d)$ is $\R^m$ with the Euclidean distance (that is separable).
\end{remark}

We observe that the Lévy-Prokhorov metric is upper-bounded by $1$ (when only probability measures are compared).\newline

A (simple) graph $G$ is a couple $G=(V,E),$ where $V=V(G)=\{v_1,\ldots,v_n\}$ is the vertex set and $E=E(G)$ is the edge set, edges are pairs of vertices. The adjacency matrix $A=A(G)$ of the simple graph $G$ is the symmetric matrix with entries
$$
A_{i,j}=A_{j,i}=\begin{cases}
    1 & \text{ if }\{v_i,v_j\}\in E\\
    0 & \text{ else.}
\end{cases}
$$
A weighted (possibly directed) graph $H$ is a triple $H=(V,E,\beta),$ where $V=V(H)=\{v_1,\ldots,v_n\}$ is the vertex set and $E=E(H)$ is the edge set, edges are ordered pairs of vertices. Moreover, $\beta=\beta(H)$ is a function
$$
\beta:V\times V\rightarrow \R
$$
$$
(v_i,v_j)\mapsto \beta(v_i,v_j)=\beta_{i,j}$$
where $\beta_{i,j}\neq 0$ if and only if $(v_i,v_j)\in E.$ We call the matrix $A=A(H)$ with entries $A_{i,j}=\beta_{i,j}$ the adjacency matrix of $H.$ This definition is compatible with the adjacency of a simple graph. In fact, the simple graph $G=(V,E)$ can be canonically represented by the weighted graph $H=(V,E,\beta)$ where $\beta_{i,j}=A_{i,j}.$ 

We will also consider measure edge-decorated graphs. A measure edge-decorated graph is a a triple $H^{\tilde \beta}=(V,E,\tilde \beta),$ where $V=V(H)=\{v_1,\ldots,v_n\}$ is the vertex set and $E=E(H)$ is the edge set, edges are ordered pairs of vertices. Moreover, $\tilde \beta=\tilde \beta(H)$ is a function
$$
\tilde\beta:V\times V\rightarrow \cP(\R)
$$
$$
(v_i,v_j)\mapsto \tilde\beta(v_i,v_j)=\tilde \beta_{i,j}$$
where $\tilde \beta_{i,j}\neq 0$ if and only if $(v_i,v_j)\in E.$ We will also consider vertex-weighted measure edge-decorated graphs. A vertex-weighted measure edge-decorated graph $H_{\alpha}^{\tilde{\beta}}$ is a measure edge-decorated graphs $H^{\tilde \beta}=(V,E,\tilde \beta)$ equipped with a weight function $\alpha$ on the vertex set $V,$ i.e.\ $\alpha:V\rightarrow \R.$ We consider now the embedding $$
in:\R \rightarrow \cP(\R)
$$
$$s \mapsto in(s)=\delta_s,$$
where $\delta_s$ denotes the Dirac measure centred at $s\in \R.$ We can naturally use the embedding $in$ to turn a weighted graph $H=(V,E,\beta)$  into a measure-decorated graph $H^{\widetilde \beta}=(V,E,\widetilde\beta)$ where $\widetilde \beta_{i,j}=in(\beta_{i,j})=\delta_{\beta_{i,j}}.$

A bipartite graph (or bigraph) is a triple $B=(U,V,E)$ where $U=\{u_1,\ldots,u_n\}$ and $V=\{v_1,\ldots,v_m\}$ are the two disjoint vertex sets and the edge set $E$ is a set of pairs containing one element of $U$ and one element of $V.$ The matrix $M\in \R^{n\times m}$ associated to the bipartite graph $B$ is the matrix with entries for $i\in [n]$ and $j\in [m]$$$
M_{i,j}=\begin{cases}
    1 & \text{ if }\{u_i,v_j\}\in E\\
    0 & \text{ else.}
\end{cases}
$$
\newline
We will call a measurable function from $\Omega_1\times \Omega_1$ to $[0,1],$ where $\Omega_1$ is a probability space, a real-valued graphon. We will typically consider $\Omega_1=[0,1]$ endowed with the Lebesge measure.

\section{Probability graphons}\label{Sec2ProbGraphons}

\subsection{Real-valued graphons}
We briefly recall here the notion of (real-valued) graphons from dense graph limit theory \cite{LovaszGraphLimits,LOVASZ2006933,BORGS20081801,borgs2011convergentAnnals}.
\begin{definition}[Real-valued graphon]
   Let $\Omega_1$ be a probability space. A \emph{real-valued graphon} is a measurable function
    $$w:\Omega_1\times \Omega_1\rightarrow [0,1].$$
\end{definition}

\begin{remark}
Sometimes, one requires real-valued graphons to be symmetric in the two variables or $\Omega_1=[0,1]$. We do not impose this condition a priori here. 
\end{remark}

Let's suppose now that $\Omega_1=[0,1].$

We denote by $\NcutRSymbol$ the real-valued cut norm defined for (linear combinations of) real-valued graphons as:
\begin{equation}
	\label{eq_def_NcutR}
\NcutR{w} = \sup_{S,T \subset [0,1]} \left\vert \int_{S\times T} w(x,y)\ \drv x\drv y \right\vert
\end{equation}
where $w$ is a (linear combination) of real-valued graphons. 

Let $u,w$ two real-valued graphons one can similarly define the real-valued cut distance 
$$\delta_{\square,\R}(u,w)=\inf_{\varphi\in \InvRelabel}\NcutR{u-w^{\varphi}}.$$
See \cite{LovaszGraphLimits} for more details about the real-valued graphons, the real-valued cut-norm and the real-valued cut metric.

\subsection{Definition of probability graphons}\label{section_def_graphons}

In this section, we introduce the theory of probability graphons developed in \cite{abraham2023probabilitygraphons}. Other references are \cite{zucal2024probabilitygraphonsrightconvergence,KUNSZENTIKOVACS2022109284,lovász2010limits}. Probability graphons are a generalization of real-valued graphons \cite{LovaszGraphLimits,LOVASZ2006933,BORGS20081801,borgs2011convergentAnnals}.

Recall that for a measurable set $\Space\subset \R$ we denote with $\mathcal{P}(\Space)$ the space of probability measures on $\Space$.

\begin{definition}[Probability graphon]\label{DefProbGraphon}
Let $\Space\subset \R$ be a Borel measurable subset of $\R$ and $\Omega_1$ be a probability space. A \emph{probability graphon} on $\Space$ (or $\mathcal{P}(\Space)$-valued kernel) is a map $\widetilde W$ from $\Omega_1\times \Omega_1$ to $\mathcal{P}(\Space)$, such that:
\begin{itemize}
    \item $\widetilde W$ is a probability measure in $\mathrm{d} z$: for every $(x_1, x_2) \in \Omega_1^2, \ W(x_1, x_2 ; \cdot)$ belongs to $\mathcal{P}(\Space)$.
    \item $\widetilde W$ is measurable in $(x_1, x_2)$: for every measurable set $E \subset \Space$, the function $(x_1, x_2) \mapsto$ $\widetilde W(x_1, x_2 ; E)$ defined on $\Omega_1\times \Omega_1$ is measurable.
\end{itemize} 
\end{definition}
We indicate with  $\Graphon$ the space of probability graphons on $\Omega_1=[0,1],$ i.e.\ probability graphons from $[0,1]\times [0,1]$ to $\mathcal{P}(\Space)$, where we identify probability graphons that are  equal almost everywhere on $[0,1]^2$ (with respect to the Lebesgue measure).  

\begin{remark}[Probability graphons $\widetilde W : \Omega_1\times \Omega_1 \to \cP(\Space)$]
One usually considers probability graphons where $\Omega_1=[0,1]$ is the unit interval endowed with the Lebesgue measure (see Definition 3.1 in \cite{abraham2023probabilitygraphons} or Definition 3.3 in \cite{zucal2024probabilitygraphonsrightconvergence}). We will also focus on this particular case in the following, i.e.\ on the space $\Graphon.$ However, the theory can be developed for the case of a general probability space $\Omega_1,$ see Remark 3.4 in \cite{abraham2023probabilitygraphons} and, for the case of real-valued graphons, \cite{jansonGraphonsCutNorm2013}.  
\end{remark}

\begin{remark}\label{RemarkgeneralCasePolish}
The definition of probability graphons used here is a special case of the probability graphons considered in \cite{abraham2023probabilitygraphons} and \cite{zucal2024probabilitygraphonsrightconvergence}. In particular, we deal only with probability measures on $\R$ (or on a measurable subset $\Space$ of $\R$). In Definition 3.1 in \cite{abraham2023probabilitygraphons} the more general situation of probability measures on a Polish space is considered. We expect all of our results to hold also for the more general probability graphons defined on Polish spaces but we focus on the case of probability measures on $\R$ to avoid technical complications for readability.
\end{remark}

\begin{remark}\label{RmkProbBigraphons}
One can also consider another important generalization (that has not appeared in the literature until now to the best of our knowledge) that we will consider only tangentially in this work. A \emph{probability bigraphon}  is a map $\widetilde W$ from $\Omega_1\times \Omega_2$ to $\mathcal{P}(\R),$ such that:
\begin{itemize}
    \item $\widetilde W$ is a probability measure in $\mathrm{d} z$: for every $(x_1, x_2) \in \Omega_1\times \Omega_2, \ \widetilde W(x_1, x_2 ; \cdot)$ belongs to $\mathcal{P}(\R)$.
    \item $\widetilde W$ is measurable in $(x_1, x_2)$: for every measurable set $E \subset \R$, the function $(x_1, x_2) \mapsto$ $\widetilde W(x_1, x_2 ; E)$ defined on $\Omega_1\times \Omega_2$ is measurable.
    
    \end{itemize}
In particular, for probability bigraphons the product probability measure $\Omega_1 \times \Omega_2$
is obtained by (possibly) different probability spaces.    
One can link probability bigraphons to (weighted random) bipartite graphs similarly to how probability graphons are related to (weighted random) graphs (see \cite{abraham2023probabilitygraphons}).  
\end{remark}

For  $\cM\subset \cP(\Space)$, we  denote  by $\cW  _\cM$  the subset  of probability graphons  $W\in  \Graphon$  which are  $\cM$-valued: $W(x,y; \cdot)\in \cM$ for every $(x,y)\in [0, 1]^2$.

\begin{example}[On real-valued kernels]\label{rem:real-valued-kernels}
Every real-valued  graphon  $w$ can be represented as a probability graphon $W$ in the following way. Let's consider $\Space =\{ 0, 1\}$ with the discrete metric and the probability graphon $W$ defined as $W(x,y;\drv z)  = w(x,y)  \delta_1(\drv z) +  (1-w(x,y)) \delta_0(\drv z)$ for every $x,y\in[0,1]$, where we recall that $\delta_s$ is the Dirac mass located at $s\in \R$. In particular, we have $$w(x,y)=W(x,y; \{1\})=\int_{\Space}W(x,y,\drv z)$$ for $x,y \in [0, 1]$. 
\end{example}

We say that a measure-valued kernel or graphon $W$ is symmetric 
if for almost every $x,y \in [0,1]$, $W(x,y;\cdot)=W(y,x;\cdot)$.

\begin{remark}[Symmetric kernels]
We will consider in general non-symmetric probability graphons to handle also directed graphs whose
adjacency matrices are thus \emph{a priori} non-symmetric.
\end{remark}
For a probability graphon $W$ and a continuous and bounded function $f\in \CbFunct$ we denote by $W[f]$ the real-valued graphon defined by 
$$
W[f](x,y)=W(x,y;f)=\int_{\Space}f(z)W(x,y,\mathrm{d}z)
$$
for each $x,y\in [0,1].$

\subsection{The cut distance}\label{section_def_dcut}
We introduce the \emph{cut  distance}, a distance on probability graphons analogous to the  cut norm for real-valued graphons, see  \cite[Chapter 8]{LovaszGraphLimits}.
For a probability graphon $W\in\Graphon$ and a measurable subset $A\subset [0,1]^2$,
we define $W(A;\cdot),$ the measure on $\Space\subset\R$
\[
  W(A;\cdot) = \int_{A} W(x,y;\cdot)\ \drv x \drv y.
\]
We have the following semidistance on the space of probability graphons.

\begin{definition}[The cut semi-distance $d_{\square}$]
Let $d_{\mathcal{LP}}$ be the Levy-Pokhorov metric on $\mathcal{P}(\mathbf{Z}).$  
The associated cut semi-distance $d_{\square}$ is the function defined on $\Graphon^2$ by:
\begin{equation}\label{cutsemi-disteq}
    d_{\square}(U, W)=\sup _{S, T \subset[0,1]} d_{\mathcal{LP}}(U(S \times T ; \cdot), W(S \times T ; \cdot)),
\end{equation}
where the supremum is taken over all measurable subsets $S$ and $T$ of $[0,1]$.
\end{definition}

\subsection{The unlabeled  cut distance}
	\label{subsection_unlabeled_cut_distance}
Recall $\Relabel$ denotes the set of measure-preserving maps from $[0, 1]$ to $[0, 1]$ equipped with the Lebesgue measure, and $\InvRelabel$ denotes the set of bijective measure-preserving maps from $[0, 1]$ to $[0, 1]$.

We denote by $W^\varphi$ the relabeling of a probability graphon $W$ by a measure-preserving map $\varphi\in \Relabel$, i.e.\
the probability graphon defined for every $x,y\in [0,1]$ and every measurable set $A\subset \Space$ by:
\[
  W^\varphi(x,y;A) = W(\varphi(x),\varphi(y);A)
  \quad \text{for $x,y\in [0,1]$ and $ A\subset \Space$ measurable}.
\]
We can now define the cut distance for unlabeled graphons.

\begin{remark}
Instead we denote with $W^{\varphi,\psi}$ the relabelling of a probability bigraphon $W$ by measure-preserving maps $\varphi,\psi \in \Relabel,$ i.e.\ the probability bigraphon defined for every $x,y\in [0,1]$ and every measurable set $A\subset \Space$ by:
\[
  W^{\varphi,\psi}(x,y;A) = W(\varphi(x),\psi(y);A)
  \quad \text{for $x,y\in [0,1]$ and $ A\subset \Space$ measurable}.
\]

\end{remark}
\begin{definition}[The unlabeled cut distance $\delta_{\square}$]\label{def:ddcut}
 The \emph{unlabelled cut distance} is the premetric $\dd$ on $\Graphon$ such that for two probability graphons $U$ and $W$:
\begin{equation}\label{def_ddcut}
\dd(U,W) = \inf_{\varphi\in\InvRelabel} d_{\square}(U,W^\varphi)	
= \inf_{\varphi\in\InvRelabel} d_{\square}\left(U^{\varphi},W\right)		 .\end{equation}
\end{definition}

The unlabelled cut distance  $\dd$ is symmetric and satisfies the triangular inequality.  Therefore, $\dd$ defines a distance (that will still be denoted by $\dd $) on the quotient space $\UGraphond = \UGraphon / \simd$ of the space of probability graphons with respect to the equivalence relation $\simd$ defined by $U\simd W$  if and only if  $\dd(U,W)=0$.

\begin{remark}\label{RKBigraphonsDistance}
We can similarly define an unlabelled cut distance $\delta_{\square,b}$ for probability bigraphons. For two probability bigraphons $U$ and $W$ 
$$
\delta_{\square,b}(U,W)= \inf_{\varphi,\psi\in\InvRelabel} d_{\square}(U,W^{\varphi,\psi})	
= \inf_{\varphi,\psi\in\InvRelabel} d_{\square}\left(U^{\varphi,\psi},W\right).	
$$
\end{remark}

\subsection{Weak isomorphism}
\label{section_weak_isomorphism}

We introduce an equivalence relation between probability graphons called weak isomorphism. Weak isomorphism is the continuum version of the relationship between adjacency matrices of isomorphic graphs. 

\begin{definition}[Weak isomorphism]\label{def_weak_isomorphism}
We say that two probability graphons $U$ and $W$ are \emph{weakly isomorphic}
(and we note $U\sim W$) if there exists two measure-preserving maps $\varphi, \psi\in\Relabel$
such that $U^\varphi(x,y; \cdot) = W^\psi(x,y;\cdot)$ for almost every $x,y\in [0,1].$

We denote by $\UGraphon=\Graphon / \sim$ the space of unlabeled probability graphons, i.e.\
the space of probability graphons where we identify probability graphons that are weakly
isomorphic. 
\end{definition}

\begin{remark}
Similarly one can define the notion of weak isomorphism for probability bigraphons.     
\end{remark}

Some additional results on the weak isomorphism of probability graphons are summarised below. 

\begin{theorem}[Lemma 3.17 in \cite{abraham2023probabilitygraphons}]
   \label{theo:Wm=W}
   Two probability graphons are weakly isomorphic, \ie $U \sim W$, if and only if
   $U \simd W$, \emph{\ie}  $\dd(U,W) = 0$.

   Furthermore, the map $\dd$ is a distance on $\UGraphon=\UGraphond.$
\end{theorem}

The unlabelled cut metric can be defined in several different ways as the following lemma states. The following lemma is a special case of Proposition 3.18 in \cite{abraham2023probabilitygraphons}.

\begin{lemma}[Proposition 3.18 in \cite{abraham2023probabilitygraphons}]\label{thm_min_dist}
 For the unlabelled cut distance $\dd$ on $\Graphon$ we have the following equality:
\begin{equation}
  \label{eq_premetric}
  \begin{aligned}
\dd(U,W) 
& = \underset{\varphi\in \InvRelabel}{\inf} d_{\square}(U,W^\varphi) 
  = \underset{\varphi\in \Relabel}{\inf} d_{\square}(U,W^\varphi)\\
& = \underset{\psi\in \InvRelabel}{\inf} d_{\square}(U^\psi,W) 
= \underset{\psi\in \Relabel}{\inf} d_{\square}(U^\psi,W) 	\\
& = \underset{\varphi, \psi\in \InvRelabel}{\inf} d_{\square}(U^\psi,W^\varphi) 
= \underset{\varphi,\psi\in \Relabel}{\min} d_{\square}(U^\psi,W^\varphi).
\end{aligned}
\end{equation}
\end{lemma}
\subsection{Regularity lemma}

 An important special case of probability graphons are step-functions probability graphons which are  often used for
approximation.

\begin{definition}[Step-functions probability graphons]
  \label{def:stepfunction}
A probability graphon $W\in \Graphon$ is a \emph{step-function} if there exists a finite partition of $[0,1]$ into measurable (possibly empty) sets, say $\mathcal{P}=\{P_1,\cdots,P_m\}$, such that $W$ is constant on the sets $P_i \times  P_j$, for $1\leq i,j\leq  k$.
  \end{definition}

 We will be particularly interested in the following step-function probability graphons obtained from a probability graphon and a given partition of the unit interval. 
 
\begin{definition}[The stepping operator]
	\label{def_stepping_operator}
Let $W\in \Graphon$ be a probability graphon and  $\mathcal{P}=\{S_1,\cdots,S_k\}$ be a finite partition of $[0,1]$.
We define the probability graphon $W_\mathcal{P}$ adapted to the
partition  $\mathcal{P}$ by averaging $W$ over the partition subsets:
\[
  W_\mathcal{P}(x,y;\cdot) = \frac{1}{\lambda(S_i)\lambda(S_j)}\,
  W(S_i \times S_j;\cdot)
  \qquad \text{for $ x\in S_i, y\in S_j$,}
\]
when $\lambda(S_i)\neq 0$ and $\lambda(S_j)\neq 0$, and $W_\mathcal{P}(x,y;\cdot) = 0$ the null measure otherwise.
\end{definition}

\begin{remark}
The value of $ W_\mathcal{P}(x,y;\cdot)$ for  $ x\in
S_i, y\in S_j$ when $\lambda(S_i)\lambda(S_j)=0$ is not relevant as probability graphons are defined up to an almost everywhere equivalence. 
\end{remark}

The following lemma gives us some regularity of the cut metric $d_{\square}$ with respect to the stepping operator.

For a finite partition $\cP$, we will denote by $\vert\cP\vert$ the number of elements of the partition $\cP.$

\begin{lemma}
[Weak regularity lemma, Corollary 4.14 in \cite{abraham2023probabilitygraphons}]
  \label{LemmaWeakRegularity}
The distance $d_{\square}$ is weakly regular, i.e.\
 whenever $\cK$ of $ \Graphon$ is  tight,  then  for  every $\varepsilon   >  0$, there exists a positive integer $m$, such that for every probability graphon $W\in\cK$, and for every  finite partition $\mathcal{Q}$  of $[0,1]$,  there exists a  finite partition $\mathcal{P}$  of $[0,1]$  that  refines $\mathcal{Q}$ such that:
  \[
    \vert\mathcal{P}\vert  \leq m  \vert\mathcal{Q}\vert
    \quad\text{and}\quad
    d_{\square}(W, W_{\mathcal{P}}) < \varepsilon.
  \]

\end{lemma}

The Lemma \ref{LemmaWeakRegularity} is a generalization of the weak regularity lemma for real-valued graphons and simple graphs (see for example \cite[Lemma~9.15]{LovaszGraphLimits} and also \cite{Lovsz2007SzemerdisLF,fox2017regularity}).

\section{P-variables and convergence notion}\label{Sec3Psets}
\subsection{P-variables}

In the previous section, we presented probability graphons that can be interpreted as families of probability measures parametrized by the unit square. For a probability measure one can always consider a random variable distributed according to the probability measure. For this reason, one could think convenient to consider families of random variables parametrized by the unit square that enjoy a similar relation with probability graphons. This point of view will turn out quite fruitful in this paper. However, one needs to make sure that all the random variables parametrized on the unit square can be glued together in a measurable way. For this reason, we will be interested in the following objects.         

\begin{definition}
    A $P-$\emph{variable} $W$ is a measurable function $$W:\Omega_1\times \Omega_1 \times \Omega_2\rightarrow \R$$
    $$(x_1,x_2,x_{12})\mapsto W(x_1,x_2,x_{12})$$
    where $\Omega_1$ and $\Omega_2$ are probability spaces and $\Omega_1\times \Omega_1 \times \Omega_2$ is the product probability space.
    We will denote with $\mathcal{B}(\Omega_1\times \Omega_1 \times \Omega_2)$ the set of $P-$variables on $\Omega_1\times \Omega_1 \times \Omega_2,$ that is the set of real-valued measurable functions on $\Omega_1\times \Omega_1 \times \Omega_2.$ 
    \end{definition}
\begin{remark}
In probabilistic language, a $P-$variable is a random variable on the product probability space $\Omega_1\times \Omega_1 \times \Omega_2.$   
\end{remark}

\begin{remark}
One could more generally consider Polish space-valued $P-$variables, i.e.\ measurable real-valued functions from $\Omega_1\times \Omega_1 \times \Omega_2$ to $\Space,$ a Polish space. Many of the results in this work do extend to this more general case but we focus on the real-valued case to keep the technicalities to a minimum for readability. 
\end{remark}
We will come back to the relationship between probability graphons and $P-$variables in Section \ref{Sec4:EquivPsetsProbGraphon}. We explain in the following examples how one can naturally associate $P-$variables with matrices, random matrices, graphs and random graphs.

\begin{example}\label{Examp:MatrixasPset}
A square matrix $M\in \R^{n\times n}$ can be naturally interpreted as a measurable function $$
\Hat{M}=M\otimes \mathbf{1}_{\Omega}:[n]\times [n]\times \Omega\rightarrow \R$$
$$
(i,j,\omega)\mapsto (\Hat{M})(i,j,\omega)=M_{i,j},
$$
where on $[n]$ we consider the uniform measure $\cU$ (the unique measure determined by $\cU(\{i\})=1/n$ for every $i\in [n]$) and $\Omega$ is a probability space (for example $[0,1]$ with the Lebesgue measure or $[m]$ with the uniform measure). Therefore, we can naturally interpret a square matrix $M$ as a $P-$variable.
\end{example}

We can also associate a $P-$variable with a graph using the adjacency matrix and Example \ref{Examp:MatrixasPset}, see Example \ref{ExGraphasPsets}.

\begin{example}\label{Examp:RandomMatrixPset}
A random matrix $R$ is a measurable function from a probability space $\Omega$ to $\R^{n\times n}.$ We can naturally represent the random matrix $R,$ if its entries are independently distributed (we will explain why this condition is needed later), as a measurable function
$$
\Hat{R}:[n]\times [n]\times \Omega\rightarrow \R
$$
  $$
(i,j,\omega)\mapsto (\Hat{R})(i,j,\omega)=R_{i,j}(\omega),
$$
where on $[n]$ we consider the uniform measure.
\end{example}

We will discuss tangentially in the following another structure similar to $P-$variables related to probability bigraphons. 

\begin{definition}\label{DefBBBiset}
A $P-$bivariable is a measurable function $$W:\Omega_1\times \Omega_2 \times \Omega_3\rightarrow \R$$
$$(x_1,x_2,x_{12})\mapsto W(x_1,x_2,x_{12})$$
where $\Omega_1,\Omega_2$ and $\Omega_3$ are probability spaces and $\Omega_1\times \Omega_2 \times \Omega_3$ is the product probability space.
\end{definition}

\begin{example}\label{Examp:RectangMatrixPsets}
A matrix $M\in \R^{n\times m}$ can be naturally interpreted as a measurable function $$
M\otimes \mathbf{1}_{\Omega}:[n]\times [m]\times \Omega\rightarrow \R$$
$$
(i,j,\omega)\mapsto (M\otimes \mathbf{1}_{\Omega})(i,j,\omega)=M_{i,j},
$$
where on $[n]$ and on $[m]$ we consider the uniform measure and $\Omega$ is a probability space (for example $[0,1]$ with the Lebesgue measure or $[s]$ for $s\in \N$ with the uniform measure).
\end{example}

\begin{remark}
Random rectangular matrices, which are measurable functions from a probability space $\Omega$ to $\R^{n\times m},$ can be interpreted as $P-$bivariables similarly to Example \ref{Examp:RandomMatrixPset}. 
\end{remark}

\subsection{Metric and convergence}

For the test functions $f_1,\ldots, f_{k}\in L^{\infty}(\Omega)$, we consider the $2k+1$ dimensional random vector
$$
(f_1 \otimes\mathbbm{1}_{\Omega_1} \otimes \mathbbm{1}_{\Omega_2},\mathbbm{1}_{\Omega_1}\otimes f_1 \otimes\mathbbm{1}_{\Omega_2},\ldots, f_k \otimes\mathbbm{1}_{\Omega_1} \otimes \mathbbm{1}_{\Omega_2},\mathbbm{1}_{\Omega_1}\otimes f_k \otimes\mathbbm{1}_{\Omega_2},W).
$$
and its distribution, which we call the measure generated by $W$ through the test functions $f_1,\ldots,f_k,$ and we denote with

\begin{equation}\label{DefSf1f2fk}
    \begin{aligned}
        &S(f_1,\ldots,f_k,W)=\\
        &\mathcal{L}(f_1 \otimes\mathbbm{1}_{\Omega_1} \otimes \mathbbm{1}_{\Omega_2},\mathbbm{1}_{\Omega_1}\otimes f_1 \otimes\mathbbm{1}_{\Omega_2},\ldots, f_k \otimes\mathbbm{1}_{\Omega_1} \otimes \mathbbm{1}_{\Omega_2},\mathbbm{1}_{\Omega_1}\otimes f_k \otimes\mathbbm{1}_{\Omega_2},W)\in \mathcal{P}(\R^{2k+1}).
    \end{aligned}
\end{equation}

\begin{definition}\label{DefSymPset}
    We will say that a $P-$variable $W$ on $\Omega_1\times \Omega_1 \times \Omega_2$ is \emph{symmetric} if for any $f_1,\ldots, f_k \in L^{\infty}(\Omega_1)$

$$
\begin{aligned}&\mathcal{L}(f_1 \otimes\mathbbm{1}_{\Omega_1} \otimes \mathbbm{1}_{\Omega_2},\mathbbm{1}_{\Omega_1}\otimes f_1 \otimes\mathbbm{1}_{\Omega_2},\ldots, f_k \otimes\mathbbm{1}_{\Omega_1} \otimes \mathbbm{1}_{\Omega_2},\mathbbm{1}_{\Omega_1}\otimes f_k \otimes\mathbbm{1}_{\Omega_2},W)=\\
&\mathcal{L}(\mathbbm{1}_{\Omega_1}\otimes f_1 \otimes\mathbbm{1}_{\Omega_2},f_1 \otimes\mathbbm{1}_{\Omega_1} \otimes \mathbbm{1}_{\Omega_2},\ldots, \mathbbm{1}_{\Omega_1}\otimes f_k \otimes\mathbbm{1}_{\Omega_2},f_k \otimes\mathbbm{1}_{\Omega_1} \otimes \mathbbm{1}_{\Omega_2},W)
\end{aligned}$$
\end{definition}

\begin{remark}
Our definition of symmetry for a $P-$variable $W$ on $\Omega_1\times \Omega_1\times \Omega_2$ is weaker than requiring $$
W(x_1,x_2,x_{12})=W(x_2,x_1,x_{12})
$$ for almost every $x_1,x_2\in \Omega_1$ and $x_{12}\in \Omega_2$.    
\end{remark}
A particular case of a $P-$variable of special interest is the following.
\begin{definition}
A \emph{graph-set} is a symmetric $P-$variable taking values in $\{0,1\}$. 
\end{definition}

In particular, graph-sets are the $P-$variable version of a graph as the following example shows.

\begin{example}\label{ExGraphasGraphset}
Let's consider a simple graph $G$ and its adjacency matrix $A.$ The adjacency matrix $A\in \R^{n\times n}$ of a simple graph $G$ is a symmetric square matrix with entries taking values in $\{0,1\}.$ A square matrix can be interpreted as a $P-$variable as explained in Example \ref{Examp:MatrixasPset}. Therefore, graphs can be interpreted as graph-sets.
\end{example}

\begin{remark}
Similarly, we can also naturally consider bipartite graphs as $\{0,1\}-$valued $P-$bivariables considering their matrix representation as $P-$bivariables, see Example \ref{Examp:RectangMatrixPsets}.
\end{remark}

\begin{example}
\label{ExGraphasPsets}
Similarly to Example \ref{ExGraphasGraphset}, we can consider a weighted graph $H$ and its adjacency matrix $A.$ The adjacency matrix $A$ can be interpreted as a $P-$variable (Example \ref{Examp:MatrixasPset}). Therefore, a weighted graph $H$ can be considered as a $P-$variable interpreting  its adjacency matrix $A$ as a $P-$variable (see again Example \ref{Examp:MatrixasPset}). 
\end{example}

\begin{remark}\label{remarkRandGraphRepPset}
The examples above can be extended to random graphs and random weighted graphs using Example \ref{Examp:RandomMatrixPset} instead of Example \ref{Examp:MatrixasPset}.
\end{remark}
We now define the set of measures generated by a $P-$variable $W.$ For technical reasons, we will consider only measures generated by functions in $L_{[-1,1]}^{\infty}(\Omega_1)$, i.e. functions taking values between $-1$ and $+1$ almost everywhere. Therefore, we define the $k-$\emph{profile} of $W$,  $S_k(W)$, as the set of measures generated by $W$ by all test functions in  $L_{[-1,1]}^{\infty}(\Omega_1)$. This is 
\begin{equation}\label{Eq:DefKProfile}
    \begin{aligned}
S_k(W)=\bigcup_{f_1\ldots,f_k\in L_{[-1,1]}^{\infty}(\Omega_1)}&\{\mathcal{S}
(f_1,\ldots, f_k,W)\}.
\end{aligned}
\end{equation}
To compare different sets of probability measures we will use the Hausdorff metric $d_H$ (see Definition	\ref{HausdorffDist} below) where the space of probability measures $\mathcal{P}(\mathbb{R}^{2k+1})$ endowed with the Levy-Prokhorov metric (recall Definition \ref{LevyProk}) $d_{\mathcal{LP}}$ is considered as the underlying metric space. 

Therefore, to compare sets of measures we introduce the following distance.

\begin{definition}\label{HausdorffDist}
Let $(M, d)$ be a metric space. For each pair of non-empty subsets $X \subset M$ and $Y \subset M$, the \emph{Hausdorff distance} between $X$ and $Y$ is defined as
$$
d^{\text {Haus }}(X, Y):=\max \left\{\sup _{x \in X} d(x, Y), \sup _{y \in Y} d(X, y)\right\},
$$

 where $d(a, B)=\inf _{b \in B} d(a, b)$ quantifies the distance from a point $a \in M$ to the subset $B \subseteq M$.
\end{definition}

As anticipated, we will consider the space of measures $\cP(\mathbb{R}^{m})$ with the Lévy-Prokhorov metric as metric space $(M,d)=(\cP(\mathbb{R}^{m}),d_{\mathrm{LP}})$. We will denote with $d_H$ the Hausdorff metric associated with the metric space $(\cP(\mathbb{R}^{m}),d_{\mathrm{LP}}).$ We will consider the Hausdorff metric again on other metric spaces later in this work.

\begin{remark}
We observe that $d_{H}(X, Y)=0$ if and only if $\overline{X}=\overline{Y}$, where with $\overline{X}$ and $\overline{Y}$ we indicate the closure with respect to $d_{\mathcal{LP}}$ of the sets $X$ and $Y.$ Thus, $d_{H}$ is a pseudometric for all sets in $\mathcal{P}\left(\mathbb{R}^{k}\right)$, and it is a metric for the closed sets.
\end{remark}

Since the Lévy-Prokhorov metric is bounded by 1, we observe that the Hausdorff metric for sets of probability measures is also bounded by $1$.\newline

We are now ready to define the pseudo-metric we are interested in. Let's consider two $P-$variables (possibly defined on different probability spaces)

 $$W:\Omega_1\times \Omega_1 \times \Omega_2\rightarrow \R$$
    $$(x_1,x_2,x_{12})\mapsto W(x_1,x_2,x_{12})$$

and  $$U:\widetilde\Omega_1\times \widetilde\Omega_1 \times \widetilde\Omega_2\rightarrow \R$$
    $$(y_1,y_2,y_{12})\mapsto W(y_1,y_2,y_{12}).$$
\begin{definition}[Metrization of $P-$variables convergence] For the two $P-$variables $W, U$ the \emph{$P-$variables metric} is
$$
d_{M}(W, U)=\sum_{k=1}^{\infty} 2^{-k} d_{H}\left(\mathcal{S}_{k}(W), \mathcal{S}_{k}(U)\right).
$$
\end{definition}

Moreover, we will say that a sequence of $P-$variables $(W_n)_{n\in\N}$ is a Cauchy sequence if the sequence is Cauchy in $d_M$. \newline

Let $\Omega_1,\Omega_2$ and $\Omega_3$ be probability spaces. For a $P-$variable $W\in \cB(\Omega_1\times \Omega_1\times \Omega_2)$ and a measure-preserving map $\varphi$ from $\Omega_3$ to $\Omega_1$ we denote with $W^{\varphi}$ the $P-$variable $W^{\varphi}\in \cB(\Omega_3\times \Omega_3\times \Omega_2)$ defined for almost every $x_1,x_2\in \Omega_3$ and $x_{12}\in \Omega_2$ as $$
W^{\varphi}(x_1,x_2,x_{12})=W(\varphi(x_1),\varphi(x_2),x_{12}).
$$
Let $\Omega, \widetilde \Omega, \Omega_1,\Omega_2,\Omega_3$ and $\Omega_4$ be probability spaces. We will say that two $P-$variables $U\in \cB(\Omega_1\times \Omega_1\times \Omega_2)$ and $W\in \cB(\Omega_3\times \Omega_3\times \Omega_4)$ are \emph{weakly isomorphic}, denoted with $U\sim W$ if there exist measure-preserving maps $\varphi,$ from $\Omega$ to $\Omega_1,$ and $\psi,$ from $\Omega$ to $\Omega_3,$ such that for almost all $x_1,x_2\in \Omega$ we have $$\mathcal{L}(W^{\varphi}(x_1,x_2,\cdot))=\mathcal{L}(U^{\psi}(x_1,x_2,\cdot)).$$
We observe that, for each positive integer $k$, the $k-$profile is invariant up to weak isomorphism, i.e. for weakly isomorphic $P-$variables $U$ and $W,$ i.e.\ $U\sim W,$ as before we have 
\begin{equation}\label{Eq:EquivKprofWeakIso}
    S_k(U)=S_k(W).
\end{equation}

Thus we also have that for a $P-$variable  $V \in \cB(\Omega_1\times \Omega_1\times \Omega_2)$
$$
S_k(V)=S_k(V^{\varphi})
$$
for any measure-preserving map $\varphi$ from $\Omega_1$ to $\Omega_1.$

\begin{remark}\label{RMK:InvarianceKProfiles}
From \eqref{Eq:EquivKprofWeakIso}, we observe that the $P-$variables metric $d_M$ is invariant up to $P-$variables weak isomorphism. 
\end{remark}

\begin{remark}\label{RmkDistanseBBBisets}
Similarly to the $P-$variables metric $d_M$ for $P-$variables we can define a metric to compare $P-$bivariables. We just have to substitute the $2k+1$ random vectors for functions $f_1,\ldots, f_{k}\in L^{\infty}(\Omega)$ and a $P-$variable $W\in \cB(\Omega\times\Omega\times \widetilde\Omega)$ 
$$
(f_1 \otimes\mathbbm{1}_{\Omega_1} \otimes \mathbbm{1}_{\Omega_2},\mathbbm{1}_{\Omega_1}\otimes f_1 \otimes\mathbbm{1}_{\Omega_2},\ldots, f_k \otimes\mathbbm{1}_{\Omega_1} \otimes \mathbbm{1}_{\Omega_2},\mathbbm{1}_{\Omega_1}\otimes f_k \otimes\mathbbm{1}_{\Omega_2},W)
$$
with the $2k+1$ dimensional random vector for functions $f_1,\ldots, f_{k}\in L^{\infty}(\Omega_1)$ and $g_1,\ldots,g_k\in L^{\infty}(\Omega_2)$ and a $P-$bivariable $U$ defined on $\Omega_1\times \Omega_2 \times \Omega_3,$ for $\Omega_1,\Omega_2, \Omega_3$ probability spaces,
$$
(f_1 \otimes\mathbbm{1}_{\Omega_1} \otimes \mathbbm{1}_{\Omega_2},\mathbbm{1}_{\Omega_1}\otimes g_1 \otimes\mathbbm{1}_{\Omega_2},\ldots, f_k \otimes\mathbbm{1}_{\Omega_1} \otimes \mathbbm{1}_{\Omega_2},\mathbbm{1}_{\Omega_1}\otimes g_k \otimes\mathbbm{1}_{\Omega_2},W).
$$

In particular, we consider the probability measures on $\R^{2k+1}$
$$
\begin{aligned}
&SB(f_1,\ldots,f_k,g_1,\ldots,g_k,W)=\\
&\mathcal{L}(f_1 \otimes\mathbbm{1}_{\Omega_2} \otimes \mathbbm{1}_{\Omega_3},\mathbbm{1}_{\Omega_1}\otimes g_1 \otimes\mathbbm{1}_{\Omega_3},\ldots, f_k \otimes\mathbbm{1}_{\Omega_2} \otimes \mathbbm{1}_{\Omega_3},\mathbbm{1}_{\Omega_1}\otimes g_k \otimes\mathbbm{1}_{\Omega_3},W)\in \mathcal{P}(\R^{2k+1}).
\end{aligned}
$$
instead of the measures $S(f_1,\ldots,f_k,W)\in \cP(\R^{2k+1})$ in the definition of the $P-$variables metric $d_M.$
\end{remark}

\subsection{An equivalent formulation}
In this section, we formulate an equivalent formulation of $P-$variables convergence. In particular, we will show that we can be more conservative than with the test functions considered in the $k-$profiles. 

\begin{definition}
A function partition of $(\Omega_1,\mathcal{A},\mu)$ is a set $\{f_i\}_{i=1}^k=\{f_1,\ldots,f_k\}$ of $\{0,1\}-$ valued measurable functions on $\Omega$ such that $\sum_{i=1}^k f_i=\mathbbm{1}_{\Omega_1}$. 
\end{definition}

We observe that $f_1,\ldots,f_k$ is a function partition if and only if there is a measurable partition $\mathcal{P}=\{P_1,P_2,\dots,P_k\}$ of $\Omega_1$ such that $f_i=\mathbbm{1}_{P_i}$. Let $M_k$ denote the set of probability measures $\mu$ on $\mathbb{R}^{2k+1}$ such that the marginal on the even coordinates of $\mu$ is concentrated on $\{e_i\}$ for a certain $i\in [k]$, where $e_i\in\mathbb{R}^{k}$ is the vector/point with $1$ at the $i$-th coordinate and $0$ everywhere else. Let $W\in\cB(\Omega_1,\cA,\mu)$ and recall the definitions \eqref{DefSf1f2fk} of the measure $\mathcal{S}(f_1,\ldots, f_k,W)$  and \eqref{Eq:DefKProfile} of the $k-$profile $\mathcal{S}_k(W)$. We have that $\{f_i\}_{i=1}^k$ is a function partition if and only if $\mathcal{S}(f_1,\ldots, f_k,W)\in M_k$. Therefore, we define 
\begin{equation}
    \mathcal{S}'_k(W)=\mathcal{S}_k(W)\cap M_k.
\end{equation} In other words, $\mathcal{S}_k'(W)$ is the set of all probability measures $\mathcal{S}(f_1,\ldots, f_k,W)$, where $f_1,\ldots,f_k$ is a function partition.

The next theorem gives a useful equivalent formulation of $P-$variables convergence for $P-$variables. 

\begin{theorem}\label{partthm} Let $(W_n)_n$ be a sequence of $P-$variables. Then $(W_n)_n$ is $P-$variables convergent if and only if for every integer $k\geq 1$ the sequence $(\mathcal{S}'_k(W_n))_{n}$ is a Cauchy sequence in the Hausdorff metric $d_H.$
\end{theorem}

Theorem \ref{partthm} follows directly from the two lemmas below.  The proofs of the two lemmas are similar in spirit to the proofs of Lemma  7.4 and 7.5 in \cite{backhausz2018action}.

\begin{lemma}\label{partlem} Let $W$ be a $P-$variable. For each $\varepsilon>0$ and integer $k\geq 1$ there exists $\delta>0$ such that if  $\mu\in\mathcal{S}_k(W)$ satisfies $d_{\mathcal{LP}}(\mu,M_k)\leq\delta$, then there exists $\mu_2\in\mathcal{S}'_k(W)$ with $d_{\mathcal{LP}}(\mu,\mu_2)\leq\varepsilon$. Furthermore, for two $P-$variables $W,U$ such that $d_H(\mathcal{S}_k(W),\mathcal{S}_k(U))\leq\delta/2$, we have $$d_H(\mathcal{S}'_k(W),~\mathcal{S}'_k(U))\leq 2\varepsilon+d_H(\mathcal{S}_k(W),\mathcal{S}_k(U)).$$
\end{lemma}

\proof
We now show the first claim of the lemma. Let $\Omega_1$ and $\Omega_2$ probability spaces, $W$ be a $P-$variable $W\in \cB(\Omega_1\times \Omega_1\times \Omega_2)$  and $\mu\in S_k(W)$ such that $d_{\mathcal{LP}}(\mu,M_k)\leq\delta\leq \varepsilon.$  Thus, $\mu=\mathcal{S}(f_1,\ldots,f_k,W)$ for some $f_1,\ldots,f_k\in L_{[-1,1]}^{\infty}(\Omega_1)$. Moreover, since $d_{\mathcal{LP}}(\mu,M_k)\leq \delta$, we can find a probability measure in $M_k$ with distance at most $\delta$ from $\mu$.  Therefore, using Lemma \ref{LemmApproxFunctPartMeas} and Lemma \ref{LemmIneqProkh} together we get that there exist $\{0,1\}$-valued functions  $g_1,g_2,\dots,g_k$ in $L^{\infty}(\Omega_1)$ such that $\sum_{i=1}^k g_i=1_{\Omega_1}$ and $\max_{i\in [k]}\|f_i-g_i\|_p\leq C_k\delta$ for $1 \leq p < \infty$. It follows that $\mu_2=\mathcal{S}(v_1,\ldots,v_k,W)$ is in $\mathcal{S}'_k(W)$ and, by Lemma \ref{coupdist},  that (if $\delta$ is small enough) $d_{\mathcal{LP}}(\mu, \mu_2)\leq \varepsilon.$

Let now $W$ and $U$ be two $P-$variables such that $d_H(\mathcal{S}_k(W),\mathcal{S}_k(U))\leq\delta/2$. We now show the second statement of the lemma. By the definition of Hausdorff distance \ref{HausdorffDist} and of supremum, for any $\nu\in \mathcal{S}'_k(W)\subset \mathcal{S}_k(W)$ there exists  $\nu_2\in\mathcal{S}_k(U)$ with $d_{\mathcal{LP}}(\nu_2,\nu)\leq d_H(\mathcal{S}_k(U),\mathcal{S}_k(W))+\delta/2 \leq\delta.$ By definition, we know that $\nu \in S^{\prime}_k(W)\subseteq M_k,$ hence $d_{\mathcal{LP}}(\nu_2,M_k)\leq \delta$. Applying the first statement of the lemma to $\nu_2$ and $U$, we conclude that there exists $\nu_3\in\mathcal{S}'_k(U)$ with $d_{\mathcal{LP}}(\nu_2,\nu_3)\leq\varepsilon$. By the triangular inequality and $\delta \leq \varepsilon$, we obtain that $d_{\rm LP}(\nu_3,\nu)\leq d_{\mathcal{LP}}(\nu_3,\nu_2)+ d_{\mathcal{LP}}(\nu_2,\nu)\leq 2\varepsilon+d_H(\mathcal{S}_k(U),\mathcal{S}_k(W)).$
\endproof

\begin{remark}
The previous Lemma can be modified to obtain a quantitative bound, i.e. an explicit inequality independent of $\delta$.
\end{remark}

\begin{lemma}\label{LemmEquivForImplMetric}
For every $\varepsilon>0$ and integer $k\geq 1$ there exists $\delta>0$ and $k'\in\mathbb{N}$ such that if $W$ and $U$ are $P-$variables and $d_H(\mathcal{S}'_{k'}(W),\mathcal{S}'_{k'}(U))\leq\delta$, then $d_H(\mathcal{S}_k(W),\mathcal{S}_k(U))\leq\varepsilon$.
\end{lemma}
\proof
For any $P-$variable $W$ there are $\Omega_1$ and $\Omega_2$ probability spaces such that $W\in \cB(\Omega_1\times \Omega_1\times \Omega_2).$ For any measure $\mu\in\mathcal{S}_k(W)$ there are functions $f_i\in L^\infty_{[-1,1]}(\Omega_1)$ for $i\in [k],$  such that $\mu=\mathcal{S}(f_1,\ldots, f_k,W)$. For $n\in \N$ and $f\in L^p(\Omega_1)$ let $[f]_n$ denote the composition of $f$ with the function $x\mapsto \lceil nx\rceil/n$, which takes values in $n^{-1}\Z$. We notice that $\|f-[f]_n\|_p\leq n^{-1}$ for any $1\leq p \leq \infty$. Moreover, for every $i\in[k]$ the level sets of $[f_i]_n$ partition $\Omega_1$ into at most $2n$ measurable sets. By taking common refinements of the level sets of the functions $[f_i]_n,$ we get a partition $\{P_i\}_{i=1}^N$ of $\Omega_1$ into $N_n\leq (2n)^k$ measurable sets. The function $[f_i]_n$ is measurable by construction in this partition for every $i\in [k]$. Therefore, there exists $\{a_{i,j}\}_{i\in [k],j\in[N]}$ real numbers between $-1$ and $1$ such that for every $i\in[k]$ we have $[f_i]_n=\sum_{j\in [N]}a_{i,j}1_{P_j}$. 

For $n$ large enough we obtain that $\mu=\mathcal{S}(f_1,\ldots,f_k,W)$ and $\mu=\mathcal{S}([f_1]_n,\ldots,[f_k]_n,W)$ are at distance $d_{\mathcal{LP}}$ smaller than $\varepsilon/2.$ Let's now consider $k^{\prime}>N_n$ and $U$ be  another $P-$variable such that $d_H(\mathcal{S}'_{k'}(W),\mathcal{S}'_{k'}(U))\leq\delta$ for $\delta>0$ small enough. Therefore, by the definition of Hausdorff distance \ref{HausdorffDist} and of the supremum, there exists a measurable partition $\{Q_i\}_{i=1}^{{N_n}}$ of $\widetilde\Omega_1$ such that 
$d_{\mathcal{LP}}(\mathcal{S}(1_{P_1},\ldots,1_{P_{N_n}},W),\mathcal{S}(1_{Q_1},\ldots,1_{Q_{N_n}},U))\leq d_H(\mathcal{S}'_{k'}(W),\mathcal{S}'_{k'}(U))+\delta\leq 2\delta,$ where in the first inequality we also used Lemma \ref{LemmIneqProkh}.

Let now $w_i=\sum_{j\in [{N_n}]}a_{i,j}1_{Q_j}$. If $\delta$ is small enough, then $\kappa=\mathcal{S}(w_1,\ldots, w_k,U)\in S_k(U)$ is arbitrarily close to $\mu=\mathcal{S}(f_1,\ldots,f_k,W)$ in $d_{\mathcal{LP}}$. We obtain that if $n$ is big enough and $\delta$ is small enough, then $d_{H}(S_k(W),S_k(U))\leq \varepsilon$. We remark that the estimates in the proof depend only on $\varepsilon$ and $k$.
\endproof

\section{Equivalence of P-variables and probability graphons}\label{Sec4:EquivPsetsProbGraphon}

\subsection{From P-variables to probability graphons}

Recall that $\Space$ is a measurable subset of $\R.$

There is a natural procedure to transform a $P-$variable $W$ on $\Omega_1\times \Omega_1 \times \Omega_2,$ the product probability space of $(\Omega_1,\cA_1,\P)$ and $(\Omega_2,\cA_2,\mu),$ into a probability graphon $\widetilde W$ from $\Omega_1 \times \Omega_1$ to $\mathcal{P}(\Space)$. 
Given a $P-$variable $W$ we can define
$$
\widetilde W (x_1,x_2,E)=\mu(\{x_{12}: \ W(x_1,x_2,x_{12})\in E\})=\mathcal{L}(W(x_1,x_2,\cdot))(E)
$$
for any $E$ measurable subset of $\Space\subset\R.$
For every fixed measurable set $E\subset \Space \subset\R$ the function $\widetilde W (x_1,x_2,E)$ is a measurable function from the $\Omega_1\times \Omega_1$ to $\R$ by Fubini-Tonelli theorem (see for Example Theorem 18.2 in \cite{Bill86}).

\begin{remark}
    We underline that for any $P-$variable $W,$ there exists a unique probability graphon $\widetilde W$ associated to $W$ in this way (up to measure zero modifications). 
\end{remark}

Integrating away the last variable we obtain the real-valued graphon
\begin{equation}\label{Eq:DefContraction}
w(x_1,x_2)=\int_{\Omega_2}W(x_1,x_2,x_{12})\mathrm{d}\mu(x_{12})=\int_{\mathbb{R}}z\widetilde W(x_1,x_2,\mathrm{d}z).
\end{equation}
that we will call the \emph{contraction} of $W.$ This is taking the expectation of the random variable $W(x_1,x_2,\cdot)$ (only with respect to the last variable). Compare also Example \ref{rem:real-valued-kernels}.

Equality \eqref{Eq:DefContraction} holds as $\widetilde W(x_1,x_2,\cdot)$ is the push-forward of $\mu$ through the function $W(x_1,x_2,\cdot)$ (or equivalently the distribution of the random variable $W(x_1,x_2,\cdot)$). 

More generally, recall \eqref{EqPushForwardInt}, for a measurable function $g:\Space \subset\R \rightarrow \R$ and a measurable set $X_1\subset \Space \subset \R$ we have the equality 
\begin{equation}\label{Eq:IdentityPushForwardProbGraphPset}
    \int_{X_2(x_1,x_2)}g(W(x_1,x_2,x_{12}))\mathrm{d}\mu(x_{12})=\int_{X_1}g(z)\widetilde W(x_1,x_2,\mathrm{d}z)
\end{equation}

where $X_2(x_1,x_2)=W(x_1,x_2,\cdot)^{-1}(X_1).$

In particular, 
$$
\int_{\Omega_2}g(W(x_1,x_2,x_{12}))\mathrm{d}\mu(x_{12})=\int_{\Space}g(z)\widetilde W(x_1,x_2,\mathrm{d}z).
$$
The special case where $g$ is a continuous and bounded function $g\in \CbFunct$ is of special interest. In this case, we have
\begin{equation}\label{Eq:ProbGraphContFunIdent}
    \widetilde W[g](x_1,x_2)=\widetilde W(x_1,x_2,g)=\int_{\Space}g(z)\widetilde W(x_1,x_2,\mathrm{d}z)=\int_{\Omega_2}g(W(x_1,x_2,x_{12}))\mathrm{d}\mu(x_{12}).
\end{equation}
These quantities play a key role in the theory of probability graphons. In fact, using these quantities, one can express the unabelled cut metric (see \cite{abraham2023probabilitygraphons,zucal2024probabilitygraphonsrightconvergence}), homomorphism densities (see \cite{abraham2023probabilitygraphons} and and Section \ref{Sec5TranslatResultsExamp}) and overlay functionals (see \cite{zucal2024probabilitygraphonsrightconvergence} and Section \ref{Sec5TranslatResultsExamp}) for probability graphons. Using identity \eqref{Eq:ProbGraphContFunIdent} we will show in Section \ref{Sec5TranslatResultsExamp} how to translate results from the theory of probability graphons to $P-$variables.

Using \eqref{Eq:IdentityPushForwardProbGraphPset} we obtain 
\begin{equation*}
\begin{aligned}
 &   \widetilde  W  (\Omega_1\times \Omega_1; E)=\\
&\int_{\Omega_1^2}  \widetilde W(x_1,x_2;E) \ \mathrm{d} \mathbb{P}(x_1)\mathrm{d} \mathbb{P}(x_2)  =\\
&\int_{\Omega_1^2}\int_{\Space}\mathbbm{1}_{E}(z)  \widetilde W(x_1,x_2;\drv z) \ \mathrm{d} \mathbb{P}(x_1)\mathrm{d} \mathbb{P}(x_2) =\\
&\int_{\Omega_1^2}\int_{\Omega_2}\mathbbm{1}_{E}(W(x_1,x_2,x_{12}))\mathrm{d} \mathbb{P}(x_1)\mathrm{d} \mathbb{P}(x_2) \mu ( \drv x_{12})\\
&(\mathbb{P}\otimes \mathbb{P}\otimes \mu)\left(W\in E\right).
\end{aligned}
\end{equation*}

Therefore, we observe that 
\begin{equation}\label{Eq:EqualityPsetProbGraphonMeasuresSets}
    \widetilde  W  (\Omega_1^2; E)=(\mathbb{P}\otimes \mathbb{P}\otimes \mu)\left(W\in E\right)=\cS(\mathbbm{1}_{\Omega_1},W)(\{1\}\times E).
\end{equation}
Another identity that will be important later is the following. From Fubini-Tonelli theorem (see for Example Theorem 18.2 in \cite{Bill86}) we have that for $S,T$ measurable disjoint subsets of $\Omega_1,$ a measurable subset $E$ of $\Space$ and the $5-$dimensional random vector $ \mathcal{S}(\mathbbm{1}_S,\mathbbm{1}_T,W)$ we have:

\begin{equation}\label{identityProbGraphPsets}
\begin{aligned}
& \mathcal{S}(\mathbbm{1}_S,\mathbbm{1}_T,W)( \{(1,0,0,1)\}\times E)\\
&\mathcal{L}\left((\mathbbm{1}_S\otimes\mathbbm{1}_{\Omega_1} \otimes \mathbbm{1}_{\Omega_2},\mathbbm{1}_{\Omega_1}\otimes \mathbbm{1}_S \otimes\mathbbm{1}_{\Omega_2}, \mathbbm{1}_T \otimes\mathbbm{1}_{\Omega_1} \otimes \mathbbm{1}_{\Omega_2},\mathbbm{1}_{\Omega_1}\otimes \mathbbm{1}_T \otimes\mathbbm{1}_{\Omega_2},W)\right)(\{(1,0,0,1)\}\times E)\\
&=\mathbb{P}\otimes \mathbb{P}\otimes \mu (\{(x_1,x_2,x_{12}): \ (\mathbbm{1}_S(x_1),\mathbbm{1}_T(x_2), W(x_1,x_2,x_{12}))\in \{(1,1)\}\times E\})=\\
& \int_{S\times T}\mu(\{x_{12}: \ W(x_1,x_2,x_{12})\in E\})\mathrm{d} \mathbb{P}(x_1)\mathrm{d} \mathbb{P}(x_2)=\\
& \int_{S\times T} \left(  \widetilde W(x_1,x_2,E)\right)\mathrm{d} \mathbb{P}(x_1)\mathrm{d} \mathbb{P}(x_2)=\widetilde W(S\times T,E).
\end{aligned}
\end{equation}

A similar formula holds for any pairwise disjoint sets $S_1,\ldots, S_k$ and $\cS(\mathbbm{1}_{S_1},\ldots, \mathbbm{1}_{S_k},W).$

\subsection{From probability graphons to P-variables}

We explain now how to go from a probability graphon $\widetilde{W}$ to a $P-$variable $W.$ First of all, we recall the definition of the distribution function $F$ of a measure $\mu$ (or equivalently of a random variable $X$ with distribution $\mu$):
$$F(x)=F_{\mu}(x)=F_{X}(x)=\mu((-\infty, x]).$$
The distribution function $F$ is a non-decreasing function.
We define also the (generalized) inverse distribution function for the measure $\mu$ (or the random variable $X$ with distribution $\mu$) with
$$
\begin{aligned}
   G(p)&=G_{\mu}(p)=G_{X}(p)=
   \inf\{x\in\R: \ F_{\mu}(x)\geq p\} = \inf\{x\in\R: \ \mu((-\infty,x])\geq p\}.
\end{aligned}$$

The inverse distribution $G$ is also a non-decreasing function.

For a probability graphon $\widetilde W$ we, therefore, define the function $W$ (denoted also as $G_{\widetilde W}$ to recall the dependency on $\widetilde W$ explicitly) for every $(x_1,x_2,x_3)\in \Omega_1\times \Omega_1 \times [0,1]$ as

$$G_{\widetilde W}(x_1,x_2,x_3)= W(x_1,x_2,x_{12})=G_{\widetilde W(x_1,x_2,\cdot)}(x_{12}),$$
where we recall that $G_{\widetilde W(x_1,x_2,\cdot)}(x_{12})$ is the inverse distribution function for the measure $\widetilde W(x_1,x_2,\cdot).$
\begin{lemma}\label{MeasProbGraphonFromPset}
For a probability graphon $\widetilde W$ on $\Omega_1,$ the associated $W$ (also denoted as $G_{\widetilde W}$) 

$$ W(x_1,x_2,x_{12})=G_{\widetilde W(x_1,x_2,\cdot)}(x_{12})$$
is a measurable function. In particular, $W=G_{\widetilde W}$ is a $P-$variable $W\in \cB(\Omega_1\times \Omega_1\times [0,1])$.
\end{lemma}
\proof

We first define the function 

$$
G:\mathcal{P(\R)}\times [0,1]\rightarrow \R
$$

$$
(\mu,p)\mapsto G(\mu,p)=\inf\{x\in\R: \ F_{\mu}(x)\geq p\} = \inf\{x\in\R: \ \mu((-\infty,x])\geq p\}= G_{\mu}(p).$$
For $r\in \R$ let $R_r=(-\infty,r].$
To show the measurability of this function it is enough to show that for any $r\in \R$  the set $G^{-1}(R_r)$ is measurable in $\mathcal{P(\R)}\times [0,1]$ with the $\sigma-$algebra of the product topology.

We observe that $$\begin{aligned}
    G^{-1}(R_r)=\{(\mu,p):\ G(\mu,p)\leq r\}=\{(\mu,p):\ F_{\mu}(r)\geq p\}=\{(\mu,p):\ \mu(R_r)-p\geq 0\},
\end{aligned}$$
where the second equality follows as $F_{\mu}(r)\geq p$ if and only if $G_{\mu}(p)\leq r.$
By definition of the topology of weak convergence on the space of measures the functional from $\mathcal{P}(\R)\rightarrow \R$ 

$$
E_{R_r}(\mu)=\mu(R_r)$$
is measurable (see Remark 2.4 in \cite{abraham2023probabilitygraphons} for example). 

Therefore, the function from $\R \times \mathcal{P}(\R)$ to $\R$

$$(p,\mu)\mapsto E_{R_r}(\mu)$$

is obviously also measurable. Also the function from $\R \times \mathcal{P}(\R)$ to $\R$

$$
(p,\mu)\mapsto -p
$$
is trivially measurable. Therefore, their sum 

$$
S:\R \times \mathcal{P}(\R) \rightarrow \R$$
$$
(p,\mu)\mapsto E_{R_r}(\mu)-p
$$
is also measurable.

Therefore, the set 

$$
S^{-1}([0,+\infty))=\{(\mu,p):\ E_{R_r}(\mu)-p\geq 0\}=\{(\mu,p):\ \mu(R_r)-p\geq 0\}=G^{-1}(R_r)$$
is measurable. Therefore, we can conclude that $G$ is measurable.

Now we define also the function 

$$
(\widetilde W \times id):\Omega_1^2\times [0,1]\rightarrow \mathcal{P}(\R)\times [0,1]
$$
$$
((x_1,x_2),p)\mapsto (\widetilde W(x_1,x_2,\cdot),p).
$$

This function is also measurable as $\widetilde W$ and the identity functions are measurable and the two functions depend on different variables.

We get that $W$ is the composition of $(\widetilde W \times id)$ with $G,$ i.e.\ 

$$
W=(G\circ (\widetilde W \times id)).$$
Thus, for all $x_1,x_2\in \Omega_1$ and $x_{12}\in [0,1]$ we have
$$W(x_1,x_2,x_{12})=G_{\widetilde W(x_1,x_2,\cdot)}(x_{12})=(G\circ (\widetilde W \times id))(x_1,x_2,x_{12}).$$

Therefore, $W$ is a measurable function as it is the composition of two measurable functions.
\endproof

\begin{remark}
    We observe that for a probability graphon $\widetilde W$  we have that $\widetilde W=\widetilde{G_{\widetilde W}},$ where $G_{\widetilde W}$ is the $P-$variable associated with $\widetilde W$ constructed above and $\widetilde{G_{\widetilde W}}$ is the probability graphon associated with the $P-$variable $G_{\widetilde W}$ constructed in the previous section.
\end{remark}

\begin{remark}
We underline that for a probability graphon $\widetilde W$ there are typically many $P-$variables $V$ such that $\widetilde W=\widetilde V.$ 
\end{remark}

\subsection{Tightness}
In this section, we will introduce a characterization of compact sets of probability graphons and $P-$variables with a tightness criterion similar to the tightness criterion for probability measures and random variables (Prokhorov theorem \cite[Theorem 5.1]{billingsley1968convergence}).

\begin{definition}
A set of probability graphons $\cK$ on the probability space $(\Omega_1,\cA_1,\P_1)$ is \emph{tight} if the set of measures $\{ M_{\widetilde W} : \widetilde W\in\mathcal{K} \} $ on $\Space\subset \R$ is tight, where $M_{\widetilde W}$ for the probability graphon $\widetilde W$ is the measure
\begin{equation}
  \label{eq:def-MW}
M_{\widetilde W}(\drv z) 
=  \widetilde W  (\Omega_1\times \Omega_1; \drv z)
= \int_{\Omega_1^2} \widetilde  W(x,y;\drv z) \ \drv x \drv y.
\end{equation}
\end{definition}

\begin{remark}
In particular, if $\Space\subset \R$ is compact then the space of $\cP(\Space)-$valued probability graphons is tight. 
\end{remark}

We present here the related notion of tightness for $P-$variables.

\begin{definition}
Let $(\Omega_1,\cA_1,\P_1)$ and $(\Omega_2,\cA_2,\mu)$ two probability spaces.  A set of $P-$variables $\cK\subset \cB(\Omega_1\times\Omega_1\times \Omega_2)$ is \emph{tight} if and only if $\mathcal{K}$ is tight considering the $P-$variables $W\in \mathcal{K}$ as random variables on the probability space $\Omega_1\times \Omega_1\times \Omega_2,$ i.e.\ for every $\varepsilon>0$ there exists $K>0$ such that for every $W\in \mathcal{K}$ we have
\begin{equation}\label{tightCondit}
    (\mathbb{P}_1\otimes \mathbb{P}_1\otimes \mu)\left(|W|>K\right)<\varepsilon.
\end{equation}
\end{definition}

\begin{remark}
One can easily generalize the definition of tightness to families of $P-$variables defined on different probability spaces as the definition of tightness of families of random variables. 
\end{remark}
We explain now how the notions of tightness for probability graphons and $P-$variables are interconnected.
For $W\in \cB(\Omega_1\times\Omega_1\times \Omega_2),$ using \eqref{Eq:EqualityPsetProbGraphonMeasuresSets} we obtain 
\begin{equation}\label{Eq:TightnessEquiv}
\begin{aligned}
  M_{\widetilde W}(\drv z) 
=\widetilde  W  (\Omega_1\times \Omega_1; [-K,K]^c)=(\P\otimes \P\otimes \mu)\left(|W|>K\right).
\end{aligned}
\end{equation}
The following lemma follows directly from \eqref{Eq:TightnessEquiv}.

\begin{lemma}\label{LemmEqTightness}
 A set of $P-$variables $\cK\subset \cB(\Omega_1\times\Omega_1\times \Omega_2)$ is \emph{tight} if and only if $\widetilde \cK=\{\widetilde W: \ W\in \cK\},$ i.e. the sets of probability graphons associated with the elements of $\cK,$ is tight. 
\end{lemma}

As anticipated, similarly to probability measures and random variables, tightness is strongly linked with compactness. We recall here Proposition 5.2 in \cite{abraham2023probabilitygraphons} which is stated for probability graphons on $[0,1].$

\begin{theorem}[Proposition 5.2 in \cite{abraham2023probabilitygraphons}]\label{ThmRelCompact}
Let $\cK$ be a set of probability graphons on $[0,1],$ i.e.\ $\cK\subset \UGraphon.$ The set $\cK$ is relatively compact for $\delta_{\square}$  if and only if it is tight.
\end{theorem}
\begin{remark}
In particular, if the probability graphons are $\cP(\Space)-$valued and $\Space\subset \R$ is compact then $\UGraphon,$ the space of probability graphons on $[0,1]$, is compact. 
\end{remark}

\subsection{Equivalence of unlabelled cut metric and P-variables metric}

We prove in this section the main result of this work: The topological equivalence of unlabelled cut-distance for probability graphons and the $P-$variables metric for tight sets of $P-$variables.

\begin{theorem}\label{TH:equivalence_dM_cut_dist}
Let $(W_n)_n$ be a tight sequence of $P-$variables $W_n\in \cB([0,1]\times [0,1]\times \Omega_2) $ and $(\widetilde W_n)_n$ be the sequence of probability graphons $\widetilde W_n$ associated with the $P-$variables $W_n.$ Convergence in the $P-$variables metric of $(W_n)_n$ and convergence in the unlabelled cut-distance $\delta_{\square}$ of $(\widetilde W_n)_n$ are equivalent.    
\end{theorem}

The proof follows directly from Lemma \ref{LemmUnlabelledCutImplies}, Theorem \ref{ThmEquivQuotients} and Corollary \ref{CorMetricImpliesUnlbCutDist} below.

\begin{lemma}\label{LemmKprofBoundCutdist}
Let $U,W\in\cB([0,1]\times [0,1]\times \Omega)$ and let $\widetilde U$ and $\widetilde W$ be the associated probability graphons and let $f_1,\ldots, f_k$ be a function partition on $[0,1]$. For every $\varepsilon>0,$ there exists $\delta>0$ such that if $d_{\square}(\widetilde U,\widetilde W)<\delta$ then 
$$d_{\mathcal{LP}}(S(f_1,\ldots,f_k,W),S(f_1,\ldots,f_k,U))<\varepsilon.$$
\end{lemma}
\proof
We will denote here $\Omega_1=[0,1]$ and $\Omega_2=\Omega.$
We have
$$
\begin{aligned}
& d_{\mathcal{LP}}(S(f_1,\ldots,f_k,W),S(f_1,\ldots,f_k,U))=\\
& d_{\mathcal{LP}}(\mathcal{L}(f_1 \otimes\mathbbm{1}_{\Omega_1} \otimes \mathbbm{1}_{\Omega_2},\mathbbm{1}_{\Omega_1}\otimes f_1 \otimes\mathbbm{1}_{\Omega_2},\ldots, f_k \otimes\mathbbm{1}_{\Omega_1} \otimes \mathbbm{1}_{\Omega_2},\mathbbm{1}_{\Omega_1}\otimes f_k \otimes\mathbbm{1}_{\Omega_2},W),\\
&\mathcal{L}(f_1 \otimes\mathbbm{1}_{\Omega_1} \otimes \mathbbm{1}_{\Omega_2},\mathbbm{1}_{\Omega_1}\otimes f_1 \otimes\mathbbm{1}_{\Omega_2},\ldots, f_k \otimes\mathbbm{1}_{\Omega_1} \otimes \mathbbm{1}_{\Omega_2},\mathbbm{1}_{\Omega_1}\otimes f_k \otimes\mathbbm{1}_{\Omega_2},U)).
\end{aligned}
$$

Moreover, for any $Q\subset \R^{2k}$ and $I\subset \R$ measurable subsets, we have 
\begin{equation}\label{eqDecompWeakMeasConv1}
 \begin{aligned}
&S(f_1,\ldots,f_k,W)(Q\times I)=\\
&\mathcal{L}(f_1 \otimes\mathbbm{1}_{\Omega_1} \otimes \mathbbm{1}_{\Omega_2},\mathbbm{1}_{\Omega_1}\otimes f_1 \otimes\mathbbm{1}_{\Omega_2},\ldots, f_k \otimes\mathbbm{1}_{\Omega_1} \otimes \mathbbm{1}_{\Omega_2},\mathbbm{1}_{\Omega_1}\otimes f_k \otimes\mathbbm{1}_{\Omega_2},W)(Q\times I)=\\
& \sum^k_{i,j=1}\mathcal{L}(f_1 \otimes\mathbbm{1}_{\Omega_1} \otimes \mathbbm{1}_{\Omega_2},\mathbbm{1}_{\Omega_1}\otimes f_1 \otimes\mathbbm{1}_{\Omega_2},\ldots, f_k \otimes\mathbbm{1}_{\Omega_1} \otimes \mathbbm{1}_{\Omega_2},\mathbbm{1}_{\Omega_1}\otimes f_k \otimes\mathbbm{1}_{\Omega_2},W)((Q\cap E_{i,j})\times I)=\\
& \sum^k_{i,j=1}\widetilde W(S_i\times S_j,I)\delta_{ij}(Q),
\end{aligned}   
\end{equation}
where 
$$\begin{aligned}
&E_{ij}=\{w\in \Omega_1\times \Omega_1 \times \Omega_2: \\
&\ f_1 \otimes\mathbbm{1}_{\Omega_1} \otimes \mathbbm{1}_{\Omega_2}=0,\mathbbm{1}_{\Omega_1}\otimes f_1 \otimes\mathbbm{1}_{\Omega_2}=0,\\
&\ldots,f_i \otimes\mathbbm{1}_{\Omega_1} \otimes \mathbbm{1}_{\Omega_2}=1,
\ldots, \mathbbm{1}_{\Omega_1}\otimes f_j \otimes \mathbbm{1}_{\Omega_2}=1,\ldots,\\
&\mathbbm{1}_{\Omega_1}\otimes f_k \otimes\mathbbm{1}_{\Omega_2}=0,\mathbbm{1}_{\Omega_1}\otimes f_k \otimes\mathbbm{1}_{\Omega_2}=0 \}
\end{aligned}$$

and 

$$\delta_{ij}(Q)=\begin{cases}
    1 & \text{ if } (e^{o}_i+e^{e}_j)\in Q \\

    0 & \text{ else },
\end{cases}$$

where $e^o_i\in \R^{2k}$  for $i\in [k]$ is the $2k-$dimensional vector with the entry $2i-1$ equal to $1$ and every other entry equal to $0$ and, similarly, $e^e_j\in \R^{2k}$ for $j\in [k] $ is the $2k-$dimensional vector with the entry $2j$ equal to $1$ and every other entry equal $0.$ Therefore $e^o_i+e_j^e$ denotes the  $2k-$dimensional vector with the entry $2i-1$ and the entry $2j$  equal to $1$ and every other entry equal to $0.$

Similarly, we obtain that 
\begin{equation}\label{eqDecompWeakMeasConv2}
    \begin{aligned}
S(f_1,\ldots,f_k,U)(Q\times I)=
 \sum^k_{i,j=1}\widetilde U(S_i\times S_j,I)\delta_{ij}(Q).
\end{aligned}
\end{equation}
Applying  Theorem \ref{BillingWeakConv} and using equalities \eqref{eqDecompWeakMeasConv1} and \eqref{eqDecompWeakMeasConv2} we obtain that if $d_{\mathcal{LP}}(\widetilde W(S_i\times S_j,\cdot),\widetilde U(S_i\times S_j,\cdot))$ is small enough for every $i,j\in [k]$ then
$$
\begin{aligned}
d_{\mathcal{LP}}(\mu_W,\mu_U)=
 d_{\mathcal{LP}}(S(f_1,\ldots,f_k,W),S(f_1,\ldots,f_k,U))\leq \varepsilon.\\
\end{aligned}
$$

In fact, for any $S,T\subset [k]$
\begin{equation}
    \begin{aligned}
&d_{\mathcal{LP}}( \sum_{i\in S,j\in T}\widetilde W(S_i\times S_j,\cdot),\sum_{i\in S,j\in T}\widetilde U(S_i\times S_j,\cdot))\leq \\
& |S||T|d_{\mathcal{LP}}\left(\frac{ \sum_{i\in S,j\in T}\widetilde W(S_i\times S_j,\cdot)}{|S||T|},\frac{\sum_{i\in S,j\in T}\widetilde U(S_i\times S_j,\cdot)}{|S||T|}\right)\leq\\
& k^2 \max_{i,j\in [k]}d_{\mathcal{LP}}\left(\widetilde W(S_i\times S_j,\cdot),\widetilde U(S_i\times S_j,\cdot)\right)
\end{aligned}
\end{equation}
where in the first inequality we used Lemma \ref{LemmaIneqScalingProkhorov} and in the second inequality we used Lemma \ref{LemmaQuasi-convProhorov}.

But from $d_{\square}(\widetilde W, \widetilde U)\leq \delta,$ we have $d_{\mathcal{LP}}(\widetilde W(S_i\times S_j,\cdot),\widetilde U(S_i\times S_j,\cdot))\leq \delta$ for every $i,j\in [k]$ as

$$
\begin{aligned}
&d_{\mathcal{LP}}(\widetilde W(S_i\times S_j,\cdot),\widetilde U(S_i\times S_j,\cdot))\leq\\
& \sup_{S,T} d_{\mathcal{LP}}(\widetilde W(S\times T,\cdot),\widetilde U(S\times T,\cdot))=\\
& d_{\square}(\widetilde W, \widetilde U)\leq \delta
\end{aligned}
$$
and this concludes the proof.

\endproof

\begin{lemma}\label{LemmUnlabelledCutImplies}
Let $(W_n)$ be a sequence of $P-$variables $W_n\in \cB([0,1]\times [0,1]\times \Omega) $ and $(\widetilde W_n)$ be the sequence of probability graphons $\widetilde W_n$ associated with the $P-$variables $W_n.$ Convergence in the unlabelled cut distance $\delta_{\square}$ of $\widetilde W_n$ implies $P-$variables convergence of $(W_n)$.
\end{lemma}
\proof
Let $W$ and $U$ be two $P-$variables $W,U\in \cB([0,1]\times [0,1]\times \Omega)$ and let $\widetilde W$ and $\widetilde U$ be the associated probability graphons. Let's assume $d_{\square}(\widetilde W, \widetilde U)\leq \delta.$ We show now that for every $\varepsilon>0$ there exists $\delta>0$ such that for $d_{\square}(\widetilde W, \widetilde U)\leq \delta$ we have
$$d_H(S^{\prime}_k(W),S^{\prime}_k(U)\leq \varepsilon$$for any positive integer $k.$

We have that for any measure $\mu_W \in S^{\prime}_k(W)$ there exist a function partition $f_1,\ldots,f_k,$ i.e.\ $f_i=\mathbbm{1}_{S_i}$ where $\{S_1,\ldots,S_k\}$ is a measurable partition of $\Omega_1,$ such that $\mu_W=S(f_1,\ldots,f_k,W)$ and let's define $\mu_U=S(f_1,\ldots,f_k,U).$

From Lemma \ref{LemmKprofBoundCutdist}, we have that choosing $\delta>0$ small enough, we obtain that as long as $d_{\square}(\widetilde W, \widetilde U)\leq \delta$ we have

$$d_H(S^{\prime}_k(W),S^{\prime}_k(U))\leq \varepsilon.$$
By Lemma \ref{LemmEquivForImplMetric} (or Theorem \ref{partthm}) this implies
$$
d_M(U,W)\leq \varepsilon^{\prime} .
$$

This concludes the proof because 

$$
d_M(U,W)\leq \varepsilon^{\prime}
$$
as long as 
$$
\inf_{\varphi \in \InvRelabel}d_{\square}(\widetilde W,\widetilde U^{\varphi})=\delta_{\square}(\widetilde W,\widetilde U)\leq \delta
$$

as the metric  $d_M$ is invariant with respect to measure preserving transformations (see \eqref{Eq:EquivKprofWeakIso} or Remark \ref{RMK:InvarianceKProfiles}).
\endproof

We now introduce an alternative characterization for convergence in the unlabelled cut metric $\delta_{\square}$ of probability graphons.

\textbf{Quotient sets of probability graphons:} Let $\lambda$ be the Lebesgue measure on $[0,1].$ Let $W$ be a probability graphon $W\in\Graphon$ and let $\mathcal{P}=\left\{S_1, \ldots, S_k\right\}$ be a partition of $[0,1]$ in $k$ measurable sets for a positive integer $k$, we define the quotient graph (or simply quotient) $W / \mathcal{P}$ as the measure edge-decorated and vertex-weighted graph on $[k]$, with node weights $\alpha_i(W / \mathcal{P})=\lambda\left(S_i\right)$ and as measure decoration of the edge $e=\{i,j\}$ the measure
$$
\beta_{i j}(W / \mathcal{P})=\frac{1}{\lambda\left(S_i\right) \lambda\left(S_j\right)} \int_{S_i \times S_j} W .
$$
To study right-convergence for probability graphons, we need to consider all quotient graphs of a given probability graphon at the same time. For a measure-valued kernel $W$ and probability distribution $\mathbf{a}$ on $[k],$ we denote by $\mathcal{Q}_{\mathbf{a}}(W)$ the set of quotients $L=W / \mathcal{P}$ with $(\alpha_i(L))_{i\in [k]}=(\alpha_1(L),\ldots,\alpha_k(L))=\mathbf{a}.$ We denote by $\mathcal{Q}_k(W)$ the set of all quotients for all $k$-partitions (measurable partition with $k$ elements).

We will consider these quotient sets as subsets of the space $$\R^k\times\cP(\Space)^{{k} \choose {2}}.$$
In particular, every quotient graph $W / \mathcal{P}$ in the quotient set $\mathcal{Q}_k(W)$ will be considered as an element of $\R^k\times\cP(\Space)^{{k} \choose {2}},$ where $\alpha=(\alpha_1(W / \mathcal{P}),\ldots,\alpha_k(W / \mathcal{P}))\in \R^k$ and $\beta(W / \mathcal{P})=(\beta_{i j}(W / \mathcal{P}))_{i,j\in [k]}\in \cP(\Space)^{{k} \choose {2}}.$

We now introduce a metric between graphs with weighted vertices and measure-decorated edges that we will employ to compare quotient graphs. 
For two vertex weighted and edge measure decorated graphs $G$ and $H$ on the vertex set $[k]$ we define the $d_1$ distance as 

$$
d_1\left(G, H\right)=\sum_{i\in[k]}|\alpha_i(G)-\alpha_i(H)|+\sum_{i, j \in[k]}d_{\mathcal{LP}}(\alpha_i(G)\alpha_j(G)\beta_{ij}(G),\alpha_i(H)\alpha_j(H)\beta_{ij}(H),\cdot)).
$$
By the definition of quotient graphs $L\in \mathcal{Q}_{\mathbf{a}}(W)\subset \mathcal{Q}_{k}(W) $ and $L^{\prime}\in \mathcal{Q}_{\mathbf{a}^{\prime}}(W) \subset\mathcal{Q}_k(U)$ we have
$$
d_1\left(L, L^{\prime}\right)=\left\|\mathbf{a}-\mathbf{a}^{\prime}\right\|_1+\sum_{i, j \in[k]}d_{\mathcal{LP}}(W(S_i \times S_j,\cdot),U(S_i^{\prime} \times S_j^{\prime},\cdot)).
$$

The convergence of probability graphons in unlabelled cut metric $\delta_{\square}$ is characterized by convergence of quotient sets in Hausdorff metric $d_{1,H},$ as proved in Theorem 4.26 in \cite{zucal2024probabilitygraphonsrightconvergence} that we recall below. This is the analogue for probability graphons of Theorem 12.12 in \cite{LovaszGraphLimits} for real-valued graphons.

\begin{theorem}[Theorem 4.26 in \cite{zucal2024probabilitygraphonsrightconvergence}]\label{ThmEquivQuotients} For any tight sequence $(\widetilde W_n)_n$ of probability graphons on $[0,1]$, the following are equivalent:
\begin{enumerate}
\item the sequence $(\widetilde W_n)_n$ is convergent in the unlabelled cut distance $\delta_{\square}$;

\item the quotient sets $\mathcal{Q}_k(\widetilde W_n)$ form a Cauchy sequence in the $d_{1,H}$ Hausdorff metric for every $k \geq 1$.
\end{enumerate}
\end{theorem}

We now show that convergence in the $P-$variables metric implies quotient convergence.

\begin{lemma}\label{LemmMetricImpliesQuotients}
Let $(W_n)$ be a tight sequence of $P-$variables $W_n\in \cB([0,1]\times [0,1]\times \Omega) $ and $(\widetilde W_n)$ be the sequence of probability graphons $\widetilde W_n$ associated with the $P-$variables $W_n.$   $P-$variables convergence of $(W_n)$ implies quotient convergence of the sequence of probability graphons $(\widetilde W_n)$. 
\end{lemma}
\proof
Let $U,W\in  \cB([0,1]\times [0,1]\times \Omega)$ and recall that $\widetilde U,\widetilde W$ denote the associated probability graphons. In the following, we will denote $\Omega_1=[0,1]$ and $\Omega_2=\Omega.$ Let's assume that $W$ and $U$ are sufficiently close in $d_M$, i.e.\ $d_M(U,W)\leq\delta$ and let's consider $R\in \mathcal{Q}_k(\widetilde W)$. We will now show that there exists an $R^\prime\in \mathcal{Q}_k(\widetilde U)$ such that $d_{1}(R,R^{\prime})\leq\varepsilon$. This is enough to prove the lemma as without loss of generality we can exchange the role of $M$ and $M^{\prime}$ in the proof.  

Let $\{f_1,\ldots,f_k\}$ be a function partition. By the definition of function partition, there exists a measurable partition $\cS=\{S_1,\ldots,S_k\}$ such that $f_1=\mathbbm{1}_{S_1}, \ldots, f_k=\mathbbm{1}_{S_k}.$ Let $R=\widetilde W/\cS,$ i.e.\ $R$ is the quotient of  $\widetilde W$ corresponding to the partition $\cS.$ We have by the definition of Hausdorff distance $d_H$ that for any $\delta^{\prime}>0$ there exists a function partition $w_1,\ldots,w_k$ such that 

$$
\begin{aligned}
&d_{\mathcal{LP}}(\mathcal{S}(f_1,\ldots,f_k,W),\mathcal{S}(w_1,\ldots,w_k,U))\leq d_H(S^{\prime}_k(W),S^{\prime}_k(U))+\delta^{\prime}.
\end{aligned}
$$
Moreover, by the definition of the $P-$variables metric $d_M,$ we have 
$$
\begin{aligned}
& d_H(S_k(W),S_k(U))\leq  2^{k+1}d_M(W,U) \leq 2^{k+1}\delta
\end{aligned}
$$
Therefore, for every $\delta^{\prime}>0$ there exists a $\delta^{\prime}>\delta >0$ small enough, such that 
$$
\begin{aligned}\label{ineq1}
&d_{\mathcal{LP}}(\mathcal{L}(f_i \otimes\mathbbm{1}_{\Omega_1} \otimes \mathbbm{1}_{\Omega_2},\mathbbm{1}_{\Omega_1}\otimes f_j \otimes\mathbbm{1}_{\Omega_2},W),\mathcal{L}(w_i \otimes\mathbbm{1}_{\Omega_1} \otimes \mathbbm{1}_{\Omega_2},\mathbbm{1}_{\Omega_1}\otimes w_j \otimes\mathbbm{1}_{\Omega_2},U))\leq\\
&d_{\mathcal{LP}}(\mathcal{S}(f_1,\ldots,f_k,W),\mathcal{S}(w_1,\ldots,w_k,U))\leq \\
&d_H(S^{\prime}_k(W),S^{\prime}_k(U))+\delta^{\prime} \leq \\
& d_H(S_k(W),S_k(U))+3\delta^{\prime}\leq  \\
& 2^{k+1}d_M(W,U)+3\delta^{\prime}\leq \\
&2^{k+1}\delta+3\delta^{\prime}\leq\\
& (2^{k+1}+3)\delta^{\prime},
\end{aligned}
$$
where the first inequality follows from Lemma \ref{LemmIneqProkh} and the third inequality follows from Lemma \ref{partlem} when $\delta^{\prime}>\delta >0$ is small enough.

Again, there exists a measurable partition $\cT=\{T_1,\ldots,T_k\}$ such that $w_1=\mathbbm{1}_{T_1},\ldots, w_k=\mathbbm{1}_{T_k}.$

In the following,  for a measure $\mu$ and for a measurable set $A\subset \R^p,$  we will denote with $\mu_A$ the restriction of the measure $\mu$ with respect to $A,$ i.e. the measure $\mu_A$ such that for any measurable set $Q\subset \R^p,$
\begin{equation}
    \mu_A(Q)=\mu(A\cap Q).
\end{equation}

Now, for $\delta^{\prime}>0$ small enough, by Lemma \ref{LemmCntinuityProhRestriction} (alternatively Portmateau theorem, Theorem \ref{PortThm}), for any $1/4>\varepsilon>0$ and  $i,j \in [k]$ we have (compare with \eqref{identityProbGraphPsets})
$$
\begin{aligned}
&d_{\mathcal{LP}}\left(\widetilde W(S_i\times S_j,\cdot),\widetilde U(T_i\times T_j,\cdot)\right)=\\
& d_{\mathcal{LP}}(\mathcal{L}(f_i \otimes\mathbbm{1}_{\Omega_1} \otimes \mathbbm{1}_{\Omega_2},\mathbbm{1}_{\Omega_1}\otimes f_j \otimes\mathbbm{1}_{\Omega_2},W )_{\{(1,1)\}\times \R},\mathcal{L}(w_i \otimes\mathbbm{1}_{\Omega_1} \otimes \mathbbm{1}_{\Omega_2},\mathbbm{1}_{\Omega_1}\otimes w_j \otimes\mathbbm{1}_{\Omega_2},U)_{\{(1,1)\}\times \R})= \\
&d_{\mathcal{LP}}(\mathcal{L}(f_i \otimes\mathbbm{1}_{\Omega_1} \otimes \mathbbm{1}_{\Omega_2},\mathbbm{1}_{\Omega_1}\otimes f_j \otimes\mathbbm{1}_{\Omega_2},W )_{E},\mathcal{L}(w_i \otimes\mathbbm{1}_{\Omega_1} \otimes \mathbbm{1}_{\Omega_2},\mathbbm{1}_{\Omega_1}\otimes w_j \otimes\mathbbm{1}_{\Omega_2},U)_{E})\leq \\
&\frac{\varepsilon}{2k^2}
\end{aligned}
$$

where $E=(3/4,5/4)\times (3/4,5/4)\times \R.$ Observe in fact that $$\mathcal{L}(f_i \otimes\mathbbm{1}_{\Omega_1} \otimes \mathbbm{1}_{\Omega_2},\mathbbm{1}_{\Omega_1}\otimes f_j \otimes\mathbbm{1}_{\Omega_2},W )(E)=\mathcal{L}(f_i \otimes\mathbbm{1}_{\Omega_1} \otimes \mathbbm{1}_{\Omega_2},\mathbbm{1}_{\Omega_1}\otimes f_j \otimes\mathbbm{1}_{\Omega_2},U)(\{(1,1)\}\times \R)$$
and 
$$\mathcal{L}(w_i \otimes\mathbbm{1}_{\Omega_1} \otimes \mathbbm{1}_{\Omega_2},\mathbbm{1}_{\Omega_1}\otimes w_j \otimes\mathbbm{1}_{\Omega_2},W )(E)=\mathcal{L}(w_i \otimes\mathbbm{1}_{\Omega_1} \otimes \mathbbm{1}_{\Omega_2},\mathbbm{1}_{\Omega_1}\otimes w_j \otimes\mathbbm{1}_{\Omega_2},U)(\{(1,1)\}\times \R).$$

Moreover, inequality \eqref{ineq1}, using Lemma \ref{LemmIneqProkh}, also implies that (choosing $\delta^{\prime}>0$ small enough)

$$\begin{aligned}
    &|\lambda(S_i)-\lambda(T_i)|=\\
    &d_{\mathcal{LP}}(\lambda(S_i)\delta_1+(1-\lambda(S_i))\delta_0,\lambda(T_i)\delta_1+(1-\lambda(T_i))\delta_0)=\\
  &d_{\mathcal{LP}}(\cL(\mathbbm{1}_{S_i}),\cL(\mathbbm{1}_{T_i}))\leq  \\
  &d_{\mathcal{LP}}(\mathcal{L}(f_i \otimes\mathbbm{1}_{\Omega_1} \otimes \mathbbm{1}_{\Omega_2},\mathbbm{1}_{\Omega_1}\otimes f_j \otimes\mathbbm{1}_{\Omega_2},W),\mathcal{L}(w_i \otimes\mathbbm{1}_{\Omega_1} \otimes \mathbbm{1}_{\Omega_2},\mathbbm{1}_{\Omega_1}\otimes w_j \otimes\mathbbm{1}_{\Omega_2},U))\leq\\
  &d_{\mathcal{LP}}(\mathcal{S}(f_1,\ldots,f_k,W),\mathcal{S}(w_1,\ldots,w_k,U))\leq \\
  &  \leq \frac{\varepsilon}{2k}.
\end{aligned}$$

Therefore, for $R^{\prime}=\widetilde U/\cT$ the quotient corresponding to the partition $\cT=\{T_1,\ldots, T_k\}$ associated with the  function partition $w_1,\ldots, w_k,$ we have $$d_1(R,R^{\prime})=\sum_{i\in[k]}|\lambda(S_i)-\lambda(T_i)|+\sum_{i,j\in [k]} d_{\mathcal{LP}}\left(\widetilde W(S_i\times S_j,\cdot),\widetilde U(T_i\times T_j,\cdot)\right)\leq \varepsilon. $$
\endproof
As a consequence, we obtain the following corollary.
\begin{corollary}
\label{CorMetricImpliesUnlbCutDist}
Let $(W_n)$ be a tight sequence of $P-$variables $W_n\in \cB([0,1]\times [0,1]\times \Omega) $ and $(\widetilde W_n)$ be the sequence of probability graphons $\widetilde W_n$ associated with the $P-$variables $W_n.$  $P-$variables convergence of $(W_n)$ implies convergence in the unlabelled cut distance of the sequence of probability graphons $(\widetilde W_n)$. 
\end{corollary}
\proof
The proof follows directly from Lemma \ref{LemmMetricImpliesQuotients} and Theorem \ref{ThmEquivQuotients} (Theorem 4.26 in \cite{zucal2024probabilitygraphonsrightconvergence}).
\endproof
We introduce here the $P-$variable step-representation of an $n\times n$ matrix $M=(M_{ij})_{i,j\in[n]},$ that is the $P-$variable $W_M\in \cB([0,1]\times [0,1]\times [0,1])$ defined as the step-function
\begin{equation}\label{Eq:P-setStepMatrix}
    W_M(x_1,x_2,x_{12})=M_{ij} 
\end{equation}
for $(x_1,x_2)\in [i-1/n,i/n] \times [j-1/n,j/n]$ and every $x_{12}\in [0,1].$

Let $H=(V(H),E(H),\beta(H))$ be a weighted graph where $V(H)=\{v_1,\ldots, v_n\}$ with adjacency matrix $A(H).$ We define $W_H\in \cB([0,1]\times [0,1]\times [0,1]),$ the $P-$variable step-representation of $H$ as  
\begin{equation*}
    W_H=W_{A(H)}
\end{equation*}
where $W_{A(H)}$ is the $P-$variable step-representation \eqref{Eq:P-setStepMatrix} of the adjacency matrix $A(H)$ of $H.$

We show now two lemmas which directly imply that convergence of weighted graphs with the number of vertices diverging (considered as $P-$sets, see Example \ref{ExGraphasPsets}) and probability graphons convergence are equivalent. The proof of lemmas \ref{InclusionskProfP-setsGraphs} and \ref{LemmApproxGraphPset} are reminiscent of the proof of Lemma 8.3 in \cite{backhausz2018action}.

\begin{lemma}\label{InclusionskProfP-setsGraphs}
For every (random) weighted graph $H=(V(H),E(H),\beta(H))$ we have the inclusions $\mathcal{S}_k(H)\subseteq \mathcal{S}_k( W_H)$ and $\mathcal{S}^{\prime}_k(H)\subseteq \mathcal{S}^{\prime}_k( W_H)$ considering $H$ as a $P-$variable.      
\end{lemma}
\proof
Assume that $V(H)=[n]$ without loss of generality. Recall that we consider $[n]$ as a probability space endowed with the discrete $\sigma-$algebra and with the uniform measure. 

We can observe that $\mathcal{S}_k(H)\subseteq \mathcal{S}_k( W_H),$ i.e.\ the $k-$profile of the $P-$variable representation of $H$ is contained in the $k-$profile of the $P-$variable step-representation $W_H$ of $H$. Indeed, if $\mu\in\mathcal{S}_k(H)$ then there exist functions/vectors $v_1,\ldots, v_k\in L_{[-1,1]}^{\infty}([n])$ such that $\mu=\mathcal{S}(v_1,\ldots,v_k,H).$ We observe that $\mu=\mathcal{S}(v_1,\ldots,v_k,H)=\mathcal{S}(f_1,\ldots,f_k,W_H)$ choosing the functions $f_i\in L_{[-1,1]}^{\infty}([0,1])$ defined for almost every $x\in [0,1]$ as $f_i(x)=v_i(\lceil xn\rceil).$ 

Similarly, we can observe that $\mathcal{S}^{\prime}_k(H)\subseteq \mathcal{S}^{\prime}_k( W_H).$
\endproof

\begin{lemma}\label{LemmApproxGraphPset} For every $\varepsilon>0,$ there exists a number $N\in \N$ such that if $H=(V(H),E(H),\beta(H))$ is a weighted graph with $|V(H)|\geq N,$ then considering $H$ as a $P-$variable (see Example \ref{ExGraphasPsets}) we have $d_M(H, W_H)\leq\varepsilon$.
\end{lemma}
\proof Assume that $V(H)=[n]$ for some $n\geq N.$ Recall that we consider $[n]$ as a probability space endowed with the discrete $\sigma-$algebra and with the uniform measure. 

From Lemma \ref{InclusionskProfP-setsGraphs} and Lemma \ref{LemmEquivForImplMetric} (Theorem \ref{partthm}), we need to show only that for $N\in \N$ big enough for $\nu\in\mathcal{S}^{\prime}_k(W_H)$ it exists $\nu^{\prime}\in \mathcal{S}^{\prime}_k(H)$ close to $\nu$ in Levy-Prokhorov distance. By the definition of $\mathcal{S}^{\prime}_k(W_H)$, there exists a partition function $f_1,\ldots,f_k\in L_{[-1,1]}^{\infty}([0,1])$ such that $\nu =\mathcal{S}(f_1,\ldots,f_k,W_H).$ 

Let's set $N=m t$ for $m$ and $t $ positive integers to be chosen later. Let's consider $\widetilde W_H$ the probability graphon associated with $W_H,$ that is constant on sets $[(i-1)/n,i/n]\times [(j-1)/n,j/n]$ for every $i,j\in [n].$
From the Weak Regularity Lemma, i.e.\ Lemma \ref{LemmaWeakRegularity} (Corollary 4.14 in \cite{abraham2023probabilitygraphons}), for each $\varepsilon>0$ there exists a positive integer $m$ and a measurable partition $\cP=\{P_1,\ldots, P_m\}$ such that $d_{\square}(\widetilde W_{H,\cP},\widetilde W_{H})<\varepsilon,$ where $\widetilde W_{H,\cP}$ is the probability graphon associated with $W_H.$ Moreover, as $\widetilde W_H$ is constant on sets $[(i-1)/n,i/n]\times [(j-1)/n,j/n]$ for every $i,j\in [n]$ we can assume without loss of generality that for every $i\in [n]$ we have $[(i-1)/n,i/n]\subset P_s$ for some $s\in [m].$

We therefore fixed the quantity $m$ in $N=mt.$

From Lemma \ref{LemmKprofBoundCutdist} and $d_{\square}(\widetilde W_{H,\cP},\widetilde W_{H})<\varepsilon,$ we obtain that $$
d_{\mathcal{LP}}(\cS(f_1,\ldots, f_k,W_{H,\cP}),\cS(f_1,\ldots, f_k,W_{H}))<\varepsilon^{\prime}.
$$
We now define the quantities $\mu_{s,i}\in [0,1]$ for $s\in [m],$ $i\in [k]$ and for $(f_1,\ldots, f_k),$ the $k-$dimensional random vector obtained from the function partition $f_1,\ldots,f_k,$ as  
\begin{equation*}
\begin{aligned}
    & \mu_{s,i}=\lambda(\{x\in P_s:\ f_1(x)=0,\ldots,f_i(x)=1,\ldots, f_k(x)=0\})=\lambda(\{x\in P_s:\ (f_1(x),\ldots, f_k)=e_i\})
\end{aligned}
   \end{equation*}
where $\lambda$ denotes the Lebesgue measure on $[0,1].$ We observe that $\sum_{i\in [k]s\in [m]}\mu_{s,i}=1.$

For each $i\in [k]$, we choose $\lfloor n\mu_{s, i}\rfloor$ intervals $((j-1) / n, j / n)$ contained in $P_s$, and define $(\bar{f}_1, \ldots, \bar{f}_k)$ to take as value $e_i$ (the $k-$dimensional vector with the $i-$th entry having value $1$ and every other entry having value $0$) on these intervals. We choose different intervals for different $i\in [k].$ This is possible since  the sum of $\mu_{s, i}$ over $s\in [m]$ and $i\in [k]$ is $1$, and we have $n$ intervals. On the rest of the intervals, which we will call reminder intervals, we let $\left(\bar{f}_1, \ldots, \bar{f}_k\right)$ be $e_1.$


Observe that as  $\widetilde{W}_{H,\cP}$ is constant on each $P_i \times P_j,$ if the quantities $\mu_{s, i}\in [0,1]$ are the same for every $s\in [m]$ and $i\in [k]$ for $\left(f_1, \ldots, f_k\right)$ and some $\left(u_1, \ldots, u_k\right)$, then
$$
\cS\left(f_1, \ldots, f_k, W_{H,\cP}\right)=\cS\left(u_1, \ldots, u_k, W_{H,\cP}\right) .
$$
However, the quantities $\mu_{s, i}$ for $(\bar{f}_1, \ldots, \bar{f}_k)$ are changed with respect to the quantities $\mu_{s, i}$ for $(f_1, \ldots, f_k)$ when we use the rounding to make the function constant on each interval $((j-1) / n, j / n)$. For every $s\in [m]$ and $i\in [k]$, we have, at most, one remainder interval, which has Lebesgue measure $1/n$ and therefore we are changing the function on a set of Lebesgue measure $km/n\leq km/N=km/mt=k/t.$ Thus, for large enough $N=mt$ (that is large enough $t$ as $m$ has been already fixed), we can find functions $u_1, \ldots, u_k:[0,1] \rightarrow \{0,1\}$ such that they have the same quantities $\mu_{s, i}$ as $f_1, \ldots, f_k$, and $\left\|u_j-\bar{f}_j\right\|_1 \leq \varepsilon_1=k/t$. By Lemma \ref{coupdist2} it follows that
$$
\begin{aligned}
&d_{\mathcal{LP}}\left(\cS(f_1, \ldots, f_k, W_{H,\cP}), \cS(\bar{f}_1, \ldots, \bar{f}_k,W_{H,\cP})\right) =\\
&d_{\mathcal{LP}}\left(\cS(u_1, \ldots, u_k, W_{H,\cP}), \cS(\bar{f}_1, \ldots, \bar{f}_k,W_{H,\cP})\right)\leq \varepsilon^{1/2}_1(k+1)^{3/4}<\varepsilon^{\prime}
\end{aligned}
$$
for $N$ large enough.

From Lemma \ref{LemmKprofBoundCutdist} and $d_{\square}(\widetilde W_{H,\cP},\widetilde W_{H})<\varepsilon,$ again we have $$
d_{\mathcal{LP}}(\cS(\bar{f}_1, \ldots, \bar{f}_k,W_{H,\cP}),\cS(\bar{f}_1, \ldots, \bar{f}_k,W_{H}))<\varepsilon^{\prime}.
$$
Since for $i\in [k]$ the function $\bar{f}_i$ is constant on each interval $((j-1)/ n,j / n)$, we can naturally associate the function $v_i:[n] \rightarrow \{0,1\}$ to it, which is defined for every $j\in [n]$ as $v_i(j)=\bar{f}_i(x)$ for $x\in((j-1)/ n,j / n) .$  Therefore, $v_1,\ldots, v_k$ is a function partition and 

$$\cS(v_1, \ldots, v_k,H)=\cS(\bar{f}_1, \ldots, \bar{f}_k,W_{H}).$$
Therefore, considering everything together, by the triangular inequality we obtain
$$
d_{\mathcal{LP}}(\cS(v_1, \ldots, v_k,H),\cS(f_1, \ldots, f_k,W_{H}))=d_{\mathcal{LP}}(\cS(\bar{f}_1, \ldots, \bar{f}_k,W_{H}),\cS(f_1, \ldots, f_k,W_{H}))<3\varepsilon^{\prime}.
$$

Lemmas \ref{InclusionskProfP-setsGraphs} and \ref{LemmApproxGraphPset}  directly imply the following theorem.

\begin{theorem} If $(H_n)_n$ is a growing sequence of finite weighted graphs, $W_{H_n}$ the $P-$variables step-representations of $H_n$ (recall \ref{Eq:P-setStepMatrix}) and $\widetilde W_{H_n}$ the associated probability graphons. Then the $P-$variables convergence of $(H_n)_n,$ interpreted as a $P-$variable, recall Example \ref{ExGraphasPsets}, is equivalent to the convergence in unlabelled cut metric of the sequence of probability graphons $(\widetilde W_{H_n})_n$. 
\end{theorem}

\section{Translating results and examples}\label{Sec5TranslatResultsExamp}

\subsection{Consequences of results for probability graphons}

In this section, we translate several results from the theory of probability graphons \cite{abraham2023probabilitygraphons,zucal2024probabilitygraphonsrightconvergence} to the theory of $P-$variables exploiting Theorem \ref{TH:equivalence_dM_cut_dist}.

First, we characterize which $P-$variables are identified by the $P-$variables metric. In particular, the following lemma is analogous to Theorem \ref{theo:Wm=W} (Theorem 3.17 in \cite{abraham2023probabilitygraphons}) for probability graphons.

\begin{lemma}
Two $P-$variables $U,W\in \cB([0,1]\times[0,1]\times \Omega)$ are identified in the $P-$variables metric $d_M,$ i.e.\ $d_M(U,W)=0,$ if and only if

$$
\cL(W(\phi(x_1),\phi(x_2),\cdot))=\cL(U(\psi(x_1),\psi(x_2),\cdot))
$$
for almost every $(x_1,x_2)\in [0,1]\times[0,1].$
\end{lemma}
\proof
The result follows directly from Theorem \ref{TH:equivalence_dM_cut_dist} and Theorem \ref{theo:Wm=W} (Theorem 3.17 in \cite{abraham2023probabilitygraphons}).
\endproof

We give now a characterization of compact sets of $P-$variables. In particular, the following theorem is analogous to Theorem \ref{ThmRelCompact} (Proposition 5.2 in \cite{abraham2023probabilitygraphons}) for probability graphons.

\begin{theorem}[Compactness]
A set $\mathcal{K}$ of $P-$variables on $\mathcal{K}\subset \cB([0,1]\times[0,1]\times \Omega)$ is relatively compact for $d_M$ if and only if $\mathcal{K}$ is 
tight. 
\end{theorem}
\proof
The result follows directly from Theorem \ref{TH:equivalence_dM_cut_dist}, Lemma \ref{LemmEqTightness} and Theorem \ref{ThmRelCompact} (Proposition 5.2 in \cite{abraham2023probabilitygraphons}).
\endproof

One can also sample random matrices from $P-$variables. This is exactly the famous construction from Aldous-Hoover theorem \cite{hoover1979relations,aldous1981representations,aldous2010exchangeability}.
We explain here briefly how to sample a weighted graph on the vertex set $[n]$, or equivalently a $n\times n$ square matrix with zero diagonal entries from a $P-$variable $W\in \cB(\Omega_1\times\Omega_1\times \Omega_2 ),$ where $(\Omega_1,\cA_1,\P)$ and $(\Omega_2, \cA_2, \mu)$ are probability spaces. First, we sample $X_1,\ldots, X_k$ in $\Omega_1$ independently with probability $\mathbb{P}.$ Additionally, we sample $Y_{ij} $ in $\Omega_2$ according to the distribution $\mu$ for $i,j\in [n]$. We can now consider the $n\times n$ random matrix $M^{(n)}=M^{(n)}(W)$ with non-diagonal entries

\begin{equation}\label{eq_RandomMatrixSamp}
    M^{(n)}_{ij}=W(X_i,X_j,Y_{ij})
\end{equation}
for $i,k\in [n]$ such that $i\neq j.$
Typically one considers $P-$variables $W\in\cB([0,1]\times [0,1]\times [0,1])$ with the Lebesgue measure. Then it is enough to sample $X_1,\ldots, X_n$ and $Y_{ij}$ for $i,j\in [n]$ i.i.d.\ uniformly at random on $[0,1]$ and define $M$ as in \eqref{eq_RandomMatrixSamp}.

\begin{remark}\label{RmkMsym}
The random graph $M^{(n)}=M^{(n)}(W)$ is in general not symmetric. There can be $M^{(n)}_{ij}\neq M^{(n)}_{ij}$ also when $W$ is symmetric (recall Definition \ref{DefSymPset}). However, one can also consider a symmetrized version of $M^{(n)}$ that we denote with $M^{(n)}_{\text{sym}}$ defined as follows $$(M^{(n)}_{\text{sym}})_{ij}=\begin{cases}
    (M^{(n)})_{ij}  \text{  if } i\leq j\\
    (M^{(n)})_{ji} \text{  if } i>j.
\end{cases}$$
\end{remark}

 Recall that we can interpret a matrix as a $P-$variable, see Example \ref{Examp:MatrixasPset}. We show that the matrices sampled from a $P-$variable (as explained above) do approximate the $P-$variable itself in the $P-$variables metric. This is the only result for which this work is not self-contained. We use results from \cite{abraham2023probabilitygraphons} that we do not summarise here for brevity. The following theorem is analogous to Theorem 6.13 in \cite{abraham2023probabilitygraphons}.

\begin{theorem}\label{ThmConVSubgraphs}
Let $U$ be a $P-$variable $U\in\cB([0,1]\times [0,1]\times \Omega).$ we have that $d_M(U,M^{(n)}(U))$ converge almost surely to $0.$ 
\end{theorem}
\proof
Let $\mathbb{G}(n,\widetilde U)$ and $\mathbb{H}(n,\widetilde U)$ be the probability graphons random graphs ($W-$random graphs) defined in $6.2$ (below Notation $6.4$) in \cite{abraham2023probabilitygraphons} where $\widetilde U$ is the probability graphon associated to the $P-$variable $U$. If we use the same random variables $X_1,\ldots,X_n$ to define $\mathbb{H}(n,\widetilde U)$ and the measure-valued matrix $\mathcal{L}(U(X_i,X_j,\cdot))$ then we obtain that $\mathbb{H}(n,\widetilde U)$ and $\mathcal{L}(U(X_i,X_j,\cdot))$ are identical. This means that $d_M(M^{(n)}(U), \mathbb{G}(n,\widetilde U))=0$ where with $\mathbb{G}(n,\widetilde U)$ here we mean the $P-$variable associated to $\mathbb{G}(n,\widetilde U).$ We observe that $$
d_M( M^{(n)}(U),U)\leq d_M(M^{(n)}(U), \mathbb{G}(n,\widetilde U))+d_M( \mathbb{G}(n,\widetilde U),U)=d_M( \mathbb{G}(n,\widetilde U),U).
$$

But by Theorem \ref{TH:equivalence_dM_cut_dist}, we obtain that $d_M( \mathbb{G}(n,\widetilde U),U)\rightarrow 0$ as $\delta_{\square}(\mathbb{G}(n,\widetilde U),\widetilde U)\rightarrow 0.$
Therefore, the result follows directly from Theorem 6.13 in \cite{abraham2023probabilitygraphons} and Theorem \ref{TH:equivalence_dM_cut_dist}.
\endproof

\begin{remark}
It is also easy to observe that symmetric $P-$variables can be approximated using $M^{(n)}_{\text{sym}}.$
\end{remark}

The previous theorem has the following important direct consequence.

\begin{corollary}
The space of (deterministic) weighted graphs is dense in the space of (weak isomorphism classes of) $P-$variables equipped with the $P-$variables metric $d_M.$
\end{corollary}

We introduce now useful invariants for probability graphons from \cite{abraham2023probabilitygraphons} and \cite{zucal2023action}, homomorphism densities and overlay functions, and we define their $P-$variables counterparts.\newline

For simplicity of notation, in the rest of this section, we will only consider $P-$variables $W$ such that $W\in \cB([0,1]\times [0,1]\times [0,1])$ and probability graphons on $[0,1],$ i.e.\ $\widetilde W\in \Graphon.$ However, it is straightforward to generalize the following quantities to $P-$variables and probability graphons defined on different probability spaces than $[0,1]$ endowed with the Lebesgue measure $\lambda.$

The first invariants we will be interested in are homomorphism densities.

\textbf{Homomorphism densities:}
Let $\CbFunct$ be the space of real-valued continuous and bounded functions defined on $\Space.$ A $\CbFunct-$graph is a triple $G^{\beta}=(G,\beta(G))=(V(G),E(G),\beta(G))$ where $G$ is a simple graph with vertex set $V=V(G),$ edge set $E(G)$ and $\beta(G)$ is a function
$$
\beta(G):V(G)\times V(G)\rightarrow \CbFunct $$
such that $\beta(G)_{v,w}\neq 0$ if and only if $\{v,w\}\in E(G).$

We define the \emph{homomorphism density} of a $\CbFunct$-graph $G^{\beta}$
in a probability graphons $\widetilde{W}\in\UGraphon$ as:
\begin{equation}\label{eq_def_hom_densProbGraph}
t(G^{\beta},\widetilde{W}) = \ \int_{[0,1]^{V(G)}}  \prod_{(i,j)\in E(G)} \widetilde{W}(x_i,x_j; \beta_{i,j})		\	\prod_{i\in V(G)} \drv x_i	 ,
\end{equation}

where we abbreviated with $\drv x$ the differential $\mathrm{d}\lambda( x ).$

We define the \emph{homomorphism density} of a $\CbFunct$-graph $G^{\beta}$
in a probability graphons $W\in\UGraphon$ as:
\begin{equation}\label{eq_def_hom_densPsets}
t(G^{\beta},W) = M_W^F(\beta) 
= \ \int_{[0,1]^{V(G)}}  \prod_{(i,j)\in E(G)} \int_{\Space}g_{i,j}(W(x_i,x_j,x_{ij}) )\drv x_{ij}		\	\prod_{i\in V(G)} \drv x_i	 .
\end{equation}

We observe the following equalities:
\begin{equation*}
\begin{aligned}
     &t(G^{\beta},W)= \ \int_{[0,1]^{V(G)}}  \prod_{(i,j)\in E(G)} \int_{\Space}\beta_{i,j}(W(x_i,x_j,x_{ij}) )\drv x_{ij}		\	\prod_{i\in V(G)} \drv x_i=	\\
& \ \int_{[0,1]^{V(G)}}  \prod_{(i,j)\in E(G)} \widetilde{W}(x_i,x_j; \beta_{i,j})		\	\prod_{i\in V(G)} \drv x_i=t(G^{\beta},\widetilde{W})
     \end{aligned}
   \end{equation*}
where we abbreviated with $\drv x$ the differential $\mathrm{d}\lambda( x ).$

One can also consider edge-decorated graphs decorated with real-valued graphons $(G,w),$ where $G=(V,E)$ is a finite simple graph and $w=(w_e:\ e\in E),$ the functional $$t(G,w)= \ \int_{[0,1]^{V(F)}}  \prod_{(i,j)\in E(F)} w_{ij}(x_i,x_j )		\	\prod_{i\in V(F)} \drv x_i,$$
see also (7.16) in \cite{LovaszGraphLimits}.

If for  $G^{\beta}=(V,E,\beta)$ we consider the simple graph given by $G=(V,E)$ and we define $w^{\beta}=\{\widetilde{W}[\beta_{e}]: \ e\in E\}$ we additionally have 
\begin{equation}
\begin{aligned}\label{Eq:HomDensEqual}
     &t(G^{\beta},W)= t(G^{\beta},\widetilde{W})=t(G,w^{\beta}).
     \end{aligned}
   \end{equation}
We can also focus on subfamilies of $\CbFunct-$graphs. We can consider subsets of the space $\CbFunct.$ 

\begin{definition}
    A sequence of $[0,1]$-valued functions $\F = (f_k )_{k\in\N}$ in $\CbFunct$, with $f_0=\un_{\mathbf{Z}}$ the constant function equal to one, is \emph{convergence determining}  if for every sequence of measures $(\mu_n)_{n}$  and  measure $\mu$ on $\Space$ such  that we have $\lim_{n\to +\infty} \mu_n(f_k)  = \mu(f_k)$ for all $k\in\N$, then $(\mu_n )_{n}$ weakly converges to $\mu$.

\end{definition}
On Polish spaces, there always exists a convergence determining sequence, see \cite[Corollary 2.2.6]{Bogachev} or the proof of Proposition 3.4.4 in \cite{Ethier} and Remark 2.3 in \cite{abraham2023probabilitygraphons}.

In particular, if $\Space$ is finite, one natural choice for $\cF$ are the indicator functions $\mathbbm{1}_{\{a\}}$ for every $a\in \Space.$ In this case equation \eqref{Eq:HomDensEqual} connects also with the definition of Homomorphism density given in Section 4.3 in \cite{HypergraphonsZhao}. 

Let $\cF$ be a convergence determining sequence. We will call an $\cF-$graph a $\CbFunct-$graph $G^{\beta}=(G,\beta(G))$ such that the decoration function $\beta(G)$ is $\cF-$valued. 

Other invariants we will be interested in are overlay functionals.

\textbf{Overlay functional:}
Let $[0,1]$ be the unit interval endowed with the Lebesgue measure $\lambda.$ For a probability measure $\alpha$ on $[k]$, we denote by $\Pi(\alpha)$ the set of partitions $\left\{S_1, \ldots, S_k\right\}$ of $[0,1]$ into $k$ measurable subsets with $\lambda\left(S_i\right)=\alpha_i$. For each probability graphon  $\widetilde W$ and vertex weighted $\CbFunct-$graph $G_{\alpha}^{\beta}$ on the vertex set $[k]$, we define the overlay functional
\begin{equation}  
\mathcal{C}(\widetilde W, G_{\alpha}^{\beta})=\sup _{\left(S_{ 1}, \ldots, S_k\right) \in \Pi(\alpha)} \sum_{i, j \in[k]}  \int_{S_i \times S_j}\int_{\Space} \beta_{i j}(G)(z)\widetilde W(x_1, x_2,\mathrm{d}z) \mathrm{d} x_1 \mathrm{d} x_2,
\end{equation}
where we abbreviated with $\mathrm{d} x$ the differential $\mathrm{d}\lambda( x ).$

For a $P-$variable $W$ we define the overlay functional as
\begin{equation}  
\mathcal{C}(W, G_{\alpha}^{\beta})=\sup _{\left(S_{ 1}, \ldots, S_k\right) \in \Pi(\alpha)} \sum_{i, j \in[k]}  \int_{S_i \times S_j}\int_{\Space} \beta_{i j}(G)\left(W(x_1, x_2,x_{12})\right) \mathrm{d} x_1 \mathrm{d} x_2 \mathrm{d}\mu(x_{12}),
\end{equation}

We have \begin{equation}\label{Eq:equalityOverlays}
    \mathcal{C}( W, G_{\alpha}^{\beta})=\mathcal{C}(\widetilde W, G_{\alpha}^{\beta}).
\end{equation}

We obtain that homomorphism densities, overlay functionals and quotient sets characterize $P-$variables convergence completely as stated in the following theorem. The following theorem is analogous to a combination of Theorem 7.11 in \cite{abraham2023probabilitygraphons} and Theorem 4.29 in \cite{zucal2024probabilitygraphonsrightconvergence} for probability graphons.

\begin{theorem}
Let $(W_n)$ be a sequence of $P-$variables $W_n\in \cB([0,1]\times [0,1]\times \Omega)$ and $W\in \cB([0,1]\times [0,1]\times \Omega)$ a $P-$variable. The following properties are equivalent:
\begin{enumerate}
\item $(W_n)$ converges to $W$ in the $P-$variables metric $d_M.$ 
\item The numerical sequence $t(G^{\beta},W_n)$ converges to $t(G^{\beta},W)$ for every $\CbFunct-$graph $G^{\beta}$.
\item The numerical sequence $t(G^{\beta},W_n)$ converges to $t(G^{\beta},W)$ for every $\cF$-graph $G^{\beta}$.
\item For all $k\geq 2,$ the sequence of sampled subgraphs $M^{(k)}(W_n)$ converges in distribution to $M^{(k)}(W).$
\item The overlay functional values $\mathcal{C}\left(W_n, G^{\beta}\right)$ converge to $\mathcal{C}\left(W, G^{\beta}\right)$ for every decorated $\CbFunct-$graph $G^{\beta}$;
\item The overlay functional values $\mathcal{C}\left(W_n, G^{\beta}\right)$ converge to $\mathcal{C}\left(W, G^{\beta}\right)$ for every decorated $\cF-$graph $H^{\beta}$;
\item  Let $\widetilde W_n$ and $\widetilde W$ be the probability graphons associated to $W_n$ and $W. $ The quotient sets (recall ) $\mathcal{Q}_k\left(\widetilde W_n\right)$ converge to $\mathcal{Q}_k\left(\widetilde W\right)$ in the $d_{\square}^{\text {Haus }}$ Hausdorff metric for every $k \geq 1$.
\end{enumerate}   
\end{theorem}

\proof 
This follows directly from Theorem 7.11 in \cite{abraham2023probabilitygraphons}, Theorem 4.29 in \cite{zucal2024probabilitygraphonsrightconvergence}, Theorem \ref{TH:equivalence_dM_cut_dist}, Theorem \ref{ThmConVSubgraphs} and \eqref{Eq:HomDensEqual} and \eqref{Eq:equalityOverlays}.
\endproof

\subsection{Examples}

Theorem \ref{ThmConVSubgraphs} provides us with a plethora of examples of $P-$variables and their limits in the $P-$variables metric $d_M.$  However, we write down explicitly some examples. 

We start with the simplest example of a random graph model, the Erd\"os-Renyi graph.
\begin{example}[Erd\"os-Renyi graph]\label{ERgraph}
Consider the vertex set $V=[n]$ and we connect every pair of the ${n \choose 2}$ possible pairs independently with probability $p$, i.e.\ following the law of independent identically distributed Bernoulli random variables with parameter $p$. This is the \emph{Erd\"os-Renyi random graph} denoted with $G(n,p)$.
\end{example}

We can consider the random graphs $G(n,p)$ as $P-$variables as in Remark \ref{remarkRandGraphRepPset} and Example \ref{ExGraphasPsets} and consider their sequence.

However, one can also consider a sequence of realizations $G_n$ sampled from $G(n,p)$ for $n\rightarrow+\infty.$ This is a sequence of deterministic simple graphs and for every $G_n$ we can consider the associated $P-$variable, see Example \ref{ExGraphasPsets}.  The $P-$variables $G_n\in \cB([n]\times [n]\times [0,1])$ or $W_{G_n}\in \cB([0,1]\times [0,1]\times [0,1])$ ( see \ref{Eq:P-setStepMatrix}), obtained in this way, are independent of the third variable. However, a $P-$variable representing the limit object of this sequence in the $P-$variables metric $d_M$ is $W=\mathbbm{1}_{[0,1]\times [0,1]\times [0,p]}\in \cB([0,1]\times [0,1]\times [0,1]).$ This can be deduced directly by the theory of probability graphons/ decorated graphons, see \cite{lovász2010limits, abraham2023probabilitygraphons,KUNSZENTIKOVACS2022109284}. However, one can also leverage Theorem \ref{ThmConVSubgraphs} noticing that $G(n,p)=M_{\text{sym}}^{(n)}(\mathbbm{1}_{[0,1]\times [0,1]\times [0,p]})$ where $M_{\text{sym}}^{(n)}(\mathbbm{1}_{[0,1]\times [0,1]\times [0,p]})$ is defined as in Remark \ref{RmkMsym}.

\begin{remark}
In general, one can consider other sequences of quasi-random graphs and obtain the same limit object $W=\mathbbm{1}_{[0,1]\times [0,1]\times [0,p]}\in \cB([0,1]\times [0,1]\times [0,1]).$
\end{remark}

\begin{remark}
We remark again that as a limit we could consider other $P-$variables, as an example we could have chosen any other function $\mathbbm{1}_{S}$ where $S=\{(x_1,x_2,x_{12})\in [0,1]^3: f(x_1,x_2)\leq x_{12}\leq g(x_1,x_2)\}$ where $f,g$ are two measurable functions from $[0,1]^2$ to $[0,1]$ such that $f(x_1,x_2)-g(x_1,x_2)=p$ for almost every $(x_1,x_2)\in [0,1]^2.$ These $P-$variables are all identified by the $d_M$ metric.
\end{remark}

The Erd\"os-Renyi graph is a \emph{local} random model, i.e.\ the presence of one edge does not affect the presence of another disjoint edge, differently from the following example.

\begin{example}[On/Off random graph]
We will call the random graph model on the vertex set $[n]$ that with probability $1\geq p\geq0$ is the complete graph and with probability $1-p$ the graph with no edges the \emph{On/Off random graph}  $OnOff(n,p)$ with parameter $p$.
\end{example}
We observe that in this case considering the $P-$variables associated with the random graph $OnOff(n,p)$ for $n\rightarrow +\infty$ and $0<p<1$ a limit $P-$variable can be again $\mathbbm{1}_{[0,1]\times [0,1]\times [0,p]}.$ However, if we consider a sequence $G_n$ of realizations of $OnOff(n,p)$ this is just a sequence in which for infinitely many $n$ the graphs $G_n$ are the complete graph on the vertex set $[n]$ and for the remaining infinitely many $n$ the graph $G_n$ is the graph on the vertex set $[n]$ with no edges. Therefore, the relative $P-$variables subsequences do converge to the two not weakly isomorphic $P-$variables $\mathbbm{1}_{[0,1]\times [0,1]\times [0,1]}$ and $0.$ This example illustrates that our notion of convergence does not capture random graphs that are non-local (or equivalently matrices with no independently distributed entries) in the correct way. This is the reason for the condition in Example \ref{Examp:MatrixasPset} and is a well-known phenomenon in the theory of real-valued graphons, see Example 11.10 in \cite{LovaszGraphLimits}. One can also compare with Aldous-Hoover theorem \cite{aldous1981representations,hoover1979relations,aldous2010exchangeability}, where measurable functions from $\Omega_0\times \Omega_1\times \Omega_1\times \Omega_2$ to $\R$ (where $\Omega_0$ is also a probability space) need to be considered to encode the non-local dependencies of an exchangeable array (random graph).

\begin{example}[Random graph with coloured edges at random]
Let $\mathbf{p}=(p_0,p_1,\ldots,p_m)$ a probability vector, i.e.\ $p_i\geq 0$ for every $i\in \{0\}\cup [m]$ and $\sum^{m}_{i=0}p_i=1.$ We consider the random weighted graph model $G(n,m,\mathbf{p})$ on the vertex set $[n]$ in which independently for every pair $\{i,j\}\in [n]\times [n]$ with probability $p_0$ we have $\{i,j\}$ is not an edge and with probability $p_s$ we have $\{i,j\}$ is an edge with colour (weight) $s$ for $s\in [m].$
\end{example}

Using similar considerations to the case of Erd\"os-Renyi random graphs, one can see that a limit $P-$variable for a sequence of  graphs sampled from $G(n,m,\mathbf{p})$ is $\sum^m_{i=0}i\mathbbm{1}_{[0,1]\times [0,1]\times [\sum^{i-1}_{r=0}p_r,\sum^{i}_{r=0}p_r]}.$

Similarly, we can consider  random matrices and their limits. 

\begin{example}[Random matrix with $+1$ and $-1$ entries]
We consider the $n\times n$ random matrix $MU(n,p)$ with every entry chosen i.i.d.\ at random choosing with probability $0\leq p\leq1$ the value  $-1$ and with probability $1-p$ the value $1.$ 
\end{example}

The limit $P-$variable for a sequence of realizations $G_n$ sampled from $MU(n,p)$ for $n\rightarrow+\infty$ can be represented as $-1\mathbbm{1}_{[0,1]^2\times [0,p]}+1\mathbbm{1}_{[0,1]^2\times [p,1]}.$

We can give a last example in which the entries of the sequence of random variables are not all contained in a compact subset of $\R$ but the tightness condition \eqref{tightCondit} is still satisfied.
\begin{example}[Random matrices with i.i.d.\ normally distributed entries]
We consider the $n\times n$ random matrix $MN(n,p)$ with i.i.d. standard normal distributed entries.  
\end{example}

The limit $P-$variable for a sequence of realizations $G_n$ sampled from $MU(n,p)$ for $n\rightarrow+\infty$ can be represented with the $P-$variable $W\in \cB([0,1]\times[0,1]\times [0,1])$ such that $W(x_1,x_2,x_{12})=Probit(x_{12}),$ for each $x_1,x_2,x_{12}\in [0,1],$ where $Probit$ is the probit function, i.e.\ the inverse of the cumulative distribution function of the standard normal distribution. 

For more examples, one can see Section 2.6 in \cite{lovász2010limits} or Section 2.3 in\cite{KUNSZENTIKOVACS2022109284} and translate the probability graphons examples in examples of $P-$variables.

\section{Implies real-valued graphon convergence}\label{Sec6ImpliesRealValuedGraphons}

Let $\Omega$ and $\widetilde \Omega$ be two probability spaces. In this section, for $f,g\in L^{1}(\Omega)$ and a $P-$variable $V\in \cB(\Omega\times \Omega\times \widetilde{\Omega}),$ we will denote
$$
(f,g)_{V}=\mathbb{E}[V(f\otimes \mathbbm{1}_{\Omega}\otimes \mathbbm{1}_{\widetilde{\Omega}}) (\mathbbm{1}_{\Omega}\otimes g \otimes \mathbbm{1}_{\widetilde{\Omega}})]=\mathbb{E}[V (g\otimes \mathbbm{1}_{\Omega}\otimes \mathbbm{1}_{\widetilde{\Omega}}) (\mathbbm{1}_{\Omega}\otimes f \otimes \mathbbm{1}_{\widetilde{\Omega}})],
$$
when well defined.

Let now $f\in L^{\infty}(\Omega)\subset L^{1}(\Omega)$ and $g\in L^{q}(\Omega)\subset L^{1}(\Omega)$ and the $P-$variable $V$ such that $V\in L^p(\Omega\times \Omega\times \widetilde\Omega)$ where $p$ and $q$ are H\"older conjugate, i.e.\ $1/p+1/q=1.$ Therefore, from H\"older inequality we have the bound
\begin{equation}\label{EqBoundHoldPsetExp}
(f,g)_V \leq \|f\|_\infty\|g\|_q \|V\|_{p}. 
\end{equation}
Recall in fact that for $f\in L^p(\Omega)$ we have $\|f\|_p=\|f \otimes\mathbbm{1}_{\Omega} \otimes \mathbbm{1}_{\widetilde\Omega}\|_p=\|\mathbbm{1}_{\Omega}\otimes f\otimes \mathbbm{1}_{\widetilde\Omega}\|_p.$
\begin{remark}
From H\"older inequality one could obtain a more general inequality but the bound \eqref{EqBoundHoldPsetExp} will be enough for our purposes. 
\end{remark}

In this section, we will consider a generalization of the notion of function partition.

\begin{definition}[Fractional function partition]
A fractional function partition is a set $\{f_i\}_{i=1}^k$ of functions in $L^\infty_{[0,1]}(\Omega)$ such that $\sum_{i=1}^k v_i=1_\Omega$. We say that $\{v_i\}_{i=1}^k$ is balanced if $\|v_i\|_1=1/k$ holds for every $i\in [k]$.
\end{definition}
In particular, differently from a function partition, in a fractional function partition the functions are not required to be $\{0,1\}-$valued but only $[0,1]-$valued.

\begin{definition}Let $k\in\mathbb{N}$ and $W\in\mathcal{B}(\Omega\times \Omega \times \widetilde \Omega)$. A balanced fractional $k\times k$ averaged quotient of $W$ is a matrix $M\in\mathbb{R}^{k\times k}$ such that there is a balanced fractional function partition $\{f_i\}_{i=1}^k=\{f_1,\ldots,f_k\}$ of  $\Omega$ with $$\begin{aligned}
    &M_{i,j}=(f_i,f_j)_W=(f_j,f_i)_W=\\
    &\mathbb{E}[W (f_i\otimes \mathbbm{1}_{\Omega}\otimes \mathbbm{1}_{\Omega_1}) (\mathbbm{1}_{\Omega}\otimes f_j \otimes \mathbbm{1}_{\Omega_1})]=\\
    &\mathbb{E}[W (f_j\otimes \mathbbm{1}_{\Omega}\otimes \mathbbm{1}_{\Omega_1}) (\mathbbm{1}_{\Omega}\otimes f_i \otimes \mathbbm{1}_{\Omega_1})]
    \end{aligned}$$
for every $i,j\in [k]$. Let $\mathcal{AQ}_{k}(W)$ denote the set of all balanced fractional $k\times k$ averaged quotients of $W$.
\end{definition}
Note that by linearity, the entry sum of any matrix $M\in\mathcal{AQ}_k(W)$ is equal to $\E[W]$ for every $k\in\mathbb{N}$.
For two $k\times k$ square matrices $A,B\in \R^{k \times k}$ we define $d_{(1,m)}(A,B):=\sum_{i,j}|A_{i,j}-B_{i,j}|$ the entry-wise $l_1$ distance. For two subsets $S_1,S_2\subseteq \mathbb{R}^{k\times k}$ let $d_{{(1,m)},H}$ denote the corresponding Hausdorff distance. 

We present here, in the language of $P-$variables, a notion of convergence considered in \cite{KUNSZENTIKOVACS20191,backhausz2018action,Lpgraphon2,local-global1,borgs2011convergentAnnals}. 

\begin{definition}[Averaged quotient convergence] A sequence of $P-$variables $\{A_i\}_{i=1}^\infty$ is averaged quotient convergent if for every $k$ we have that $\mathcal{AQ}_k(A_i)$ is convergent in $d_{(1,m),H}$. 
\end{definition}
\begin{remark}
This is quotient convergence for real-valued graphons, in this work, we renamed it averaged quotient convergence to differentiate it from the quotient convergence for probability graphons \cite{zucal2024probabilitygraphonsrightconvergence}.    
\end{remark}

The proof of Lemma \ref{LemmadMAverageQuotconv} is similar in spirit to Lemma 7.1 in \cite{backhausz2018action}.
\begin{lemma}\label{LemmadMAverageQuotconv}
For every $k\in \N$ and $\varepsilon >0$ there exists $\delta$ such that for any two $P-$variables $U$ and $W$ with $\|U\|_p,\|W\|_p\leq C\leq \infty $ for $\infty \geq p>1$ and $d_M(U,W)\leq\delta$ we have that $d_{1,H}(\mathcal{AQ}_k(U),\mathcal{AQ}_k(W))\leq\varepsilon.$
\end{lemma}
\proof
Let's assume that $U\in\mathcal{B}(\Omega_1\times \Omega_1 \times \widetilde \Omega_1)$ and $W\in\mathcal{B}(\Omega_2\times \Omega_2 \times \widetilde \Omega_2)$ are sufficiently close in $d_M$, i.e.\ $d_M(U,W)\leq\delta$ and let's consider $M\in \mathcal{AQ}_k(U)$. We will now show that there exists an $M^\prime\in \mathcal{AQ}_k(W)$ such that $d_1(M,M')\leq\varepsilon$. This proves the lemma because without loss of generality we can exchange the role of $M$ and $M^{\prime}$ in the proof.  

Let $f_1,f_2,\dots,f_k$ be a balanced fractional function partition of $\Omega_1$ such that the corresponding balanced fractional $k\times k$ averaged quotient of $U$ is $M$. We have by the definition of $d_M$ that there are functions $w_1,\ldots,w_k$ in $L^\infty_{[-1,1]}(\Omega_2)$ such that $$d_{\mathcal{LP}}(\cS(f_1,\dots,f_k,U),\cS(w_1,\ldots, w_k,W)))\leq 2^{k+1}d_M(U,W).$$

Using Lemma \ref{LemmIneqProkh} and a modification of Lemma \ref{LemmApproxFunctPartMeas} for balanced fractional partitions, it is easy to see that depending on $k$ and an arbitrary constant $\varepsilon_2>0$ if $d_M(U,W)$ is small enough then there is a balanced fractional function partition $w_1^{\prime},\ldots,w_k^{\prime}$ on $\Omega_2$ such that $\|w_i-w'_i\|_{q}\leq\varepsilon_2$ holds for every $i\in [k].$ For such functions, we have for every $i,j\in [k]$ that 
$$\begin{aligned}&|(w_i,w_j)_W-(w_i',w'_j)_W|\leq |(w_i,w_j)_W-(w_i',w_j)_W|+|(w_i',w_j)_W-(w_i',w_j')_W|=\\
&|(w_i-w_i',w_j)_W|+|(w_i',w_j-w_j')_W|\leq 2\varepsilon_2C,\end{aligned}$$
where the last inequality is by \eqref{EqBoundHoldPsetExp}. Let $M^{\prime }\in\mathcal{AQ}_k(W)$ be a balanced fractional $k\times k$ averaged quotient of $W$ corresponding to $w_1^{\prime},\ldots,w_k^{\prime}$ and let $M^{\prime \prime}$ be the $k\times k$ matrix defined by $M^{\prime \prime}_{i,j}=(w_i,w_j)_W$. Therefore, by the definition of the metric $d_{(1,m)}$ we obtain that $d_{(1,m)}(M^{\prime},M^{\prime \prime})\leq 2\varepsilon_2k^2C\leq\varepsilon/2$ for $\varepsilon_2>0$ small enough. Moreover, from Lemma \ref{closedlem2}, we obtain that $|(f_i,f_j)_U-(w_i,w_j)_W|$ is small for every $i,j\in [k]$ if $d_M(U,W)$ is small enough. Hence, we also have $d_{(1,m)}(M,M^{\prime})\leq\varepsilon/2$ and thus, by the triangular inequality, we obtain also $d_{(1,m)}(M,M^{\prime\prime})\leq\varepsilon$ and this concludes the proof.
\endproof

The next lemma follows directly from the previous lemma.

\begin{lemma}\label{quprop} Let $(W_n)_{n}$ be a convergent sequence of  $P-$variables in the $P-$variables metric $d_M$ and such that $\|W_n\|_{p}\leq C$ for every $n\in \N$. Then $(W_n)_{n}$ is averaged quotient convergent. 
\end{lemma}

We remark that the assumption of a uniform bound on the $p-$moments of a sequence of $P-$variables $W_n$ is restrictive when considering sparse graphs. We illustrate this by adapting the argument from Page 3023 in \cite{LpGraphons1} to $P-$variables. 
Let $G$ be a simple graph  with vertex set $V(G)=[m],$ edge set $E(G)$ and let's consider it as a $P-$variable, recall Exercise \ref{ExGraphasGraphset}. Let $A(G)$ be its adjacency matrix.
We have that 
$$
\|G\|_{1}=\|A(G)\|_{1}=\frac{1}{m^2}\sum^m_{i,j=1}|A(G)_{i,j}|,
$$
corresponds to the density of $G$, i.e.\ the quantity $E(G)/V(G)^2=E(G)/m^2$ which is a quantity between $0$ and $1.$ One would like to define for a sequence of graphs $G_n$ the associated normalized $P-$variables sequence 
\begin{equation}\label{Eq:PsetFromNormAdjac}
W_n=\frac{G_n}{\|G_n\|_{1}},
\end{equation}
where $G_n$ is considered as a $P-$variable.
In such a way, the sequence of $P-$variables $W_n$ has at least uniformly bounded $1-$moments but this is not enough, see Example \ref{CounterexamoUnifL1bound}, we want $\|W_n\|_p\leq C$ for $1<p\leq \infty.$
However, $\|W_n\|_p\leq C$ for every $n$ implies $\|G_n\|_{1}\geq c>0,$ i.e.\ the graph has to be dense. Observe in fact that for simple graphs:

$$
C\geq\|W_n\|_p=\left\|\frac{G_n}{\|G_n\|_{1}}\right\|_p=\left\|\frac{A(G_n)}{\|G_n\|_{1}}\right\|_p=\frac{1}{\|G_n\|_{1}}\left\|A(G_n)\right\|^{\frac{1}{p}}_1=\|G_n\|^{\frac{1}{p}-1}_{1}.
$$
We remark that the third equality holds only because we consider the adjacency matrix of a simple graph.
Therefore, from this inequality, the density $E(G_n)/V(G_n)^2=\|G_n\|_{1}$ is uniformly upper-bounded by a number $c>0.$ This implies that the sequence of simple graphs $G_n$ has to be dense. 

Recall the definition of a contraction of a $P-$variable \eqref{Eq:DefContraction}. We have the following result that relates averaged quotient convergence of a sequence of $P-$variables and convergence in real-valued cut distance their contractions.

\begin{lemma}\label{ThmAvQuotCutMet}
A sequence of $P-$variables $(W_n)_n$ such that $W_n\in \cB([0,1]\times [0,1]\times \Omega)$ is averaged quotient convergent if and only if the sequence $(w_n)_n$ of the contractions $w_n$ of the $P-$variables $W_n$ is convergent in the real-valued cut-norm $\delta_{\square,\R}$.
\end{lemma}
\proof
The result follows directly from \cite{Lpgraphon2} (actually just a slight modification of Theorem 12.12 in \cite{LovaszGraphLimits} is already enough).
\endproof
Combining Lemma \ref{quprop} and Lemma \ref{ThmAvQuotCutMet} we directly obtain the following.
\begin{corollary}\label{CorollAverQuot}
Let $(W_n)_n$ be a sequence of $P-$variables such that $W_n\in \cB([0,1]\times [0,1]\times \Omega)$ and $\|W_n\|_{p}\leq C$ for every $n\in \N.$ If $(W_n)_n$ is convergent in the $P-$variables metric $d_M$ then the sequence $(w_n)_n$ of the contractions $w_n$ of the  $W_n$ is convergent in the real-valued cut-metric $\delta_{\square,\R}.$
\end{corollary}
The following example shows that there are different $P-$variables (in the $P-$variables metric $d_M$) that have the same contractions. 
\begin{example}
Let's consider the two $P-$variables $U,W\in \cB([0,1]\times [0,1]\times [0,1])$ such that for $U$ we have $U(x_1,x_2,x_{12})=1/2$ for almost every $x_1,x_2,x_{12}\in [0,1]$ and for $W=\mathbbm{1}_{[0,1]\times[0,1]\times[0,1/2]}.$ Moreover, let's denote with $u$ the contraction of $U$ and $w$ the contraction of $W.$ It is easy to observe that $d_M(U,W)\neq 0$ but $u=w=1/2$ almost everywhere on $[0,1]^2$ (and therefore $\delta_{\square,\R}(u,w)=0$).
\end{example}

The previous example shows that the converse of Lemma \ref{quprop} and Corollary \ref{CorollAverQuot} are false.

We now also give a counterexample to Corollary \ref{CorollAverQuot} and Lemma \ref{quprop} when one has only a uniform bound on the $1-$moments of a sequence of $P-$variables.

\begin{example}\label{CounterexamoUnifL1bound}
 Let's consider the sequence of (sparse but not uniformly bounded degree) Erd\"os-Renyi random graphs $G(n,p_n)$ where $p_n\rightarrow 0$ but $np_n \rightarrow +\infty$ for $n$ going to infinity. For the simple deterministic graphs $G_n$ obtained as realizations sampling from $G(n,p_n)$ we consider the associated $P-$variable that we still denote with $G_n.$ In particular, we are interested in the $P-$variables $G_n/p_n.$ We remark that the sequence $G_n/p_n$ has a uniform bound on the $p-$moments only for $p=1.$ We observe that the sequence $G_n/p_n$ converges in $P-$variables metric $d_M$ to the null $P-$variable $0.$ However, the sequence of contractions of $G_n/p_n$ converges in real-valued cut distance $\delta_{\square,\R}$ to $1,$ this follows from the theory of $L^p-$graphons, see Section 3.3.1 in\cite{Lpgraphon2}. Therefore, the statement of Lemma \ref{quprop} is false if a sequence of $P-$variables $W_n$ has no uniform bound on the $p-$moments for $p>1.$
\end{example}

We expect uniform integrability to be a sufficient (minimal) condition for Lemma \ref{quprop} to still hold, differently from the uniform bound on the $1-$moments, see for comparison also \cite{LpGraphons1,Lpgraphon2}. However, we do not explore this here for brevity.

We can also translate Corollary \ref{CorollAverQuot} in the language of probability graphons obtaining a new result for the theory of probability graphons.
\begin{corollary}
Let $(\widetilde W_n)_n$ be a sequence of probability graphons such that for some $1< p\leq \infty$ we have $\int_{[0,1]^2}\int_{\Space}|z|^p\widetilde W_n(x_1,x_2,\mathrm{d}z)\mathrm{d}x_1 \mathrm{d}x_2\leq C$ for all $n\in \N.$  If $(\widetilde W_n)_n$ is convergent in the unlabelled cut distance $\delta_{\square}$ then the sequence $(w_n)_n$ of the contractions $w_n$ of $\widetilde W_n$ is convergent in the real-valued cut-metric $\delta_{\square,\R}.$   
\end{corollary}
\proof
The result follows directly from Corollary \ref{CorollAverQuot} and Theorem \ref{TH:equivalence_dM_cut_dist} and \eqref{Eq:IdentityPushForwardProbGraphPset}.
\endproof

This shows the potential of the right convergence point of view for probability graphons developed in \cite{zucal2024probabilitygraphonsrightconvergence} and of $P-$variables convergence developed in this work for obtaining new results for probability graphons.

The results in this section also connect $P-$variable convergence with action convergence developed in \cite{backhausz2018action}.   

\section{Other properties of P-variables and probability graphons convergence}\label{Sec7OtherProp}

In this section, we show several properties of $P-$variables convergence that can be translated to probability graphons obtaining new results also for this theory.

We start showing that taking values in a closed subset of $\R$ (for example a finite set or the nonnegative reals $\R_+$) is closed under $P-$variables convergence.

\begin{lemma}\label{LemmPsetsClosedSetsClosed}
Let $\mathbf{C}$ be a closed subset of $\R.$ Let $(W_n)_n$ be a sequence of $\mathbf{C}-$valued $P-$variables and $W$ a $P-$variable. If  $(W_n)_n$ is $P-$variables convergent to $W$ then $W$ is also $\mathbf{C}-$valued.
\end{lemma}
\proof
From the definition of $P-$variables convergence follows directly that the measures $\mathcal{L}(W_n)$ converge weakly to $\mathcal{L}(W).$ As $W_n$ are $\mathbf{C}-$valued we have $\mathcal{L}(W_n)(\mathbf{C})=1.$  From Portmanteau theorem (Theorem \ref{PortThm}) we have $1=\limsup_n\mathcal{L}(W_n)(\mathbf{C})\leq \mathcal{L}(W)(\mathbf{C})\leq 1.$ From $\mathcal{L}(W)= 1$ follows directly that $W$ is $\mathbf{C}-$valued (almost everywhere).
\endproof

A direct consequence of the previous lemma is the following.
\begin{corollary}
Let's consider a closed set $\mathbf{C}\subset \R$ and a sequence of probability graphons $(\widetilde W_n)_n$ such that for every positive integer $n$ and every $x_1,x_2\in [0,1]$ the support of the probability measure $\widetilde W_n(x_1,x_2,\cdot)$ is contained in $\mathbf{C}.$  We have that if $(\widetilde W_n)_n$ converges in unlabelled cut metric to a probability graphon $\widetilde W$ then the probability measures $\widetilde W_n(x_1,x_2,\cdot)$ are all supported in $\mathbf{C}.$
\end{corollary}

We show also that uniform bounds on the $L^p-$norm for $p>1$ are preserved under $P-$variables convergence.

\begin{lemma}
Let $C>0$ and $1<p\leq \infty.$ Let $(W_n)_n$ be a sequence of $P-$variables such that $\|W_n\|_p\leq C$ for every positive integer $n$ and $W$ a $P-$variable. If  $(W_n)_n$ is $P-$variables convergent to $W$ then $\|W_n\|_p\leq C$ for $1<p\leq \infty.$
\end{lemma}
\proof
From the definition of $P-$variables convergence follows directly that the measures $\mathcal{L}(W_n)$ converge weakly to $\mathcal{L}(W).$ The result follows now directly from the properties of weak convergence.
\endproof

\begin{remark}
The previous result can be easily translated to the case of probability graphons.
\end{remark}

We show also that if a sequence of $P-$variables is defined on an atomless probability space (the unit interval $[0,1]$ with the Lebesgue measure for example) also the limit $P-$variable is defined on an atomless probability spaces.

\begin{lemma}\label{LemmPsetsClosSym}
Let $\Omega_{1,n}$ and $\Omega_{2,n}$ probability spaces and in addition assume that $\Omega_{1,n}$ is atomless for every positive integer $n$. Let $(W_n)_n$ be a sequence of $P-$variables such that $W_n\in \cB(\Omega_{1,n}\times\Omega_{1,n} \times \Omega_{2,n} )$ and $W$ a $P-$variable. If $(W_n)_n$ is $P-$variables convergent to $W$ then $W\in \cB(\Omega_1\times\Omega_1\times\Omega_2)$ for some atomless probability space $\Omega_1$ and any probability space $\Omega_2.$
\end{lemma}
\proof
Let $\Omega$ be an atomless $P-$variable. Let $W$ be a $P-$variable $W\in \mathcal{B}(\Omega\times\Omega\times \widetilde \Omega)$ and let $U$ be another $P-$variable $U \in \mathcal{B}\left(\Omega_1\times \Omega_1 \times  \Omega_2\right)$ such that $d_M(U, W)=d.$ As $\Omega$ is atomless there exists a function $f \in L_{[-1,1]}^{\infty}(\Omega)$ such that the distribution $\cL(f)$ is uniform on $[-1,1],$ i.e.\ $\cL(f)=\text{Unif}_{[-1,1]}.$ Let $\mu=\mathcal{S}(f,W).$ From $d_H\left(S_1(W), S_1(U)\right) \leq 2d$ we have that there exists $g\in L_{[-1,1]}^{\infty}(\Omega_1) $ such that $\nu=\mathcal{S}(g,W) \in \mathcal{S}_1(W)$ with $d_{\mathcal{LP}}(\mu, \nu) \leq 3 d.$ Therefore, for $\mu_1=\mathcal{L}(f)=\text{Unif}_{[-1,1]}$ and $\nu_1=\mathcal{L}(g),$the marginals of $\mu$ and $\nu$ on the first coordinate, we have $d_{\mathcal{LP}}\left(\mu_1, \nu_1\right) \leq 3 d.$ Thus we obtained that $\nu_1$ is at Levy-Prokhorov distance $3 d$ from the uniform distribution. Therefore, the largest atom of $\nu_1$ is at most $10d$ as by the definition of Levy-Prokhorov distance
$$
\inf\{\delta:\ \nu_1(\{x_0\})\leq \mu_1(B_{\delta}(x_0))+\delta\}\leq d_{\mathcal{LP}}(\mu_1,\nu_1)\leq 3d
$$
and $\mu_1(B_{\delta}(x_0))=2\delta$. 

We obtained that if $W$ is the limit of a sequence of $P-$variables $W_n\in \cB(\Omega_{1,n}\times\Omega_{1,n} \times \Omega_{2,n} )$ where $\Omega_{1,n}$ are atomless, then $W$ is atomless.
\endproof

Symmetry is another property preserved under $P-$variables convergence and probability graphons convergence.

\begin{lemma}
Let $(W_n)_n$ be a sequence of symmetric $P-$variables and $W$ a $P-$variable. If  $(W_n)_n$ is $P-$variables convergent to $W$ then $W$ is also symmetric.
\end{lemma}
\proof
Let $W_n\in \cB(\Omega_{1,n}\times \Omega_{1,n}\times \Omega_{2,n})$ and $W_n\in \cB(\Omega_{1}\times \Omega_{1}\times \Omega_{2})$ and we assume that 
For $\mu_n\in \cS_k(W_n)$ we have that there exists a $\mu\in S_k(W)$ such that $\mu_n$ weakly converges to $\mu.$
Moreover, we know there exist functions $f_{1,n},\ldots, f_{k,n}\in L^{\infty }_{[-1,1]}(\Omega_{1,n})$ and $f_{1},\ldots, f_{k}\in L^{\infty }_{[-1,1]}(\Omega_{1})$ such that $\mu_n=\cS(f_{1,n},\ldots, f_{k,n},W)$ and $\mu=\cS(f_{1},\ldots, f_{k},W).$

As $W_n$ is symmetric we have
$$
\begin{aligned}&\mu_n=\cS(f_{1,n},\ldots, f_{k,n},W)=\\
&\mathcal{L}(f_{1,n} \otimes\mathbbm{1}_{\Omega_1} \otimes \mathbbm{1}_{\Omega_2},\mathbbm{1}_{\Omega_1}\otimes f_{1,n} \otimes\mathbbm{1}_{\Omega_2},\ldots, f_{k,n} \otimes\mathbbm{1}_{\Omega_1} \otimes \mathbbm{1}_{\Omega_2},\mathbbm{1}_{\Omega_1}\otimes f_{k,n} \otimes\mathbbm{1}_{\Omega_2},W)=\\
&\mathcal{L}(\mathbbm{1}_{\Omega_1}\otimes f_{1,n} \otimes\mathbbm{1}_{\Omega_2},f_{1,n} \otimes\mathbbm{1}_{\Omega_1} \otimes \mathbbm{1}_{\Omega_2},\ldots, \mathbbm{1}_{\Omega_1}\otimes f_{k,n} \otimes\mathbbm{1}_{\Omega_2},f_{k,n} \otimes\mathbbm{1}_{\Omega_1} \otimes \mathbbm{1}_{\Omega_2},W).
\end{aligned}$$

The result follows now as $\mu=\cS_k(f_{1},\ldots, f_{k},W)$ is the weak limit of these measures.
\endproof

As a direct consequence, we get the following result.

\begin{corollary}
The set of symmetric probability graphons is closed in the topology of the unlabelled cut metric.  
\end{corollary}

Another direct corollary of the previous results is the following corollary.

\begin{corollary}
Graph-sets are closed in the space of $P-$variables.
\end{corollary}
\proof
This follows directly from Lemma \ref{LemmPsetsClosedSetsClosed} and Lemma \ref{LemmPsetsClosSym}.
\endproof

\section{Generalizations of P-variables convergence}\label{Sec8GeneralizationsAndDirect}

We consider in this section generalizations and alternative versions of $P-$variables and $P-$variables convergence showing the flexibility and the generality of this mathematical framework. 
\subsection{Multiple P-variables at the same time:}

One can consider the convergence of multiple $P-$variables at the same time. In particular, let $\Omega_1$ and $\Omega_2$ be two probability spaces and for a positive integer $m$ let's consider the vector $(W_1,\ldots,W_m)$ obtained by the $P-$variables $W_1,\ldots,W_m\in \cB(\Omega_1\times \Omega_1\times \Omega_2).$ 

We can now consider for functions $f_1,\ldots , f_k \in L_{[-1,1]}^{\infty}(\Omega_1)$ the quantities
$$
\begin{aligned}&\cS(f_{1},\ldots, f_{k},W_1,\ldots ,W_m)=\\
&\mathcal{L}(f_{1} \otimes\mathbbm{1}_{\Omega_1} \otimes \mathbbm{1}_{\Omega_2},\mathbbm{1}_{\Omega_1}\otimes f_{1} \otimes\mathbbm{1}_{\Omega_2},\ldots, f_{k} \otimes\mathbbm{1}_{\Omega_1} \otimes \mathbbm{1}_{\Omega_2},\mathbbm{1}_{\Omega_1}\otimes f_{k} \otimes\mathbbm{1}_{\Omega_2},W_1,\ldots ,W_m),
\end{aligned}$$
which are the $P-$variables vectors variant of \ref{DefSf1f2fk}.

Using these quantities we can define $k-$profiles $\cS_k(W_1,\ldots,W_m)$ and therefore a metric and a convergence notion for vectors of $P-$variables. This type of convergence is useful for comparing limits of graphs and weighted graphs.

We explain this with a motivating example.
\begin{example}
Recall the Erd\"os-Renyi random graph $G(n,1/2)$ (Example Erd\"os-Renyi random graph). Let for every positive integer $n$ the graph $G_n$ be a realization of $G(n,1/2).$ The sequence $(G_n)_n$ considered as a sequence of $P-$variables (recall Example \ref{ExGraphasGraphset}) converges to the $P-$variable $\mathbbm{1}_{[0,1]\times [0,1]\times [0,1/2] }.$ Let consider for every $G_n$ its complement graph $G_n^{\prime}$, i.e.\ the graph $G_n^{\prime}$ on the same set of vertices as of $G_n$ such that there is an edge between two vertices $\{v, w\}$ in $G_n^{\prime}$, if and only if there is no edge between $\{v, w\}$ in $G_n.$  Also the sequence $(G^{\prime}_n)_n$ considered as a sequence of $P-$variables converges to the $P-$variable $\mathbbm{1}_{[0,1]\times [0,1]\times [0,1/2] }.$ However, if we consider the sequence of vectors of $P-$variables $((G_n,G_n^{\prime}))_n$ it is easy to observe that it converges to the vector of $P-$variables \\$(\mathbbm{1}_{[0,1]\times [0,1]\times [0,1/2] },\mathbbm{1}_{[0,1]\times [0,1]\times [1/2,1] }),$ which is not identified with $(\mathbbm{1}_{[0,1]\times [0,1]\times [0,1/2] },\mathbbm{1}_{[0,1]\times [0,1]\times [0,1/2] })$ in the $P-$variables metric for $P-$variable vectors. Observe that $\mathbbm{1}_{[0,1]\times [0,1]\times [0,1/2] }+\mathbbm{1}_{[0,1]\times [0,1]\times [1/2,1] }=\mathbbm{1}_{[0,1]\times [0,1]\times [0,1] }$ that is what one wants because $G_n+G^{\prime}_n=\mathbbm{1}_{[n]\times [n]\times [0,1]}.$
\end{example}

The previous example shows that, even if multiple sequences of graphs $(G^1_n)_n,\ldots,(G^m_n)_n,$ where for every positive integer $n$ the graphs $G^1_n,\ldots,G^m_n$ are on the same vertex set, have the same limit in $P-$variables convergence we can still distinguish them in the limit considering the convergence of the sequence of $P-$variable vectors $((G^1_n,\ldots,G^m_n))_n.$  

Observe that, similarly to random vectors and their distributions, a vector $P-$variable can be associated with a probability graphon taking values in $\cP(\Space )$ where $\Space \subset \R^m$ for $m\geq 1,$ see Remark \ref{RemarkgeneralCasePolish}.

\subsection{Probability bigraphons and P-bivariables:}

Recall the notions of probability bigraphons (Remark \ref{RmkProbBigraphons}) and $P-$bivariable (Definition \ref{DefBBBiset}). We already saw that one can define for probability bigraphons an unlabelled cut distance (recall Remark \ref{RKBigraphonsDistance}) and for $P-$bivariables, one can define a $P-$bivariables metric similar to the $P-$variables metric, see Remark \ref{RmkDistanseBBBisets}. We did not explore the connection between probability graphons and $P-$bivariables and their respective metrics here in detail but we expect the two formulations to be equivalent also in this case, recall Theorem \ref{TH:equivalence_dM_cut_dist} for $P-$variables and probability graphons. In particular, we expect all the proof techniques and results in this work to carry out to probability $P-$bivariables and probability bigraphons.

\subsection{Changing probability space}

In this work, we considered sequences of graphs considered as $P-$variables as explained in Examples \ref{ExGraphasPsets} and \ref{ExGraphasGraphset}. However, this is not the only possibility. In particular, when considering a graph as a $P-$variable we considered its adjacency matrix and the uniform measure of the vertex set. However, one can consider a different matrix related to the graph (the Laplacian matrix for example) and a different measure on the vertex set (the stationary measure of the random walk on the graph for example).

One can also consider the product of the vertex set with itself equipped with a measure that is not a product measure. This might be interesting for obtaining meaningful convergence notions for graphs with self-loops (when we consider the convergence notion considered in this work the diagonal has measure zero in the limit losing all the possible information about self-loops) and for sparse graphs.      

\

\subsection{On a path towards hypergraph limits}

An important advantage of $P-$variables and their metric over probability graphons and the associated unlabelled cut-metric is that it is more transparent how to generalize these notions to cover hypergraphs. To represent hypergraph limits for uniform hypergraphs with edge cardinality $3$, for example, one uses measurable functions $$
W:\Omega_1\times \Omega_1 \times \Omega_1 \times \Omega_2\times \Omega_2 \times \Omega_2 \times \Omega_3 \rightarrow \R
$$
or measurable functions
$$
w:\Omega_1\times \Omega_1 \times \Omega_1 \times \Omega_2\times \Omega_2 \times \Omega_2  \rightarrow \R
$$
where $\Omega_1,\Omega_2 $ and $\Omega_3$ are probability spaces.

We will explore this direction and the connections of this convergence with the convergence notions for hypergraphs \cite{hypergrELEK20121731,HypergraphonsZhao,zucal2023action} and Aldous-Hoover theorem \cite{aldous1981representations,hoover1979relations,austin2008exchangeable,diaconis2007graph} in future work.
\newline

\textbf{Acknowledgements:}  The author would like to thank Rostislav Matveev for a useful discussion that led to the proof of Lemma \ref{MeasProbGraphonFromPset}, Florentin M\"unch for useful conversations about the Levy-Prokhorov distance and Julien Weibel for helpful discussions, in particular on probability graphons, and for pointing out several references.
\medskip

\section*{Appendix}
In this appendix, we present results from probability theory and measure theory which we used throughout our work.

Let's consider a metric space $(X,d).$  Portmanteau theorem characterizes weak convergence of measures.

\begin{theorem}[Portmanteau theorem, Theorem 2.1 in \cite{billingsley1968convergence}]\label{PortThm}
    Let $(\mu_n)_n$ a sequence of measures $\mu_n$ on $X$ and $\mu$ a measure on $X$ such that $\lim_n\mu_n(X)=\mu(X)$. Then the following statements are equivalent.
    \begin{itemize}
        \item $\mu_n$ weakly converges to $\mu$, i.e. $\lim \mu_n[f]=\mu[f]$ for all $f\in \CbFunct$.
        \item $\liminf_n \mu_n[ f]\geq \mu[f] $ for all lower semi-continuous functions $f: X \rightarrow \mathbb{R}$ that are bounded from below.
        \item  $\limsup_n\mu_n[f] \leq  \mu[f]$ for all upper semi-continuous functions $f: X \rightarrow \mathbb{R}$ that are bounded from above.
        \item $\liminf _n \mu_n(O) \geq \mu(O)$ for all open subsets $O$ of $X$.
        \item  $\limsup _n \mu_n(C) \leq \mu(C)$ for all closed subsets $C$ of $X$.
        \item  $\lim _n \mu_n(A)=\mu(A)$ for all Borel subsets $A$ of $\Space$ such that $\mu(\partial A)=0$, where $\partial A$ denotes the topological boundary of $A$.
        \item $\lim _n \mu_n[f]=\mu[f]$ for all bounded Lipschitz-continuous functions $f: X \rightarrow \mathbb{R}$.
    \end{itemize}
\end{theorem}

The next lemma is a general probabilistic result about the weak convergence of random variables, products and expectations.

\begin{lemma}[Lemma 13.4 in \cite{backhausz2018action}]\label{closedlem2} Let $q\in (1,\infty)$.  Let $((X_n,Y_n))_n$ be a sequence of pairs of jointly distributed random variables (2-random vectors) such that $X_n$ is $[-1,1]-$valued and $\mathbb{E}(|Y_n|^q)=\|Y_n\|_q^q\leq c<\infty$ for some $c\in\mathbb{R}^+$. Assume that the distributions of $((X_n,Y_n))_n$ weakly converge to some probability distribution $(X,Y)$ as $n$ goes to infinity. Then $\mathbb{E}(|Y|^q)=\|Y\|_q^q\leq c$ and $$\lim_{i\to\infty} \mathbb{E}(X_iY_i)=\mathbb{E}(XY).$$
\end{lemma}

Let's now $X$ be the product space $X=S^{\prime}\times S^{\prime \prime}$ for $S^{\prime}$ and $ S^{\prime \prime}$ metric spaces and $\pi_1$ and $\pi_2$ the projection maps to $S^{\prime}$ and $ S^{\prime \prime}$  respectively. For a measure $\mu$ on $X$ we define with ${(\pi_1)}_{\#}\mu$ and ${(\pi_2)}_{\#}\mu$ the marginal measures, i.e.\ the pushforwards of $\mu$ through $\pi_1$ and $\pi_2$ respectively. We have the following result.

\begin{theorem}[Theorem 2.8 in \cite{billingsley1968convergence}] \label{BillingWeakConv}
Let $X=S^{\prime}\times S^{\prime \prime}$ to be separable and $(\mu_n)$ a sequence of measures $\mu_n$ on $X$ and $\mu$ a measure on $X.$ The sequence $(\mu_n)_n$ weakly converges to $\mu$ if and only if $\mu_n(Q^{\prime }\times Q^{\prime \prime})\rightarrow \mu(Q^{\prime }\times Q^{\prime \prime})$ for any measurable set $Q^{\prime}\subset S^{\prime}$ such that ${(\pi_1)}_{\#}\mu(\partial Q^{\prime})$ and any measurable set $Q^{\prime \prime}\subset S^{\prime \prime}$ such that ${(\pi_2)}_{\#}\mu(\partial Q^{\prime\prime}).$
\end{theorem}

We conclude this appendix with some useful properties of the Levy-Prokhorov distance $d_{\mathcal{LP}}$ (Definition \ref{LevyProk}). We start with a useful bound.

\begin{lemma}[Lemma 2.4 in \cite{zucal2024probabilitygraphonsrightconvergence}]\label{LemmaIneqScalingProkhorov}
  Let $(X,d)$ be a metric space. Let $\mu,\nu$ two measures on $X$ and $\alpha \geq1$. Then 
$$
    d_{\mathcal{LP}}(\mu,\nu)\leq d_{\mathcal{LP}}(\alpha\mu,\alpha\nu)\leq \alpha d_{\mathcal{LP}}(\mu,\nu).$$ Moreover, the previous inequalities are sharp, i.e.\ there exist measures $\mu_1,\nu_1$ such that the first inequality is an equality and measures $\mu_2,\nu_2$ such that the second inequality is an equality.
\end{lemma}

The Levy-Prokhorov distance is also a quasi-convex metric.

\begin{lemma}[Lemma 3.21 in \cite{abraham2023probabilitygraphons}]\label{LemmaQuasi-convProhorov}
    The Levy-Prokhorov distance $d_{\mathcal{LP}}$ is quasi-convex on $\mathcal{M}_{+}$, i.e.\ for any measures $\mu_1,\mu_2,\nu_1,\nu_2 \in \mathcal{M}_{+} $ and any $\alpha \in [0,1]$ we have $$
    d_{\mathcal{LP}}(\alpha \mu_1+(1-\alpha)\mu_2,\alpha \nu_1+(1-\alpha)\nu_2)\leq \max(d_{\mathcal{LP}}(\mu_1, \nu_1),d_{\mathcal{LP}}( \mu_2,\nu_2)).    
    $$
\end{lemma}

Let's consider again the product space $X=S^{\prime}\times S^{\prime \prime}$ for $S^{\prime}$ and $ S^{\prime \prime}$ metric spaces and let $\pi=\pi_1$ (or $\pi=\pi_2$) the projection map to $S^{\prime}$ (or to $ S^{\prime \prime}$). We have that the Levy-Prokhorov distance between two measures is an upper bound to the Levy-Prokhorov distance of the marginals.

\begin{lemma}\label{LemmIneqProkh}
Let $\mu$ and $\nu$ be two measures. The following inequality holds:
$$
d_{\mathcal{LP}}({(\pi)}_{\#}\mu,{(\pi)}_{\#}\nu)\leq d_{\mathcal{LP}}(\mu,\nu)
$$
\end{lemma}
\proof
The result follows from the inequality 
$$\begin{aligned}
d_{\mathcal{LP}}({(\pi)}_{\#}\mu,{(\pi)}_{\#}\nu)=&\inf \left\{\varepsilon>0:{(\pi)}_{\#}\mu(U) \leq {(\pi)}_{\#}\nu\left(U^{\varepsilon}\right)+\varepsilon \text{ and } \right.\\
&\left.{(\pi)}_{\#}\mu(U) \leq {(\pi)}_{\#}\nu\left(U^{\varepsilon}\right)+\varepsilon  \text{ for all } U \in \mathcal{B}(S^{\prime})\right\}\\
=&\inf \left\{\varepsilon>0:\mu(S^{\prime} \times U) \leq \nu\left((S^{\prime}\times U)^{\varepsilon}\right)+\varepsilon \text{ and } \right.\\
&\left.\mu(S^{\prime} \times U) \leq \nu\left((Y \times U)^{\varepsilon}\right)+\varepsilon  \text{ for all } U \in \mathcal{B}(S^{\prime})\right\}\\
=&\inf \left\{\varepsilon>0:\mu(S^{\prime} \times U) \leq \nu\left(S^{\prime} \times U^{\varepsilon}\right)+\varepsilon \text{ and } \right.\\
&\left.\mu(S^{\prime}\times U) \leq \nu\left(S^{\prime}\times U^{\varepsilon}\right)+\varepsilon  \text{ for all } U \in \mathcal{B}(S^{\prime})\right\}\\
\leq &\inf \left\{\varepsilon>0:\mu(V) \leq \nu\left(V^{\varepsilon}\right)+\varepsilon \text{ and } \right.\\
&\left.\mu(V) \leq \nu\left(V^{\varepsilon}\right)+\varepsilon  \text{ for all } V \in \mathcal{B}(X))\right\}\\
=& d_{\mathcal{LP}}\left(\mu, \nu\right),
\end{aligned}$$
where the inequality holds because the infimum on a smaller set is smaller and we denoted with $\cB(S^{\prime})$ and $\cB(X)$ the Borel $\sigma-$algebras of $S^{\prime}$ and $X$ respectively.
\endproof

In the following,  for a measure $\mu$ and a measurable set $E\subset X,$  we will denote with $\mu_E$ the restriction of the measure $\mu$ with respect to $E,$ i.e. the measure $\mu_E$ such that for any measurable set $Q$
\begin{equation}
    \mu_E(Q)=\mu(E\cap Q).
\end{equation}
Moreover, we recall that for a measurable set $S$ and $\alpha>0$ we denote with $S^{\alpha}$ the set of points at a distance smaller than $\alpha $ from $S$.

\begin{lemma}\label{LemmCntinuityProhRestriction}
Let $E$ be a Borel subset of $X$ and $S$ be a measurable subset of $ E$ such that there exists an $\alpha>0$ such that $S^{\alpha}\subset E \subset X.$  Let's consider $\cK$ the set of measures $\eta$ on $X$ such that $\eta(S)=\eta(E).$ For any two measures $\mu,\nu \in \cK$ 
such that $
d_{\mathcal{LP}}(\mu,\nu)<\alpha
$ we have 
$$
d_{\mathcal{LP}}(\mu_E,\nu_E)\leq d_{\mathcal{LP}}(\mu,\nu).
$$
\end{lemma}
\proof
Let $\alpha>\delta >d_{\mathcal{LP}}(\mu,\nu) ,$ then for all measurable subsets $A$ of $X$ we have $\mu(A)\leq \nu(A^{\delta})+\delta.$ In particular, choosing $A=B\cap S$  we obtain $\mu(B\cap S)\leq \nu((B\cap S)^{\delta})+\delta.$ Therefore, we obtain
$$\begin{aligned}
&\mu_E(B)=\\
   & \mu(B\cap E)=\\
   &\mu(B\cap E \cap S)+\mu(B\cap E \cap S^c)\leq \\
   &\mu(B\cap S)+\mu( E \cap S^c)\leq \\
   &\mu(B\cap S)\leq \\
   & \nu((B\cap S)^{\delta})+\delta\leq \\
  & \nu(B^{\delta}\cap S^{\delta})+\delta\leq \\
  & \nu(B^{\delta}\cap E)+\delta= \\
  & \nu_E(B^{\delta})+\delta.
\end{aligned}
$$
Sending $\delta \rightarrow d_{\mathcal{LP}}(\mu,\nu),$ as we can also exchange the role of $\mu $ and $\nu,$ we obtain the inequality.
\endproof

For a measure $\mu$ on $\R^k$ we denote with $\tau(\mu)$ the quantity
\begin{equation}\label{eqn:tau}\tau(\mu)=\max_{1\leq i\leq k} \int_{(x_1,x_2,\dots,x_k)\in\mathbb{R}^k}|x_i|~d\mu.\end{equation}
Moreover, for a random vector $Z=(Z_1,\ldots, Z_k)$ we have its distribution $\tau(\mathcal{L}(Z))=\tau(\mathcal{L}((Z_1,\ldots, Z_k)))$ and we denote with $\tau(Z)$ the quantity $\tau(\mathcal{L}(Z))=\max_{i\in [k]}\|Z_i\|_1.$ We have the following bound.

\begin{lemma}[Lemma 13.1 in \cite{backhausz2018action}]\label{coupdist} Let $X,Y$ be two jointly distributed 
$k-$random vectors. Then $$d_{\mathcal{LP}}(\mathcal{L}(X),\mathcal{L}(Y))\leq \tau(X-Y)^{1/2}k^{3/4},$$

where $\tau$ is defined as in \eqref{eqn:tau}.
\end{lemma}

As a direct consequence, one has the following lemma.

\begin{lemma}[Lemma 13.2 in \cite{backhausz2018action}] \label{coupdist2} 
Let $v_1,v_2,\dots,v_k$ and $w_1,w_2,\dots,w_k$ be in $L^1(\Omega)$ for some probability space $\Omega$. Let $m:=\max_{i\in [k]} \|v_i-w_i\|_1$. Then
$$d_{\mathcal{LP}}(\mathcal{L}(v_1,v_2,\dots,v_k),\mathcal{L}(w_1,w_2,\dots,w_k))\leq m^{1/2}k^{3/4}.$$
\end{lemma}

We prove now a technical lemma about the Levy-Prokhorov metric that is useful throughout this work.

Let's consider $\R^k$ with the Euclidean distance and its subset $\bigcup^k_{i=1}\{e_i\}\subset \R^k.$
We denote with $M_k$ the set of measures $\nu $ such that the support of $\nu$ is contained in $\bigcup^k_{i=1}\{e_i\}, $ i.e.\

$$
M_k=\left\{\nu\in \cP(\R^k):\ \text{supp}(\nu)\subset \bigcup^k_{i=1}\{e_i\}\right\}\subset \cP(\R^k).
$$
\begin{lemma}\label{LemmApproxFunctPartMeas}
Let $f_1,\ldots,f_k\in L_{[-1,1]}^{\infty}(\Omega)$ such that the distribution $\mu=\mathcal{L}(f_1,\ldots, f_k)\in \cP(\R^k)$ of the random vector $(f_1,\ldots,f_k)$ is at Levy-Prokhorov distance $d_{\mathcal{LP}}$ from $M_k\subset \cP(\R^k)$ at most $\delta>0,$ i.e.\ $d_{\mathcal{LP}}(\mu, M_k)<\delta.$  Then there exists a measurable partition $\cP=\{P_1,\ldots, P_k\}$ such that for $g_i=\mathbbm{1}_{P_i}$ we have $\|f_i-g_i\|_p<C_k\delta$ and such that the distribution $\nu=\mathcal{L}(g_1,\ldots,g_k)$ is such that $d_{\mathcal{LP}}(\mu, \nu)\leq \delta.$   
\end{lemma}
\proof
 Let's define for $\delta>0$ the set $S_{\delta}=[-1,\delta)\cup (\delta,1-\delta).$  Moreover, let's define $\mu_i=\pi_{i\#}\mu$ and $\nu_i=\pi_{i\#}\nu$ where $\pi_i$ denotes the projection on the $i-$th coordinate and $\tilde{g_i}=\mathbbm{1}_{f_i^{-1}(S_{\delta})}.$ 

We have $$\|f_i-\tilde{g}_i\|_p=\left(\int_{\Omega}\left|f_i(\omega)-\mathbbm{1}_{f_i^{-1}(S_{\delta})}(\omega) \right|^p\P(\drv \omega)\right)^{1/p}\leq  \delta+ \P(f_i^{-1}(S_{\delta}))=\delta+\mu_i(S_{\delta})$$
Moreover, from $d_{\mathcal{LP}}(\mu, M_k)\leq \delta,$ there exists a measure $\nu \in M_k$ such that $d_{\mathcal{LP}}(\mu, \nu)\leq \delta.$ From Lemma \ref{LemmIneqProkh} we have that $d_{\mathcal{LP}}(\mu_i, \nu_i)\leq d_{\mathcal{LP}}(\mu, \nu)\leq \delta.$

Therefore, from the definition of Levy-Prokhorov distance, it follows that $\mu_i(S_{\delta})\leq \nu_i((-1,1)\setminus \{0\})+\delta=\delta$ because $\nu_i((-1,1)\setminus \{0\})=0$ as $\nu\in M_k.$

Thus we obtain 
$$\|f_i-\tilde{g}_i\|_p=\delta+\mu_i(S_{\delta})\leq 2 \delta.$$

Moreover, from this last inequality and the triangular inequality we obtain $$\|\sum^k_{i=1}f_i-\sum^k_{i=1}\tilde{g}_i\|_p=\|\sum^k_{i=1}(f_i-\tilde{g}_i)\|_p\leq\sum^k_{i=1}\|f_i-\tilde{g}_i\|_p\leq 2k \delta.$$
We have also the chain of inequalities: $$\begin{aligned}
   & \|\sum^k_{i=1}f_i-\mathbbm{1}_{\Omega}\|_p\leq (k+1) \P\left(\ \left\{\omega: \ |\sum^k_{i=1}f_i(\omega)-\mathbbm{1}_{\Omega}(\omega)|>k^2\delta\right\}\right)+k^2\delta\leq\\
   &(k+1)\P(\{\omega:\ (f_1 ,\ldots ,f_k)(\omega)\notin \bigcup^k_{i=1} R(e_i,k\delta) \})+k^2\delta=\\
   &(k+1)\mu\left(\R^k\setminus\bigcup^k_{i=1} R(e_i,k\delta) \right)+k^2\delta\leq (k+1)\mu\left(\R^k\setminus\bigcup^k_{i=1} B(e_i,\delta)\right)+k^2\delta ,
   \end{aligned}$$
where with $R(p,r)$ and $B(p,r)$ we denote respectively the closed $L^1-$ball and the closed $L^2-$ball (Euclidean ball) in $\R^k$ of radius $r>0$ centred at point $p\in \R^k$ and we used monotonicity of probability measures in the second and third inequality and $B(p,\delta)\subset R(p,k\delta)$.

But again from $d_{\mathcal{LP}}(\mu, M_k)<\delta,$ there exists a measure $\nu \in M_k$ such that $d_{\mathcal{LP}}(\mu, \nu)< \delta.$ Therefore, from the definition of Levy-Prokhorov distance, it follows that $$
\mu\left(\R^k\setminus\bigcup^k_{i=1} B(e_i,\delta)\right)\leq\nu\left(\left(\R^k\setminus\bigcup^k_{i=1} B(e_i,\delta)\right)^{\delta}\right)+\delta=\nu\left(\R^k\setminus\bigcup^k_{i=1} \{e_i\}\right)+\delta=\delta,
$$
because $\nu\left(\R^k\setminus\bigcup^k_{i=1} \{e_i\}\right)=0$ as $\nu\in M_k.$ Therefore, we have $
\|\sum^k_{i=1}f_i-\mathbbm{1}_{\Omega}\|_p\leq (k^2+k+1)\delta
$

Putting everything together and using the triangular inequality we obtain 
$$\begin{aligned}   
   &\|\sum^k_{i=1}\tilde{g}_i-\mathbbm{1}_{\Omega}\|_p\leq \|\sum^k_{i=1}f_i-\mathbbm{1}_{\Omega}\|_p+\|\sum^k_{i=1}f_i-\sum^k_{i=1}\tilde{g}_i\|_p\leq \\
   &(k^2+k+1)\delta +2k\delta=(k^2+3k+1)\delta
   \end{aligned}$$
From the previous inequality we obtain 

$$\begin{aligned}
   & \P(\{\omega:\ \sum^k_{i=1}\mathbbm{1}_{f_i^{-1}((1-\delta,1])}(\omega)-\mathbbm{1}_{\Omega}(\omega)\neq 0 \})=\\
   &\P(\{\omega:\ \sum^k_{i=1}\tilde{g}_i(\omega)-\mathbbm{1}_{\Omega}(\omega)\neq 0 \})\leq \|\sum^k_{i=1}\tilde{g}_i-\mathbbm{1}_{\Omega}\|_p\leq (k^2+k+1)\delta.
\end{aligned}
$$
Therefore, we can define $$
g_1(\omega)=\begin{cases}
    1 &  \text{ if } \omega\in\{\omega:\ \sum^k_{i=1}\tilde{g}_i(\omega)-\mathbbm{1}_{\Omega}(\omega)\neq 0 \}\\
    \tilde{g}_1(\omega) & \text{else}
\end{cases}
$$

and for $i\in [k]\setminus \{1\}$
$$
g_i(\omega)=\begin{cases}
    0&  \text{ if } \omega\in\{\omega:\ \sum^k_{i=1}\tilde{g}_i(\omega)-\mathbbm{1}_{\Omega}(\omega)\neq 0 \}\\
    \tilde{g}_i(\omega) & \text{else}.
\end{cases}
$$
We observe that the functions $g_1,\ldots,g_k$ are $\{0,1\}-$valued and $\sum^k_{i=1}g_i=\mathbbm{1}_{\Omega}$ and thus there exists a partition $\cP=\{P_1,\ldots, P_k\}$ such that $g_i=\mathbbm{1}_{P_i}$ for every $i\in [k]. $ Moreover, we have $\|g_i-\tilde{g}_i\|_p\leq (k+1)\P(\{\omega:\ \sum^k_{i=1}\tilde{g}_i(\omega)-\mathbbm{1}_{\Omega}(\omega)\neq 0 \})\leq (k+1)(k^2+3k+1)\delta.$ 

Thus by the triangular inequality, we obtain
$$\begin{aligned}
   & \|f_i-g_i\|_p\leq \|f_i-\tilde{g}_i\|_p+\|g_i-\tilde{g}_i\|_p\leq 2\delta+(k+1)(k^2+3k+1)\delta=(k^3+3k^2+5k+1)\delta.
   \end{aligned}$$
This concludes the proof.
\endproof
\begin{remark}
In the language used in this work, the set of functions $\{g_1,\ldots,g_k\}$ in Lemma \ref{LemmApproxFunctPartMeas} is called a function partition.
\end{remark}
\section*{References}

\bibliographystyle{plain} 
\bibliography{biblio}
\end{document}